\documentstyle[amssymb,11pt]{amsart}
\title{Stacks of Stable Maps and Gromov-Witten Invariants}
\dedicatory{To Professor Friedrich Hirzebruch}
\date{October 18, 1995}
\author{K. Behrend}
\address{University of British Columbia\\
Mathematics Department\\
121--1984 Mathematics Road\\
Vancouver, B.C.\ Canada V6T 1Z2}
\email{behrend@@math.ubc.ca}
\author{Yu. Manin}
\address{Max-Planck-Institut f\"ur Mathematik\\
Gottfried-Claren-Str.\ 26\\
D--53225 Bonn}
\email{manin@@mpim-bonn.mpg.de}

\catcode`@=11 \catcode`!=11

\expandafter\ifx\csname fiverm\endcsname\relax
  \let\fiverm\fivrm
\fi

\let\!latexendpicture=\endpicture
\let\!latexframe=\frame
\let\!latexlinethickness=\linethickness
\let\!latexmultiput=\multiput
\let\!latexput=\put

\def\@picture(#1,#2)(#3,#4){%
  \@picht #2\unitlength
  \setbox\@picbox\hbox to #1\unitlength\bgroup
  \let\endpicture=\!latexendpicture
  \let\frame=\!latexframe
  \let\linethickness=\!latexlinethickness
  \let\multiput=\!latexmultiput
  \let\put=\!latexput
  \hskip -#3\unitlength \lower #4\unitlength \hbox\bgroup}

\catcode`@=12 \catcode`!=12

\newcommand{\fiverm}{}

\catcode`!=11 

\def\PiC{P\kern-.12em\lower.5ex\hbox{I}\kern-.075emC}
\def\PiCTeX{\PiC\kern-.11em\TeX}

\def\!ifnextchar#1#2#3{%
  \let\!testchar=#1%
  \def\!first{#2}%
  \def\!second{#3}%
  \futurelet\!nextchar\!testnext}
\def\!testnext{%
  \ifx \!nextchar \!spacetoken
    \let\!next=\!skipspacetestagain
  \else
    \ifx \!nextchar \!testchar
      \let\!next=\!first
    \else
      \let\!next=\!second
    \fi
  \fi
  \!next}
\def\\{\!skipspacetestagain}
  \expandafter\def\\ {\futurelet\!nextchar\!testnext}
\def\\{\let\!spacetoken= } \\  

\def\!tfor#1:=#2\do#3{%
  \edef\!fortemp{#2}%
  \ifx\!fortemp\!empty
    \else
    \!tforloop#2\!nil\!nil\!!#1{#3}%
  \fi}
\def\!tforloop#1#2\!!#3#4{%
  \def#3{#1}%
  \ifx #3\!nnil
    \let\!nextwhile=\!fornoop
  \else
    #4\relax
    \let\!nextwhile=\!tforloop
  \fi
  \!nextwhile#2\!!#3{#4}}

\def\!etfor#1:=#2\do#3{%
  \def\!!tfor{\!tfor#1:=}%
  \edef\!!!tfor{#2}%
  \expandafter\!!tfor\!!!tfor\do{#3}}

\def\!cfor#1:=#2\do#3{%
  \edef\!fortemp{#2}%
  \ifx\!fortemp\!empty
  \else
    \!cforloop#2,\!nil,\!nil\!!#1{#3}%
  \fi}
\def\!cforloop#1,#2\!!#3#4{%
  \def#3{#1}%
  \ifx #3\!nnil
    \let\!nextwhile=\!fornoop
  \else
    #4\relax
    \let\!nextwhile=\!cforloop
  \fi
  \!nextwhile#2\!!#3{#4}}

\def\!ecfor#1:=#2\do#3{%
  \def\!!cfor{\!cfor#1:=}%
  \edef\!!!cfor{#2}%
  \expandafter\!!cfor\!!!cfor\do{#3}}

\def\!empty{}
\def\!nnil{\!nil}
\def\!fornoop#1\!!#2#3{}

\def\!ifempty#1#2#3{%
  \edef\!emptyarg{#1}%
  \ifx\!emptyarg\!empty
    #2%
  \else
    #3%
  \fi}

\def\!getnext#1\from#2{%
  \expandafter\!gnext#2\!#1#2}%
\def\!gnext\\#1#2\!#3#4{%
  \def#3{#1}%
  \def#4{#2\\{#1}}%
  \ignorespaces}

\def\!getnextvalueof#1\from#2{%
  \expandafter\!gnextv#2\!#1#2}%
\def\!gnextv\\#1#2\!#3#4{%
  #3=#1%
  \def#4{#2\\{#1}}%
  \ignorespaces}

\def\!copylist#1\to#2{%
  \expandafter\!!copylist#1\!#2}
\def\!!copylist#1\!#2{%
  \def#2{#1}\ignorespaces}

\def\!wlet#1=#2{%
  \let#1=#2
  \wlog{\string#1=\string#2}}

\def\!listaddon#1#2{%
  \expandafter\!!listaddon#2\!{#1}#2}
\def\!!listaddon#1\!#2#3{%
  \def#3{#1\\#2}}

\def\!rightappend#1\withCS#2\to#3{\expandafter\!!rightappend#3\!#2{#1}#3}
\def\!!rightappend#1\!#2#3#4{\def#4{#1#2{#3}}}

\def\!leftappend#1\withCS#2\to#3{\expandafter\!!leftappend#3\!#2{#1}#3}
\def\!!leftappend#1\!#2#3#4{\def#4{#2{#3}#1}}

\def\!lop#1\to#2{\expandafter\!!lop#1\!#1#2}
\def\!!lop\\#1#2\!#3#4{\def#4{#1}\def#3{#2}}

\def\!loop#1\repeat{\def\!body{#1}\!iterate}
\def\!iterate{\!body\let\!next=\!iterate\else\let\!next=\relax\fi\!next}

\def\!!loop#1\repeat{\def\!!body{#1}\!!iterate}
\def\!!iterate{\!!body\let\!!next=\!!iterate\else\let\!!next=\relax\fi\!!next}

\def\!removept#1#2{\edef#2{\expandafter\!!removePT\the#1}}
{\catcode`p=12 \catcode`t=12 \gdef\!!removePT#1pt{#1}}

\def\placevalueinpts of <#1> in #2 {%
  \!removept{#1}{#2}}

\def\!mlap#1{\hbox to 0pt{\hss#1\hss}}
\def\!vmlap#1{\vbox to 0pt{\vss#1\vss}}

\def\!not#1{%
  #1\relax
    \!switchfalse
  \else
    \!switchtrue
  \fi
  \if!switch
  \ignorespaces}

\let\!!!wlog=\wlog              
\def\wlog#1{}

\newdimen\headingtoplotskip     
\newdimen\linethickness         
\newdimen\longticklength        
\newdimen\plotsymbolspacing     
\newdimen\shortticklength       
\newdimen\stackleading          
\newdimen\tickstovaluesleading  
\newdimen\totalarclength        
\newdimen\valuestolabelleading  

\newbox\!boxA                   
\newbox\!boxB                   
\newbox\!picbox                 
\newbox\!plotsymbol             
\newbox\!putobject              
\newbox\!shadesymbol            

\newcount\!countA               
\newcount\!countB               
\newcount\!countC               
\newcount\!countD               
\newcount\!countE               
\newcount\!countF               
\newcount\!countG               
\newcount\!fiftypt              
\newcount\!intervalno           
\newcount\!npoints              
\newcount\!nsegments            
\newcount\!ntemp                
\newcount\!parity               
\newcount\!scalefactor          
\newcount\!tfs                  
\newcount\!tickcase             

\newdimen\!Xleft                
\newdimen\!Xright               
\newdimen\!Xsave                
\newdimen\!Ybot                 
\newdimen\!Ysave                
\newdimen\!Ytop                 
\newdimen\!angle                
\newdimen\!arclength            
\newdimen\!areabloc             
\newdimen\!arealloc             
\newdimen\!arearloc             
\newdimen\!areatloc             
\newdimen\!bshrinkage           
\newdimen\!checkbot             
\newdimen\!checkleft            
\newdimen\!checkright           
\newdimen\!checktop             
\newdimen\!dimenA               
\newdimen\!dimenB               
\newdimen\!dimenC               
\newdimen\!dimenD               
\newdimen\!dimenE               
\newdimen\!dimenF               
\newdimen\!dimenG               
\newdimen\!dimenH               
\newdimen\!dimenI               
\newdimen\!distacross           
\newdimen\!downlength           
\newdimen\!dp                   
\newdimen\!dshade               
\newdimen\!dxpos                
\newdimen\!dxprime              
\newdimen\!dypos                
\newdimen\!dyprime              
\newdimen\!ht                   
\newdimen\!leaderlength         
\newdimen\!lshrinkage           
\newdimen\!midarclength         
\newdimen\!offset               
\newdimen\!plotheadingoffset    
\newdimen\!plotsymbolxshift     
\newdimen\!plotsymbolyshift     
\newdimen\!plotxorigin          
\newdimen\!plotyorigin          
\newdimen\!rootten              
\newdimen\!rshrinkage           
\newdimen\!shadesymbolxshift    
\newdimen\!shadesymbolyshift    
\newdimen\!tenAa                
\newdimen\!tenAc                
\newdimen\!tenAe                
\newdimen\!tshrinkage           
\newdimen\!uplength             
\newdimen\!wd                   
\newdimen\!wmax                 
\newdimen\!wmin                 
\newdimen\!xB                   
\newdimen\!xC                   
\newdimen\!xE                   
\newdimen\!xM                   
\newdimen\!xS                   
\newdimen\!xaxislength          
\newdimen\!xdiff                
\newdimen\!xleft                
\newdimen\!xloc                 
\newdimen\!xorigin              
\newdimen\!xpivot               
\newdimen\!xpos                 
\newdimen\!xprime               
\newdimen\!xright               
\newdimen\!xshade               
\newdimen\!xshift               
\newdimen\!xtemp                
\newdimen\!xunit                
\newdimen\!xxE                  
\newdimen\!xxM                  
\newdimen\!xxS                  
\newdimen\!xxloc                
\newdimen\!yB                   
\newdimen\!yC                   
\newdimen\!yE                   
\newdimen\!yM                   
\newdimen\!yS                   
\newdimen\!yaxislength          
\newdimen\!ybot                 
\newdimen\!ydiff                
\newdimen\!yloc                 
\newdimen\!yorigin              
\newdimen\!ypivot               
\newdimen\!ypos                 
\newdimen\!yprime               
\newdimen\!yshade               
\newdimen\!yshift               
\newdimen\!ytemp                
\newdimen\!ytop                 
\newdimen\!yunit                
\newdimen\!yyE                  
\newdimen\!yyM                  
\newdimen\!yyS                  
\newdimen\!yyloc                
\newdimen\!zpt                  

\newif\if!axisvisible           
\newif\if!gridlinestoo          
\newif\if!keepPO                
\newif\if!placeaxislabel        
\newif\if!switch                
\newif\if!xswitch               

\newtoks\!axisLaBeL             
\newtoks\!keywordtoks           

\newwrite\!replotfile           

\newhelp\!keywordhelp{The keyword mentioned in the error message in unknown.
Replace NEW KEYWORD in the indicated response by the keyword that
should have been specified.}    

\!wlet\!!origin=\!xM                   
\!wlet\!!unit=\!uplength               
\!wlet\!Lresiduallength=\!dimenG       
\!wlet\!Rresiduallength=\!dimenF       
\!wlet\!axisLength=\!distacross        
\!wlet\!axisend=\!ydiff                
\!wlet\!axisstart=\!xdiff              
\!wlet\!axisxlevel=\!arclength         
\!wlet\!axisylevel=\!downlength        
\!wlet\!beta=\!dimenE                  
\!wlet\!gamma=\!dimenF                 
\!wlet\!shadexorigin=\!plotxorigin     
\!wlet\!shadeyorigin=\!plotyorigin     
\!wlet\!ticklength=\!xS                
\!wlet\!ticklocation=\!xE              
\!wlet\!ticklocationincr=\!yE          
\!wlet\!tickwidth=\!yS                 
\!wlet\!totalleaderlength=\!dimenE     
\!wlet\!xone=\!xprime                  
\!wlet\!xtwo=\!dxprime                 
\!wlet\!ySsave=\!yM                    
\!wlet\!ybB=\!yB                       
\!wlet\!ybC=\!yC                       
\!wlet\!ybE=\!yE                       
\!wlet\!ybM=\!yM                       
\!wlet\!ybS=\!yS                       
\!wlet\!ybpos=\!yyloc                  
\!wlet\!yone=\!yprime                  
\!wlet\!ytB=\!xB                       
\!wlet\!ytC=\!xC                       
\!wlet\!ytE=\!downlength               
\!wlet\!ytM=\!arclength                
\!wlet\!ytS=\!distacross               
\!wlet\!ytpos=\!xxloc                  
\!wlet\!ytwo=\!dyprime                 

\!zpt=0pt                              
\!xunit=1pt
\!yunit=1pt
\!arearloc=\!xunit
\!areatloc=\!yunit
\!dshade=5pt
\!leaderlength=24in
\!tfs=256                              
\!wmax=5.3pt                           
\!wmin=2.7pt                           
\!xaxislength=\!xunit
\!xpivot=\!zpt
\!yaxislength=\!yunit
\!ypivot=\!zpt
  \!dimenA=50pt \!fiftypt=\!dimenA     

\!rootten=3.162278pt                   
\!tenAa=8.690286pt                     
\!tenAc=2.773839pt                     
\!tenAe=2.543275pt                     

\def\!cosrotationangle{1}      
\def\!sinrotationangle{0}      
\def\!xpivotcoord{0}           
\def\!xref{0}                  
\def\!xshadesave{0}            
\def\!ypivotcoord{0}           
\def\!yref{0}                  
\def\!yshadesave{0}            
\def\!zero{0}                  

\let\wlog=\!!!wlog
\def\normalgraphs{%
  \longticklength=.4\baselineskip
  \shortticklength=.25\baselineskip
  \tickstovaluesleading=.25\baselineskip
  \valuestolabelleading=.8\baselineskip
  \linethickness=.4pt
  \stackleading=.17\baselineskip
  \headingtoplotskip=1.5\baselineskip
  \visibleaxes
  \ticksout
  \nogridlines
  \unloggedticks}
\def\setplotarea x from #1 to #2, y from #3 to #4 {%
  \!arealloc=\!M{#1}\!xunit \advance \!arealloc -\!xorigin
  \!areabloc=\!M{#3}\!yunit \advance \!areabloc -\!yorigin
  \!arearloc=\!M{#2}\!xunit \advance \!arearloc -\!xorigin
  \!areatloc=\!M{#4}\!yunit \advance \!areatloc -\!yorigin
  \!initinboundscheck
  \!xaxislength=\!arearloc  \advance\!xaxislength -\!arealloc
  \!yaxislength=\!areatloc  \advance\!yaxislength -\!areabloc
  \!plotheadingoffset=\!zpt
  \!dimenput {{\setbox0=\hbox{}\wd0=\!xaxislength\ht0=\!yaxislength\box0}}
     [bl] (\!arealloc,\!areabloc)}
\def\visibleaxes{%
  \def\!axisvisibility{\!axisvisibletrue}}

\def\!fixkeyword#1{%
  \errhelp=\!keywordhelp
  \errmessage{Unrecognized keyword `#1': \the\!keywordtoks{NEW KEYWORD}'}}

\!keywordtoks={enter `i\fixkeyword}

\def\fixkeyword#1{%
  \!nextkeyword#1 }

\def\axis {%
  \def\!nextkeyword##1 {%
    \expandafter\ifx\csname !axis##1\endcsname \relax
      \def\!next{\!fixkeyword{##1}}%
    \else
      \def\!next{\csname !axis##1\endcsname}%
    \fi
    \!next}%
  \!offset=\!zpt
  \!axisvisibility
  \!placeaxislabelfalse
  \!nextkeyword}

\def\!axisbottom{%
  \!axisylevel=\!areabloc
  \def\!tickxsign{0}%
  \def\!tickysign{-}%
  \def\!axissetup{\!axisxsetup}%
  \def\!axislabeltbrl{t}%
  \!nextkeyword}

\def\!axistop{%
  \!axisylevel=\!areatloc
  \def\!tickxsign{0}%
  \def\!tickysign{+}%
  \def\!axissetup{\!axisxsetup}%
  \def\!axislabeltbrl{b}%
  \!nextkeyword}

\def\!axisleft{%
  \!axisxlevel=\!arealloc
  \def\!tickxsign{-}%
  \def\!tickysign{0}%
  \def\!axissetup{\!axisysetup}%
  \def\!axislabeltbrl{r}%
  \!nextkeyword}

\def\!axisright{%
  \!axisxlevel=\!arearloc
  \def\!tickxsign{+}%
  \def\!tickysign{0}%
  \def\!axissetup{\!axisysetup}%
  \def\!axislabeltbrl{l}%
  \!nextkeyword}

\def\!axisshiftedto#1=#2 {%
  \if 0\!tickxsign
    \!axisylevel=\!M{#2}\!yunit
    \advance\!axisylevel -\!yorigin
  \else
    \!axisxlevel=\!M{#2}\!xunit
    \advance\!axisxlevel -\!xorigin
  \fi
  \!nextkeyword}

\def\!axisvisible{%
  \!axisvisibletrue
  \!nextkeyword}

\def\!axisinvisible{%
  \!axisvisiblefalse
  \!nextkeyword}

\def\!axislabel#1 {%
  \!axisLaBeL={#1}%
  \!placeaxislabeltrue
  \!nextkeyword}

\expandafter\def\csname !axis/\endcsname{%
  \!axissetup 
  \if!placeaxislabel
    \!placeaxislabel
  \fi
  \if +\!tickysign 
    \!dimenA=\!axisylevel
    \advance\!dimenA \!offset 
    \advance\!dimenA -\!areatloc 
    \ifdim \!dimenA>\!plotheadingoffset
      \!plotheadingoffset=\!dimenA 
    \fi
  \fi}

\def\grid #1 #2 {%
  \!countA=#1\advance\!countA 1
  \axis bottom invisible ticks length <\!zpt> andacross quantity {\!countA} /
  \!countA=#2\advance\!countA 1
  \axis left   invisible ticks length <\!zpt> andacross quantity {\!countA} / }

\def\plotheading#1 {%
  \advance\!plotheadingoffset \headingtoplotskip
  \!dimenput {#1} [B] <.5\!xaxislength,\!plotheadingoffset>
    (\!arealloc,\!areatloc)}

\def\!axisxsetup{%
  \!axisxlevel=\!arealloc
  \!axisstart=\!arealloc
  \!axisend=\!arearloc
  \!axisLength=\!xaxislength
  \!!origin=\!xorigin
  \!!unit=\!xunit
  \!xswitchtrue
  \if!axisvisible
    \!makeaxis
  \fi}

\def\!axisysetup{%
  \!axisylevel=\!areabloc
  \!axisstart=\!areabloc
  \!axisend=\!areatloc
  \!axisLength=\!yaxislength
  \!!origin=\!yorigin
  \!!unit=\!yunit
  \!xswitchfalse
  \if!axisvisible
    \!makeaxis
  \fi}

\def\!makeaxis{%
  \setbox\!boxA=\hbox{
    \beginpicture
      \!setdimenmode
      \setcoordinatesystem point at {\!zpt} {\!zpt}
      \putrule from {\!zpt} {\!zpt} to
        {\!tickysign\!tickysign\!axisLength}
        {\!tickxsign\!tickxsign\!axisLength}
    \endpicturesave <\!Xsave,\!Ysave>}%
    \wd\!boxA=\!zpt
    \!placetick\!axisstart}

\def\!placeaxislabel{%
  \advance\!offset \valuestolabelleading
  \if!xswitch
    \!dimenput {\the\!axisLaBeL} [\!axislabeltbrl]
      <.5\!axisLength,\!tickysign\!offset> (\!axisxlevel,\!axisylevel)
    \advance\!offset \!dp  
    \advance\!offset \!ht  
  \else
    \!dimenput {\the\!axisLaBeL} [\!axislabeltbrl]
      <\!tickxsign\!offset,.5\!axisLength> (\!axisxlevel,\!axisylevel)
  \fi
  \!axisLaBeL={}}

\def\arrow <#1> [#2,#3]{%
  \!ifnextchar<{\!arrow{#1}{#2}{#3}}{\!arrow{#1}{#2}{#3}<\!zpt,\!zpt> }}

\def\!arrow#1#2#3<#4,#5> from #6 #7 to #8 #9 {%
  \!xloc=\!M{#8}\!xunit
  \!yloc=\!M{#9}\!yunit
  \!dxpos=\!xloc  \!dimenA=\!M{#6}\!xunit  \advance \!dxpos -\!dimenA
  \!dypos=\!yloc  \!dimenA=\!M{#7}\!yunit  \advance \!dypos -\!dimenA
  \let\!MAH=\!M
  \!setdimenmode
  \!xshift=#4\relax  \!yshift=#5\relax
  \!reverserotateonly\!xshift\!yshift
  \advance\!xshift\!xloc  \advance\!yshift\!yloc
  \!xS=-\!dxpos  \advance\!xS\!xshift
  \!yS=-\!dypos  \advance\!yS\!yshift
  \!start (\!xS,\!yS)
  \!ljoin (\!xshift,\!yshift)
  \!Pythag\!dxpos\!dypos\!arclength
  \!divide\!dxpos\!arclength\!dxpos
  \!dxpos=32\!dxpos  \!removept\!dxpos\!!cos
  \!divide\!dypos\!arclength\!dypos
  \!dypos=32\!dypos  \!removept\!dypos\!!sin
  \!halfhead{#1}{#2}{#3}
  \!halfhead{#1}{-#2}{-#3}
  \let\!M=\!MAH
  \ignorespaces}
  \def\!halfhead#1#2#3{%
    \!dimenC=-#1%
    \divide \!dimenC 2 
    \!dimenD=#2\!dimenC
    \!rotate(\!dimenC,\!dimenD)by(\!!cos,\!!sin)to(\!xM,\!yM)
    \!dimenC=-#1
    \!dimenD=#3\!dimenC
    \!dimenD=.5\!dimenD
    \!rotate(\!dimenC,\!dimenD)by(\!!cos,\!!sin)to(\!xE,\!yE)
    \!start (\!xshift,\!yshift)
    \advance\!xM\!xshift  \advance\!yM\!yshift
    \advance\!xE\!xshift  \advance\!yE\!yshift
    \!qjoin (\!xM,\!yM) (\!xE,\!yE)
    \ignorespaces}

\def\betweenarrows #1#2 from #3 #4 to #5 #6 {%
  \!xloc=\!M{#3}\!xunit  \!xxloc=\!M{#5}\!xunit%
  \!yloc=\!M{#4}\!yunit  \!yyloc=\!M{#6}\!yunit%
  \!dxpos=\!xxloc  \advance\!dxpos by -\!xloc
  \!dypos=\!yyloc  \advance\!dypos by -\!yloc
  \advance\!xloc .5\!dxpos
  \advance\!yloc .5\!dypos
  \let\!MBA=\!M
  \!setdimenmode
  \ifdim\!dypos=\!zpt
    \ifdim\!dxpos<\!zpt \!dxpos=-\!dxpos \fi
    \put {\!lrarrows{\!dxpos}{#1}}#2{} at {\!xloc} {\!yloc}
  \else
    \ifdim\!dxpos=\!zpt
      \ifdim\!dypos<\!zpt \!dypos=-\!zpt \fi
      \put {\!udarrows{\!dypos}{#1}}#2{} at {\!xloc} {\!yloc}
    \fi
  \fi
  \let\!M=\!MBA
  \ignorespaces}

\def\!lrarrows#1#2{
  {\setbox\!boxA=\hbox{$\mkern-2mu\mathord-\mkern-2mu$}%
   \setbox\!boxB=\hbox{$\leftarrow$}\!dimenE=\ht\!boxB
   \setbox\!boxB=\hbox{}\ht\!boxB=2\!dimenE
   \hbox to #1{$\mathord\leftarrow\mkern-6mu
     \cleaders\copy\!boxA\hfil
     \mkern-6mu\mathord-$%
     \kern.4em $\vcenter{\box\!boxB}$$\vcenter{\hbox{#2}}$\kern.4em
     $\mathord-\mkern-6mu
     \cleaders\copy\!boxA\hfil
     \mkern-6mu\mathord\rightarrow$}}}

\def\!udarrows#1#2{
  {\setbox\!boxB=\hbox{#2}%
   \setbox\!boxA=\hbox to \wd\!boxB{\hss$\vert$\hss}%
   \!dimenE=\ht\!boxA \advance\!dimenE \dp\!boxA \divide\!dimenE 2
   \vbox to #1{\offinterlineskip
      \vskip .05556\!dimenE
      \hbox to \wd\!boxB{\hss$\mkern.4mu\uparrow$\hss}\vskip-\!dimenE
      \cleaders\copy\!boxA\vfil
      \vskip-\!dimenE\copy\!boxA
      \vskip\!dimenE\copy\!boxB\vskip.4em
      \copy\!boxA\vskip-\!dimenE
      \cleaders\copy\!boxA\vfil
      \vskip-\!dimenE \hbox to \wd\!boxB{\hss$\mkern.4mu\downarrow$\hss}
      \vskip .05556\!dimenE}}}

\def\putbar#1breadth <#2> from #3 #4 to #5 #6 {%
  \!xloc=\!M{#3}\!xunit  \!xxloc=\!M{#5}\!xunit%
  \!yloc=\!M{#4}\!yunit  \!yyloc=\!M{#6}\!yunit%
  \!dypos=\!yyloc  \advance\!dypos by -\!yloc
  \!dimenI=#2
  \ifdim \!dimenI=\!zpt 
    \putrule#1from {#3} {#4} to {#5} {#6} 
  \else 
    \let\!MBar=\!M
    \!setdimenmode 
    \divide\!dimenI 2
    \ifdim \!dypos=\!zpt
      \advance \!yloc -\!dimenI 
      \advance \!yyloc \!dimenI
    \else
      \advance \!xloc -\!dimenI 
      \advance \!xxloc \!dimenI
    \fi
    \putrectangle#1corners at {\!xloc} {\!yloc} and {\!xxloc} {\!yyloc}
    \let\!M=\!MBar 
  \fi
  \ignorespaces}

\def\setbars#1breadth <#2> baseline at #3 = #4 {%
  \edef\!barshift{#1}%
  \edef\!barbreadth{#2}%
  \edef\!barorientation{#3}%
  \edef\!barbaseline{#4}%
  \def\!bardobaselabel{\!bardoendlabel}%
  \def\!bardoendlabel{\!barfinish}%
  \let\!drawcurve=\!barcurve
  \!setbars}
\def\!setbars{%
  \futurelet\!nextchar\!!setbars}
\def\!!setbars{%
  \if b\!nextchar
    \def\!!!setbars{\!setbarsbget}%
  \else
    \if e\!nextchar
      \def\!!!setbars{\!setbarseget}%
    \else
      \def\!!!setbars{\relax}%
    \fi
  \fi
  \!!!setbars}
\def\!setbarsbget baselabels (#1) {%
  \def\!barbaselabelorientation{#1}%
  \def\!bardobaselabel{\!!bardobaselabel}%
  \!setbars}
\def\!setbarseget endlabels (#1) {%
  \edef\!barendlabelorientation{#1}%
  \def\!bardoendlabel{\!!bardoendlabel}%
  \!setbars}

\def\!barcurve #1 #2 {%
  \if y\!barorientation
    \def\!basexarg{#1}%
    \def\!baseyarg{\!barbaseline}%
  \else
    \def\!basexarg{\!barbaseline}%
    \def\!baseyarg{#2}%
  \fi
  \expandafter\putbar\!barshift breadth <\!barbreadth> from {\!basexarg}
    {\!baseyarg} to {#1} {#2}
  \def\!endxarg{#1}%
  \def\!endyarg{#2}%
  \!bardobaselabel}

\def\!!bardobaselabel "#1" {%
  \put {#1}\!barbaselabelorientation{} at {\!basexarg} {\!baseyarg}
  \!bardoendlabel}

\def\!!bardoendlabel "#1" {%
  \put {#1}\!barendlabelorientation{} at {\!endxarg} {\!endyarg}
  \!barfinish}

\def\!barfinish{%
  \!ifnextchar/{\!finish}{\!barcurve}}

\def\putrectangle{%
  \!ifnextchar<{\!putrectangle}{\!putrectangle<\!zpt,\!zpt> }}
\def\!putrectangle<#1,#2> corners at #3 #4 and #5 #6 {%
  \!xone=\!M{#3}\!xunit  \!xtwo=\!M{#5}\!xunit%
  \!yone=\!M{#4}\!yunit  \!ytwo=\!M{#6}\!yunit%
  \ifdim \!xtwo<\!xone
    \!dimenI=\!xone  \!xone=\!xtwo  \!xtwo=\!dimenI
  \fi
  \ifdim \!ytwo<\!yone
    \!dimenI=\!yone  \!yone=\!ytwo  \!ytwo=\!dimenI
  \fi
  \!dimenI=#1\relax  \advance\!xone\!dimenI  \advance\!xtwo\!dimenI
  \!dimenI=#2\relax  \advance\!yone\!dimenI  \advance\!ytwo\!dimenI
  \let\!MRect=\!M
  \!setdimenmode
  \!shaderectangle
  \!dimenI=.5\linethickness
  \advance \!xone  -\!dimenI
  \advance \!xtwo   \!dimenI
  \putrule from {\!xone} {\!yone} to {\!xtwo} {\!yone}
  \putrule from {\!xone} {\!ytwo} to {\!xtwo} {\!ytwo}
  \advance \!xone   \!dimenI
  \advance \!xtwo  -\!dimenI%
  \advance \!yone  -\!dimenI
  \advance \!ytwo   \!dimenI
  \putrule from {\!xone} {\!yone} to {\!xone} {\!ytwo}
  \putrule from {\!xtwo} {\!yone} to {\!xtwo} {\!ytwo}
  \let\!M=\!MRect
  \ignorespaces}

\def\shaderectangleson{%
  \def\!shaderectangle{\!!shaderectangle}%
  \ignorespaces}
\def\shaderectanglesoff{%
  \def\!shaderectangle{}%
  \ignorespaces}

\shaderectanglesoff

\def\!!shaderectangle{%
  \!dimenA=\!xtwo  \advance \!dimenA -\!xone
  \!dimenB=\!ytwo  \advance \!dimenB -\!yone
  \ifdim \!dimenA<\!dimenB
    \!startvshade (\!xone,\!yone,\!ytwo)
    \!lshade      (\!xtwo,\!yone,\!ytwo)
  \else
    \!starthshade (\!yone,\!xone,\!xtwo)
    \!lshade      (\!ytwo,\!xone,\!xtwo)
  \fi
  \ignorespaces}

\def\frame{%
  \!ifnextchar<{\!frame}{\!frame<\!zpt> }}
\long\def\!frame<#1> #2{%
  \beginpicture
    \setcoordinatesystem units <1pt,1pt> point at 0 0
    \put {#2} [Bl] at 0 0
    \!dimenA=#1\relax
    \!dimenB=\!wd \advance \!dimenB \!dimenA
    \!dimenC=\!ht \advance \!dimenC \!dimenA
    \!dimenD=\!dp \advance \!dimenD \!dimenA
    \let\!MFr=\!M
    \!setdimenmode
    \putrectangle corners at {-\!dimenA} {-\!dimenD} and {\!dimenB} {\!dimenC}
    \!setcoordmode
    \let\!M=\!MFr
  \endpicture
  \ignorespaces}

\def\rectangle <#1> <#2> {%
  \setbox0=\hbox{}\wd0=#1\ht0=#2\frame {\box0}}

\def\plot{%
  \!ifnextchar"{\!plotfromfile}{\!drawcurve}}
\def\!plotfromfile"#1"{%
  \expandafter\!drawcurve \input #1 /}

\def\setquadratic{%
  \let\!drawcurve=\!qcurve
  \let\!!Shade=\!!qShade
  \let\!!!Shade=\!!!qShade}

\def\setlinear{%
  \let\!drawcurve=\!lcurve
  \let\!!Shade=\!!lShade
  \let\!!!Shade=\!!!lShade}

\def\sethistograms{%
  \let\!drawcurve=\!hcurve}

\def\!qcurve #1 #2 {%
  \!start (#1,#2)
  \!Qjoin}
\def\!Qjoin#1 #2 #3 #4 {%
  \!qjoin (#1,#2) (#3,#4)             
  \!ifnextchar/{\!finish}{\!Qjoin}}

\def\!lcurve #1 #2 {%
  \!start (#1,#2)
  \!Ljoin}
\def\!Ljoin#1 #2 {%
  \!ljoin (#1,#2)                    
  \!ifnextchar/{\!finish}{\!Ljoin}}

\def\!finish/{\ignorespaces}

\def\!hcurve #1 #2 {%
  \edef\!hxS{#1}%
  \edef\!hyS{#2}%
  \!hjoin}
\def\!hjoin#1 #2 {%
  \putrectangle corners at {\!hxS} {\!hyS} and {#1} {#2}
  \edef\!hxS{#1}%
  \!ifnextchar/{\!finish}{\!hjoin}}

\def\vshade #1 #2 #3 {%
  \!startvshade (#1,#2,#3)
  \!Shadewhat}

\def\hshade #1 #2 #3 {%
  \!starthshade (#1,#2,#3)
  \!Shadewhat}

\def\!Shadewhat{%
  \futurelet\!nextchar\!Shade}
\def\!Shade{%
  \if <\!nextchar
    \def\!nextShade{\!!Shade}%
  \else
    \if /\!nextchar
      \def\!nextShade{\!finish}%
    \else
      \def\!nextShade{\!!!Shade}%
    \fi
  \fi
  \!nextShade}
\def\!!lShade<#1> #2 #3 #4 {%
  \!lshade <#1> (#2,#3,#4)                 
  \!Shadewhat}
\def\!!!lShade#1 #2 #3 {%
  \!lshade (#1,#2,#3)
  \!Shadewhat}
\def\!!qShade<#1> #2 #3 #4 #5 #6 #7 {%
  \!qshade <#1> (#2,#3,#4) (#5,#6,#7)      
  \!Shadewhat}
\def\!!!qShade#1 #2 #3 #4 #5 #6 {%
  \!qshade (#1,#2,#3) (#4,#5,#6)
  \!Shadewhat}

\setlinear

\def\setdashpattern <#1>{%
  \def\!Flist{}\def\!Blist{}\def\!UDlist{}%
  \!countA=0
  \!ecfor\!item:=#1\do{%
    \!dimenA=\!item\relax
    \expandafter\!rightappend\the\!dimenA\withCS{\\}\to\!UDlist%
    \advance\!countA  1
    \ifodd\!countA
      \expandafter\!rightappend\the\!dimenA\withCS{\!Rule}\to\!Flist%
      \expandafter\!leftappend\the\!dimenA\withCS{\!Rule}\to\!Blist%
    \else
      \expandafter\!rightappend\the\!dimenA\withCS{\!Skip}\to\!Flist%
      \expandafter\!leftappend\the\!dimenA\withCS{\!Skip}\to\!Blist%
    \fi}%
  \!leaderlength=\!zpt
  \def\!Rule##1{\advance\!leaderlength  ##1}%
  \def\!Skip##1{\advance\!leaderlength  ##1}%
  \!Flist%
  \ifdim\!leaderlength>\!zpt
  \else
    \def\!Flist{\!Skip{24in}}\def\!Blist{\!Skip{24in}}\ignorespaces
    \def\!UDlist{\\{\!zpt}\\{24in}}\ignorespaces
    \!leaderlength=24in
  \fi
  \!dashingon}

\def\!dashingon{%
  \def\!advancedashing{\!!advancedashing}%
  \def\!drawlinearsegment{\!lineardashed}%
  \def\!puthline{\!putdashedhline}%
  \def\!putvline{\!putdashedvline}%
  \ignorespaces}%
\def\!dashingoff{%
  \def\!advancedashing{\relax}%
  \def\!drawlinearsegment{\!linearsolid}%
  \def\!puthline{\!putsolidhline}%
  \def\!putvline{\!putsolidvline}%
  \ignorespaces}

\def\setdots{%
  \!ifnextchar<{\!setdots}{\!setdots<5pt>}}
\def\!setdots<#1>{%
  \!dimenB=#1\advance\!dimenB -\plotsymbolspacing
  \ifdim\!dimenB<\!zpt
    \!dimenB=\!zpt
  \fi
\setdashpattern <\plotsymbolspacing,\!dimenB>}

\def\setdotsnear <#1> for <#2>{%
  \!dimenB=#2\relax  \advance\!dimenB -.05pt
  \!dimenC=#1\relax  \!countA=\!dimenC
  \!dimenD=\!dimenB  \advance\!dimenD .5\!dimenC  \!countB=\!dimenD
  \divide \!countB  \!countA
  \ifnum 1>\!countB
    \!countB=1
  \fi
  \divide\!dimenB  \!countB
  \setdots <\!dimenB>}

\def\setdashes{%
  \!ifnextchar<{\!setdashes}{\!setdashes<5pt>}}
\def\!setdashes<#1>{\setdashpattern <#1,#1>}

\def\setdashesnear <#1> for <#2>{%
  \!dimenB=#2\relax
  \!dimenC=#1\relax  \!countA=\!dimenC
  \!dimenD=\!dimenB  \advance\!dimenD .5\!dimenC  \!countB=\!dimenD
  \divide \!countB  \!countA
  \ifodd \!countB
  \else
    \advance \!countB  1
  \fi
  \divide\!dimenB  \!countB
  \setdashes <\!dimenB>}

\def\setsolid{%
  \def\!Flist{\!Rule{24in}}\def\!Blist{\!Rule{24in}}%
  \def\!UDlist{\\{24in}\\{\!zpt}}%
  \!dashingoff}
\setsolid

\def\!divide#1#2#3{%
  \!dimenB=#1
  \!dimenC=#2
  \!dimenD=\!dimenB
  \divide \!dimenD \!dimenC
  \!dimenA=\!dimenD
  \multiply\!dimenD \!dimenC
  \advance\!dimenB -\!dimenD
  \!dimenD=\!dimenC
    \ifdim\!dimenD<\!zpt \!dimenD=-\!dimenD
  \fi
  \ifdim\!dimenD<64pt
    \!divstep[\!tfs]\!divstep[\!tfs]%
  \else
    \!!divide
  \fi
  #3=\!dimenA\ignorespaces}

\def\!!divide{%
  \ifdim\!dimenD<256pt
    \!divstep[64]\!divstep[32]\!divstep[32]%
  \else
    \!divstep[8]\!divstep[8]\!divstep[8]\!divstep[8]\!divstep[8]%
    \!dimenA=2\!dimenA
  \fi}

\def\!divstep[#1]{
  \!dimenB=#1\!dimenB
  \!dimenD=\!dimenB
    \divide \!dimenD by \!dimenC
  \!dimenA=#1\!dimenA
    \advance\!dimenA by \!dimenD%
  \multiply\!dimenD by \!dimenC
    \advance\!dimenB by -\!dimenD}

\def\Divide <#1> by <#2> forming <#3> {%
  \!divide{#1}{#2}{#3}}

\def\circulararc{%
  \ellipticalarc axes ratio 1:1 }

\def\ellipticalarc axes ratio #1:#2 #3 degrees from #4 #5 center at #6 #7 {%
  \!angle=#3pt\relax
  \ifdim\!angle>\!zpt
    \def\!sign{}
  \else
    \def\!sign{-}\!angle=-\!angle
  \fi
  \!xxloc=\!M{#6}\!xunit
  \!yyloc=\!M{#7}\!yunit
  \!xxS=\!M{#4}\!xunit
  \!yyS=\!M{#5}\!yunit
  \advance\!xxS -\!xxloc
  \advance\!yyS -\!yyloc
  \!divide\!xxS{#1pt}\!xxS 
  \!divide\!yyS{#2pt}\!yyS 
  \let\!MC=\!M
  \!setdimenmode
  \!xS=#1\!xxS  \advance\!xS\!xxloc
  \!yS=#2\!yyS  \advance\!yS\!yyloc
  \!start (\!xS,\!yS)%
  \!loop\ifdim\!angle>14.9999pt
    \!rotate(\!xxS,\!yyS)by(\!cos,\!sign\!sin)to(\!xxM,\!yyM)
    \!rotate(\!xxM,\!yyM)by(\!cos,\!sign\!sin)to(\!xxE,\!yyE)
    \!xM=#1\!xxM  \advance\!xM\!xxloc  \!yM=#2\!yyM  \advance\!yM\!yyloc
    \!xE=#1\!xxE  \advance\!xE\!xxloc  \!yE=#2\!yyE  \advance\!yE\!yyloc
    \!qjoin (\!xM,\!yM) (\!xE,\!yE)
    \!xxS=\!xxE  \!yyS=\!yyE
    \advance \!angle -15pt
  \repeat
  \ifdim\!angle>\!zpt
    \!angle=100.53096\!angle
    \divide \!angle 360 
    \!sinandcos\!angle\!!sin\!!cos
    \!rotate(\!xxS,\!yyS)by(\!!cos,\!sign\!!sin)to(\!xxM,\!yyM)
    \!rotate(\!xxM,\!yyM)by(\!!cos,\!sign\!!sin)to(\!xxE,\!yyE)
    \!xM=#1\!xxM  \advance\!xM\!xxloc  \!yM=#2\!yyM  \advance\!yM\!yyloc
    \!xE=#1\!xxE  \advance\!xE\!xxloc  \!yE=#2\!yyE  \advance\!yE\!yyloc
    \!qjoin (\!xM,\!yM) (\!xE,\!yE)
  \fi
  \let\!M=\!MC
  \ignorespaces}

\def\!rotate(#1,#2)by(#3,#4)to(#5,#6){%
  \!dimenA=#3#1\advance \!dimenA -#4#2
  \!dimenB=#3#2\advance \!dimenB  #4#1
  \divide \!dimenA 32  \divide \!dimenB 32
  #5=\!dimenA  #6=\!dimenB
  \ignorespaces}
\def\!sin{4.17684}
\def\!cos{31.72624}

\def\!sinandcos#1#2#3{%
 \!dimenD=#1
 \!dimenA=\!dimenD
 \!dimenB=32pt
 \!removept\!dimenD\!value
 \!dimenC=\!dimenD
 \!dimenC=\!value\!dimenC \divide\!dimenC by 64 
 \advance\!dimenB by -\!dimenC
 \!dimenC=\!value\!dimenC \divide\!dimenC by 96 
 \advance\!dimenA by -\!dimenC
 \!dimenC=\!value\!dimenC \divide\!dimenC by 128 
 \advance\!dimenB by \!dimenC%
 \!removept\!dimenA#2
 \!removept\!dimenB#3
 \ignorespaces}

\def\putrule#1from #2 #3 to #4 #5 {%
  \!xloc=\!M{#2}\!xunit  \!xxloc=\!M{#4}\!xunit%
  \!yloc=\!M{#3}\!yunit  \!yyloc=\!M{#5}\!yunit%
  \!dxpos=\!xxloc  \advance\!dxpos by -\!xloc
  \!dypos=\!yyloc  \advance\!dypos by -\!yloc
  \ifdim\!dypos=\!zpt
    \def\!!Line{\!puthline{#1}}\ignorespaces
  \else
    \ifdim\!dxpos=\!zpt
      \def\!!Line{\!putvline{#1}}\ignorespaces
    \else
       \def\!!Line{}
    \fi
  \fi
  \let\!ML=\!M
  \!setdimenmode
  \!!Line%
  \let\!M=\!ML
  \ignorespaces}

\def\!putsolidhline#1{%
  \ifdim\!dxpos>\!zpt
    \put{\!hline\!dxpos}#1[l] at {\!xloc} {\!yloc}
  \else
    \put{\!hline{-\!dxpos}}#1[l] at {\!xxloc} {\!yyloc}
  \fi
  \ignorespaces}

\def\!putsolidvline#1{%
  \ifdim\!dypos>\!zpt
    \put{\!vline\!dypos}#1[b] at {\!xloc} {\!yloc}
  \else
    \put{\!vline{-\!dypos}}#1[b] at {\!xxloc} {\!yyloc}
  \fi
  \ignorespaces}

\def\!hline#1{\hbox to #1{\leaders \hrule height\linethickness\hfill}}
\def\!vline#1{\vbox to #1{\leaders \vrule width\linethickness\vfill}}

\def\!putdashedhline#1{%
  \ifdim\!dxpos>\!zpt
    \!DLsetup\!Flist\!dxpos
    \put{\hbox to \!totalleaderlength{\!hleaders}\!hpartialpattern\!Rtrunc}
      #1[l] at {\!xloc} {\!yloc}
  \else
    \!DLsetup\!Blist{-\!dxpos}
    \put{\!hpartialpattern\!Ltrunc\hbox to \!totalleaderlength{\!hleaders}}
      #1[r] at {\!xloc} {\!yloc}
  \fi
  \ignorespaces}

\def\!putdashedvline#1{%
  \!dypos=-\!dypos
  \ifdim\!dypos>\!zpt
    \!DLsetup\!Flist\!dypos
    \put{\vbox{\vbox to \!totalleaderlength{\!vleaders}
      \!vpartialpattern\!Rtrunc}}#1[t] at {\!xloc} {\!yloc}
  \else
    \!DLsetup\!Blist{-\!dypos}
    \put{\vbox{\!vpartialpattern\!Ltrunc
      \vbox to \!totalleaderlength{\!vleaders}}}#1[b] at {\!xloc} {\!yloc}
  \fi
  \ignorespaces}

\def\!DLsetup#1#2{
  \let\!RSlist=#1
  \!countB=#2
  \!countA=\!leaderlength
  \divide\!countB by \!countA
  \!totalleaderlength=\!countB\!leaderlength
  \!Rresiduallength=#2%
  \advance \!Rresiduallength by -\!totalleaderlength
  \!Lresiduallength=\!leaderlength
  \advance \!Lresiduallength by -\!Rresiduallength
  \ignorespaces}

\def\!hleaders{%
  \def\!Rule##1{\vrule height\linethickness width##1}%
  \def\!Skip##1{\hskip##1}%
  \leaders\hbox{\!RSlist}\hfill}

\def\!hpartialpattern#1{%
  \!dimenA=\!zpt \!dimenB=\!zpt
  \def\!Rule##1{#1{##1}\vrule height\linethickness width\!dimenD}%
  \def\!Skip##1{#1{##1}\hskip\!dimenD}%
  \!RSlist}

\def\!vleaders{%
  \def\!Rule##1{\hrule width\linethickness height##1}%
  \def\!Skip##1{\vskip##1}%
  \leaders\vbox{\!RSlist}\vfill}

\def\!vpartialpattern#1{%
  \!dimenA=\!zpt \!dimenB=\!zpt
  \def\!Rule##1{#1{##1}\hrule width\linethickness height\!dimenD}%
  \def\!Skip##1{#1{##1}\vskip\!dimenD}%
  \!RSlist}

\def\!Rtrunc#1{\!trunc{#1}>\!Rresiduallength}
\def\!Ltrunc#1{\!trunc{#1}<\!Lresiduallength}

\def\!trunc#1#2#3{%
  \!dimenA=\!dimenB
  \advance\!dimenB by #1%
  \!dimenD=\!dimenB  \ifdim\!dimenD#2#3\!dimenD=#3\fi
  \!dimenC=\!dimenA  \ifdim\!dimenC#2#3\!dimenC=#3\fi
  \advance \!dimenD by -\!dimenC}

\def\!start (#1,#2){%
  \!plotxorigin=\!xorigin  \advance \!plotxorigin by \!plotsymbolxshift
  \!plotyorigin=\!yorigin  \advance \!plotyorigin by \!plotsymbolyshift
  \!xS=\!M{#1}\!xunit \!yS=\!M{#2}\!yunit
  \!rotateaboutpivot\!xS\!yS
  \!copylist\!UDlist\to\!!UDlist
  \!getnextvalueof\!downlength\from\!!UDlist
  \!distacross=\!zpt
  \!intervalno=0 
  \global\totalarclength=\!zpt
  \ignorespaces}

\def\!ljoin (#1,#2){%
  \advance\!intervalno by 1
  \!xE=\!M{#1}\!xunit \!yE=\!M{#2}\!yunit
  \!rotateaboutpivot\!xE\!yE
  \!xdiff=\!xE \advance \!xdiff by -\!xS
  \!ydiff=\!yE \advance \!ydiff by -\!yS
  \!Pythag\!xdiff\!ydiff\!arclength
  \global\advance \totalarclength by \!arclength%
  \!drawlinearsegment
  \!xS=\!xE \!yS=\!yE
  \ignorespaces}

\def\!linearsolid{%
  \!npoints=\!arclength
  \!countA=\plotsymbolspacing
  \divide\!npoints by \!countA
  \ifnum \!npoints<1
    \!npoints=1
  \fi
  \divide\!xdiff by \!npoints
  \divide\!ydiff by \!npoints
  \!xpos=\!xS \!ypos=\!yS
  \loop\ifnum\!npoints>-1
    \!plotifinbounds
    \advance \!xpos by \!xdiff
    \advance \!ypos by \!ydiff
    \advance \!npoints by -1
  \repeat
  \ignorespaces}

\def\!lineardashed{%
  \ifdim\!distacross>\!arclength
    \advance \!distacross by -\!arclength  
  \else
    \loop\ifdim\!distacross<\!arclength
      \!divide\!distacross\!arclength\!dimenA
      \!removept\!dimenA\!t
      \!xpos=\!t\!xdiff \advance \!xpos by \!xS
      \!ypos=\!t\!ydiff \advance \!ypos by \!yS
      \!plotifinbounds
      \advance\!distacross by \plotsymbolspacing
      \!advancedashing
    \repeat
    \advance \!distacross by -\!arclength
  \fi
  \ignorespaces}

\def\!!advancedashing{%
  \advance\!downlength by -\plotsymbolspacing
  \ifdim \!downlength>\!zpt
  \else
    \advance\!distacross by \!downlength
    \!getnextvalueof\!uplength\from\!!UDlist
    \advance\!distacross by \!uplength
    \!getnextvalueof\!downlength\from\!!UDlist
  \fi}

\def\inboundscheckoff{%
  \def\!plotifinbounds{\!plot(\!xpos,\!ypos)}%
  \def\!initinboundscheck{\relax}\ignorespaces}

\inboundscheckoff

\def\!!plotifinbounds{%
  \ifdim \!xpos<\!checkleft
  \else
    \ifdim \!xpos>\!checkright
    \else
      \ifdim \!ypos<\!checkbot
      \else
         \ifdim \!ypos>\!checktop
         \else
           \!plot(\!xpos,\!ypos)
         \fi
      \fi
    \fi
  \fi}

\def\!!initinboundscheck{%
  \!checkleft=\!arealloc     \advance\!checkleft by \!xorigin
  \!checkright=\!arearloc    \advance\!checkright by \!xorigin
  \!checkbot=\!areabloc      \advance\!checkbot by \!yorigin
  \!checktop=\!areatloc      \advance\!checktop by \!yorigin}

\def\!logten#1#2{%
  \expandafter\!!logten#1\!nil
  \!removept\!dimenF#2%
  \ignorespaces}

\def\!!logten#1#2\!nil{%
  \if -#1%
    \!dimenF=\!zpt
    \def\!next{\ignorespaces}%
  \else
    \if +#1%
      \def\!next{\!!logten#2\!nil}%
    \else
      \if .#1%
        \def\!next{\!!logten0.#2\!nil}%
      \else
        \def\!next{\!!!logten#1#2..\!nil}%
      \fi
    \fi
  \fi
  \!next}

\def\!!!logten#1#2.#3.#4\!nil{%
  \!dimenF=1pt 
  \if 0#1%
    \!!logshift#3pt 
  \else 
    \!logshift#2/
    \!dimenE=#1.#2#3pt 
  \fi 
  \ifdim \!dimenE<\!rootten
    \multiply \!dimenE 10 
    \advance  \!dimenF -1pt
  \fi
  \!dimenG=\!dimenE
    \advance\!dimenG 10pt
  \advance\!dimenE -10pt 
  \multiply\!dimenE 10 
  \!divide\!dimenE\!dimenG\!dimenE
  \!removept\!dimenE\!t
  \!dimenG=\!t\!dimenE
  \!removept\!dimenG\!tt
  \!dimenH=\!tt\!tenAe
    \divide\!dimenH 100
  \advance\!dimenH \!tenAc
  \!dimenH=\!tt\!dimenH
    \divide\!dimenH 100
  \advance\!dimenH \!tenAa
  \!dimenH=\!t\!dimenH
    \divide\!dimenH 100 
  \advance\!dimenF \!dimenH}

\def\!logshift#1{%
  \if #1/%
    \def\!next{\ignorespaces}%
  \else
    \advance\!dimenF 1pt
    \def\!next{\!logshift}%
  \fi
  \!next}

 \def\!!logshift#1{%
   \advance\!dimenF -1pt
   \if 0#1%
     \def\!next{\!!logshift}%
   \else
     \if p#1%
       \!dimenF=1pt
       \def\!next{\!dimenE=1p}%
     \else
       \def\!next{\!dimenE=#1.}%
     \fi
   \fi
   \!next}

\def\beginpicture{%
  \setbox\!picbox=\hbox\bgroup%
  \!xleft=\maxdimen
  \!xright=-\maxdimen
  \!ybot=\maxdimen
  \!ytop=-\maxdimen}

\def\endpicture{%
  \ifdim\!xleft=\maxdimen
    \!xleft=\!zpt \!xright=\!zpt \!ybot=\!zpt \!ytop=\!zpt
  \fi
  \global\!Xleft=\!xleft \global\!Xright=\!xright
  \global\!Ybot=\!ybot \global\!Ytop=\!ytop
  \egroup%
  \ht\!picbox=\!Ytop  \dp\!picbox=-\!Ybot
  \ifdim\!Ybot>\!zpt
  \else
    \ifdim\!Ytop<\!zpt
      \!Ybot=\!Ytop
    \else
      \!Ybot=\!zpt
    \fi
  \fi
  \hbox{\kern-\!Xleft\lower\!Ybot\box\!picbox\kern\!Xright}}

\def\endpicturesave <#1,#2>{%
  \endpicture \global #1=\!Xleft \global #2=\!Ybot \ignorespaces}

\def\setcoordinatesystem{%
  \!ifnextchar{u}{\!getlengths }
    {\!getlengths units <\!xunit,\!yunit>}}
\def\!getlengths units <#1,#2>{%
  \!xunit=#1\relax
  \!yunit=#2\relax
  \!ifcoordmode
    \let\!SCnext=\!SCccheckforRP
  \else
    \let\!SCnext=\!SCdcheckforRP
  \fi
  \!SCnext}
\def\!SCccheckforRP{%
  \!ifnextchar{p}{\!cgetreference }
    {\!cgetreference point at {\!xref} {\!yref} }}
\def\!cgetreference point at #1 #2 {%
  \edef\!xref{#1}\edef\!yref{#2}%
  \!xorigin=\!xref\!xunit  \!yorigin=\!yref\!yunit
  \!initinboundscheck 
  \ignorespaces}
\def\!SCdcheckforRP{%
  \!ifnextchar{p}{\!dgetreference}%
    {\ignorespaces}}
\def\!dgetreference point at #1 #2 {%
  \!xorigin=#1\relax  \!yorigin=#2\relax
  \ignorespaces}

\long\def\put#1#2 at #3 #4 {%
  \!setputobject{#1}{#2}%
  \!xpos=\!M{#3}\!xunit  \!ypos=\!M{#4}\!yunit
  \!rotateaboutpivot\!xpos\!ypos%
  \advance\!xpos -\!xorigin  \advance\!xpos -\!xshift
  \advance\!ypos -\!yorigin  \advance\!ypos -\!yshift
  \kern\!xpos\raise\!ypos\box\!putobject\kern-\!xpos%
  \!doaccounting\ignorespaces}

\long\def\multiput #1#2 at {%
  \!setputobject{#1}{#2}%
  \!ifnextchar"{\!putfromfile}{\!multiput}}
\def\!putfromfile"#1"{%
  \expandafter\!multiput \input #1 /}
\def\!multiput{%
  \futurelet\!nextchar\!!multiput}
\def\!!multiput{%
  \if *\!nextchar
    \def\!nextput{\!alsoby}%
  \else
    \if /\!nextchar
      \def\!nextput{\!finishmultiput}%
    \else
      \def\!nextput{\!alsoat}%
    \fi
  \fi
  \!nextput}
\def\!finishmultiput/{%
  \setbox\!putobject=\hbox{}%
  \ignorespaces}

\def\!alsoat#1 #2 {%
  \!xpos=\!M{#1}\!xunit  \!ypos=\!M{#2}\!yunit
  \!rotateaboutpivot\!xpos\!ypos%
  \advance\!xpos -\!xorigin  \advance\!xpos -\!xshift
  \advance\!ypos -\!yorigin  \advance\!ypos -\!yshift
  \kern\!xpos\raise\!ypos\copy\!putobject\kern-\!xpos%
  \!doaccounting
  \!multiput}

\def\!alsoby*#1 #2 #3 {%
  \!dxpos=\!M{#2}\!xunit \!dypos=\!M{#3}\!yunit
  \!rotateonly\!dxpos\!dypos
  \!ntemp=#1%
  \!!loop\ifnum\!ntemp>0
    \advance\!xpos by \!dxpos  \advance\!ypos by \!dypos
    \kern\!xpos\raise\!ypos\copy\!putobject\kern-\!xpos%
    \advance\!ntemp by -1
  \repeat
  \!doaccounting
  \!multiput}

\def\accountingon{\def\!doaccounting{\!!doaccounting}\ignorespaces}

\accountingon
\def\!!doaccounting{%
  \!xtemp=\!xpos
  \!ytemp=\!ypos
  \ifdim\!xtemp<\!xleft
     \!xleft=\!xtemp
  \fi
  \advance\!xtemp by  \!wd
  \ifdim\!xright<\!xtemp
    \!xright=\!xtemp
  \fi
  \advance\!ytemp by -\!dp
  \ifdim\!ytemp<\!ybot
    \!ybot=\!ytemp
  \fi
  \advance\!ytemp by  \!dp
  \advance\!ytemp by  \!ht
  \ifdim\!ytemp>\!ytop
    \!ytop=\!ytemp
  \fi}

\long\def\!setputobject#1#2{%
  \setbox\!putobject=\hbox{#1}%
  \!ht=\ht\!putobject  \!dp=\dp\!putobject  \!wd=\wd\!putobject
  \wd\!putobject=\!zpt
  \!xshift=.5\!wd   \!yshift=.5\!ht   \advance\!yshift by -.5\!dp
  \edef\!putorientation{#2}%
  \expandafter\!SPOreadA\!putorientation[]\!nil%
  \expandafter\!SPOreadB\!putorientation<\!zpt,\!zpt>\!nil\ignorespaces}

\def\!SPOreadA#1[#2]#3\!nil{\!etfor\!orientation:=#2\do\!SPOreviseshift}

\def\!SPOreadB#1<#2,#3>#4\!nil{\advance\!xshift by -#2\advance\!yshift by -#3}

\def\!SPOreviseshift{%
  \if l\!orientation
    \!xshift=\!zpt
  \else
    \if r\!orientation
      \!xshift=\!wd
    \else
      \if b\!orientation
        \!yshift=-\!dp
      \else
        \if B\!orientation
          \!yshift=\!zpt
        \else
          \if t\!orientation
            \!yshift=\!ht
          \fi
        \fi
      \fi
    \fi
  \fi}

\long\def\!dimenput#1#2(#3,#4){%
  \!setputobject{#1}{#2}%
  \!xpos=#3\advance\!xpos by -\!xshift
  \!ypos=#4\advance\!ypos by -\!yshift
  \kern\!xpos\raise\!ypos\box\!putobject\kern-\!xpos%
  \!doaccounting\ignorespaces}

\def\!setdimenmode{%
  \let\!M=\!M!!\ignorespaces}
\def\!setcoordmode{%
  \let\!M=\!M!\ignorespaces}
\def\!ifcoordmode{%
  \ifx \!M \!M!}
\def\!ifdimenmode{%
  \ifx \!M \!M!!}
\def\!M!#1#2{#1#2}
\def\!M!!#1#2{#1}
\!setcoordmode
\let\setdimensionmode=\!setdimenmode
\let\setcoordinatemode=\!setcoordmode

\def\!stack[#1]{%
  \let\!lglue=\hfill \let\!rglue=\hfill
  \expandafter\let\csname !#1glue\endcsname=\relax
  \!ifnextchar<{\!!stack}{\!!stack<\stackleading>}}
\def\!!stack<#1>#2{%
  \vbox{\def\!valueslist{}\!ecfor\!value:=#2\do{%
    \expandafter\!rightappend\!value\withCS{\\}\to\!valueslist}%
    \!lop\!valueslist\to\!value
    \let\\=\cr\lineskiplimit=\maxdimen\lineskip=#1%
    \baselineskip=-1000pt\halign{\!lglue##\!rglue\cr \!value\!valueslist\cr}}%
  \ignorespaces}

\def\!lines[#1]#2{%
  \let\!lglue=\hfill \let\!rglue=\hfill
  \expandafter\let\csname !#1glue\endcsname=\relax
  \vbox{\halign{\!lglue##\!rglue\cr #2\crcr}}%
  \ignorespaces}

\def\!Lines[#1]#2{%
  \let\!lglue=\hfill \let\!rglue=\hfill
  \expandafter\let\csname !#1glue\endcsname=\relax
  \vtop{\halign{\!lglue##\!rglue\cr #2\crcr}}%
  \ignorespaces}

\def\setplotsymbol(#1#2){%
  \!setputobject{#1}{#2}
  \setbox\!plotsymbol=\box\!putobject%
  \!plotsymbolxshift=\!xshift
  \!plotsymbolyshift=\!yshift
  \ignorespaces}

\setplotsymbol({\fiverm .})

\def\!!plot(#1,#2){%
  \!dimenA=-\!plotxorigin \advance \!dimenA by #1
  \!dimenB=-\!plotyorigin \advance \!dimenB by #2
  \kern\!dimenA\raise\!dimenB\copy\!plotsymbol\kern-\!dimenA%
  \ignorespaces}

\def\!!!plot(#1,#2){%
  \!dimenA=-\!plotxorigin \advance \!dimenA by #1
  \!dimenB=-\!plotyorigin \advance \!dimenB by #2
  \kern\!dimenA\raise\!dimenB\copy\!plotsymbol\kern-\!dimenA%
  \!countE=\!dimenA
  \!countF=\!dimenB
  \immediate\write\!replotfile{\the\!countE,\the\!countF.}%
  \ignorespaces}

\def\savelinesandcurves on "#1" {%
  \immediate\closeout\!replotfile
  \immediate\openout\!replotfile=#1%
  \let\!plot=\!!!plot}

\def\dontsavelinesandcurves {%
  \let\!plot=\!!plot}
\dontsavelinesandcurves

{\catcode`\%=11\xdef\!Commentsignal{
\def\writesavefile#1 {%
  \immediate\write\!replotfile{\!Commentsignal #1}%
  \ignorespaces}

\def\replot"#1" {%
  \expandafter\!replot\input #1 /}
\def\!replot#1,#2. {%
  \!dimenA=#1sp
  \kern\!dimenA\raise#2sp\copy\!plotsymbol\kern-\!dimenA
  \futurelet\!nextchar\!!replot}
\def\!!replot{%
  \if /\!nextchar
    \def\!next{\!finish}%
  \else
    \def\!next{\!replot}%
  \fi
  \!next}

\def\!Pythag#1#2#3{%
  \!dimenE=#1\relax
  \ifdim\!dimenE<\!zpt
    \!dimenE=-\!dimenE
  \fi
  \!dimenF=#2\relax
  \ifdim\!dimenF<\!zpt
    \!dimenF=-\!dimenF
  \fi
  \advance \!dimenF by \!dimenE
  \ifdim\!dimenF=\!zpt
    \!dimenG=\!zpt
  \else
    \!divide{8\!dimenE}\!dimenF\!dimenE
    \advance\!dimenE by -4pt
      \!dimenE=2\!dimenE
    \!removept\!dimenE\!!t
    \!dimenE=\!!t\!dimenE
    \advance\!dimenE by 64pt
    \divide \!dimenE by 2
    \!dimenH=7pt
    \!!Pythag\!!Pythag\!!Pythag
    \!removept\!dimenH\!!t
    \!dimenG=\!!t\!dimenF
    \divide\!dimenG by 8
  \fi
  #3=\!dimenG
  \ignorespaces}

\def\!!Pythag{
  \!divide\!dimenE\!dimenH\!dimenI
  \advance\!dimenH by \!dimenI
    \divide\!dimenH by 2}

\def\placehypotenuse for <#1> and <#2> in <#3> {%
  \!Pythag{#1}{#2}{#3}}

\def\!qjoin (#1,#2) (#3,#4){%
  \advance\!intervalno by 1
  \!ifcoordmode
    \edef\!xmidpt{#1}\edef\!ymidpt{#2}%
  \else
    \!dimenA=#1\relax \edef\!xmidpt{\the\!dimenA}%
    \!dimenA=#2\relax \edef\!xmidpt{\the\!dimenA}%
  \fi
  \!xM=\!M{#1}\!xunit  \!yM=\!M{#2}\!yunit   \!rotateaboutpivot\!xM\!yM
  \!xE=\!M{#3}\!xunit  \!yE=\!M{#4}\!yunit   \!rotateaboutpivot\!xE\!yE
  \!dimenA=\!xM  \advance \!dimenA by -\!xS
  \!dimenB=\!xE  \advance \!dimenB by -\!xM
  \!xB=3\!dimenA \advance \!xB by -\!dimenB
  \!xC=2\!dimenB \advance \!xC by -2\!dimenA
  \!dimenA=\!yM  \advance \!dimenA by -\!yS%
  \!dimenB=\!yE  \advance \!dimenB by -\!yM%
  \!yB=3\!dimenA \advance \!yB by -\!dimenB%
  \!yC=2\!dimenB \advance \!yC by -2\!dimenA%
  \!xprime=\!xB  \!yprime=\!yB
  \!dxprime=.5\!xC  \!dyprime=.5\!yC
  \!getf \!midarclength=\!dimenA
  \!getf \advance \!midarclength by 4\!dimenA
  \!getf \advance \!midarclength by \!dimenA
  \divide \!midarclength by 12
  \!arclength=\!dimenA
  \!getf \advance \!arclength by 4\!dimenA
  \!getf \advance \!arclength by \!dimenA
  \divide \!arclength by 12
  \advance \!arclength by \!midarclength
  \global\advance \totalarclength by \!arclength
  \ifdim\!distacross>\!arclength
    \advance \!distacross by -\!arclength
  \else
    \!initinverseinterp
    \loop\ifdim\!distacross<\!arclength
      \!inverseinterp
      \!xpos=\!t\!xC \advance\!xpos by \!xB
        \!xpos=\!t\!xpos \advance \!xpos by \!xS
      \!ypos=\!t\!yC \advance\!ypos by \!yB
        \!ypos=\!t\!ypos \advance \!ypos by \!yS
      \!plotifinbounds
      \advance\!distacross \plotsymbolspacing
      \!advancedashing
    \repeat
    \advance \!distacross by -\!arclength
  \fi
  \!xS=\!xE
  \!yS=\!yE
  \ignorespaces}

\def\!getf{\!Pythag\!xprime\!yprime\!dimenA%
  \advance\!xprime by \!dxprime
  \advance\!yprime by \!dyprime}

\def\!initinverseinterp{%
  \ifdim\!arclength>\!zpt
    \!divide{8\!midarclength}\!arclength\!dimenE
    \ifdim\!dimenE<\!wmin \!setinverselinear
    \else
      \ifdim\!dimenE>\!wmax \!setinverselinear
      \else
        \def\!inverseinterp{\!inversequad}\ignorespaces
         \!removept\!dimenE\!Ew
         \!dimenF=-\!Ew\!dimenE
         \advance\!dimenF by 32pt
         \!dimenG=8pt
         \advance\!dimenG by -\!dimenE
         \!dimenG=\!Ew\!dimenG
         \!divide\!dimenF\!dimenG\!beta
         \!gamma=1pt
         \advance \!gamma by -\!beta
      \fi
    \fi
  \fi
  \ignorespaces}

\def\!inversequad{%
  \!divide\!distacross\!arclength\!dimenG
  \!removept\!dimenG\!v
  \!dimenG=\!v\!gamma
  \advance\!dimenG by \!beta
  \!dimenG=\!v\!dimenG
  \!removept\!dimenG\!t}

\def\!setinverselinear{%
  \def\!inverseinterp{\!inverselinear}%
  \divide\!dimenE by 8 \!removept\!dimenE\!t
  \!countC=\!intervalno \multiply \!countC 2
  \!countB=\!countC     \advance \!countB -1
  \!countA=\!countB     \advance \!countA -1
  \wlog{\the\!countB th point (\!xmidpt,\!ymidpt) being plotted
    doesn't lie in the}%
  \wlog{ middle third of the arc between the \the\!countA th
    and \the\!countC th points:}%
  \wlog{ [arc length \the\!countA\space to \the\!countB]/[arc length
    \the \!countA\space to \the\!countC]=\!t.}%
  \ignorespaces}

\def\!inverselinear{%
  \!divide\!distacross\!arclength\!dimenG
  \!removept\!dimenG\!t}

\def\startrotation{%
  \let\!rotateaboutpivot=\!!rotateaboutpivot
  \let\!rotateonly=\!!rotateonly
  \!ifnextchar{b}{\!getsincos }%
    {\!getsincos by {\!cosrotationangle} {\!sinrotationangle} }}
\def\!getsincos by #1 #2 {%
  \edef\!cosrotationangle{#1}%
  \edef\!sinrotationangle{#2}%
  \!ifcoordmode
    \let\!ROnext=\!ccheckforpivot
  \else
    \let\!ROnext=\!dcheckforpivot
  \fi
  \!ROnext}
\def\!ccheckforpivot{%
  \!ifnextchar{a}{\!cgetpivot}%
    {\!cgetpivot about {\!xpivotcoord} {\!ypivotcoord} }}
\def\!cgetpivot about #1 #2 {%
  \edef\!xpivotcoord{#1}%
  \edef\!ypivotcoord{#2}%
  \!xpivot=#1\!xunit  \!ypivot=#2\!yunit
  \ignorespaces}
\def\!dcheckforpivot{%
  \!ifnextchar{a}{\!dgetpivot}{\ignorespaces}}
\def\!dgetpivot about #1 #2 {%
  \!xpivot=#1\relax  \!ypivot=#2\relax
  \ignorespaces}

\def\stoprotation{%
  \let\!rotateaboutpivot=\!!!rotateaboutpivot
  \let\!rotateonly=\!!!rotateonly
  \ignorespaces}

\def\!!rotateaboutpivot#1#2{%
  \!dimenA=#1\relax  \advance\!dimenA -\!xpivot
  \!dimenB=#2\relax  \advance\!dimenB -\!ypivot
  \!dimenC=\!cosrotationangle\!dimenA
    \advance \!dimenC -\!sinrotationangle\!dimenB
  \!dimenD=\!cosrotationangle\!dimenB
    \advance \!dimenD  \!sinrotationangle\!dimenA
  \advance\!dimenC \!xpivot  \advance\!dimenD \!ypivot
  #1=\!dimenC  #2=\!dimenD
  \ignorespaces}

\def\!!rotateonly#1#2{%
  \!dimenA=#1\relax  \!dimenB=#2\relax
  \!dimenC=\!cosrotationangle\!dimenA
    \advance \!dimenC -\!rotsign\!sinrotationangle\!dimenB
  \!dimenD=\!cosrotationangle\!dimenB
    \advance \!dimenD  \!rotsign\!sinrotationangle\!dimenA
  #1=\!dimenC  #2=\!dimenD
  \ignorespaces}
\def\!rotsign{}
\def\!!!rotateaboutpivot#1#2{\relax}
\def\!!!rotateonly#1#2{\relax}
\stoprotation

\def\!reverserotateonly#1#2{%
  \def\!rotsign{-}%
  \!rotateonly{#1}{#2}%
  \def\!rotsign{}%
  \ignorespaces}

\def\setshadegrid{%
  \!ifnextchar{s}{\!getspan }
    {\!getspan span <\!dshade>}}
\def\!getspan span <#1>{%
  \!dshade=#1\relax
  \!ifcoordmode
    \let\!GRnext=\!GRccheckforAP
  \else
    \let\!GRnext=\!GRdcheckforAP
  \fi
  \!GRnext}
\def\!GRccheckforAP{%
  \!ifnextchar{p}{\!cgetanchor }
    {\!cgetanchor point at {\!xshadesave} {\!yshadesave} }}
\def\!cgetanchor point at #1 #2 {%
  \edef\!xshadesave{#1}\edef\!yshadesave{#2}%
  \!xshade=\!xshadesave\!xunit  \!yshade=\!yshadesave\!yunit
  \ignorespaces}
\def\!GRdcheckforAP{%
  \!ifnextchar{p}{\!dgetanchor}%
    {\ignorespaces}}
\def\!dgetanchor point at #1 #2 {%
  \!xshade=#1\relax  \!yshade=#2\relax
  \ignorespaces}

\def\setshadesymbol{%
  \!ifnextchar<{\!setshadesymbol}{\!setshadesymbol<,,,> }}

\def\!setshadesymbol <#1,#2,#3,#4> (#5#6){%
  \!setputobject{#5}{#6}%
  \setbox\!shadesymbol=\box\!putobject%
  \!shadesymbolxshift=\!xshift \!shadesymbolyshift=\!yshift
  \!dimenA=\!xshift \advance\!dimenA \!smidge
  \!override\!dimenA{#1}\!lshrinkage%
  \!dimenA=\!wd \advance \!dimenA -\!xshift
    \advance\!dimenA \!smidge
    \!override\!dimenA{#2}\!rshrinkage
  \!dimenA=\!dp \advance \!dimenA \!yshift
    \advance\!dimenA \!smidge
    \!override\!dimenA{#3}\!bshrinkage
  \!dimenA=\!ht \advance \!dimenA -\!yshift
    \advance\!dimenA \!smidge
    \!override\!dimenA{#4}\!tshrinkage
  \ignorespaces}
\def\!smidge{-.2pt}%

\def\!override#1#2#3{%
  \edef\!!override{#2}%
  \ifx \!!override\empty
    #3=#1\relax
  \else
    \if z\!!override
      #3=\!zpt
    \else
      \ifx \!!override\!blankz
        #3=\!zpt
      \else
        #3=#2\relax
      \fi
    \fi
  \fi
  \ignorespaces}
\def\!blankz{ z}

\setshadesymbol ({\fiverm .})

\def\!startvshade#1(#2,#3,#4){%
  \let\!!xunit=\!xunit%
  \let\!!yunit=\!yunit%
  \let\!!xshade=\!xshade%
  \let\!!yshade=\!yshade%
  \def\!getshrinkages{\!vgetshrinkages}%
  \let\!setshadelocation=\!vsetshadelocation%
  \!xS=\!M{#2}\!!xunit
  \!ybS=\!M{#3}\!!yunit
  \!ytS=\!M{#4}\!!yunit
  \!shadexorigin=\!xorigin  \advance \!shadexorigin \!shadesymbolxshift
  \!shadeyorigin=\!yorigin  \advance \!shadeyorigin \!shadesymbolyshift
  \ignorespaces}

\def\!starthshade#1(#2,#3,#4){%
  \let\!!xunit=\!yunit%
  \let\!!yunit=\!xunit%
  \let\!!xshade=\!yshade%
  \let\!!yshade=\!xshade%
  \def\!getshrinkages{\!hgetshrinkages}%
  \let\!setshadelocation=\!hsetshadelocation%
  \!xS=\!M{#2}\!!xunit
  \!ybS=\!M{#3}\!!yunit
  \!ytS=\!M{#4}\!!yunit
  \!shadexorigin=\!xorigin  \advance \!shadexorigin \!shadesymbolxshift
  \!shadeyorigin=\!yorigin  \advance \!shadeyorigin \!shadesymbolyshift
  \ignorespaces}

\def\!lattice#1#2#3#4#5{%
  \!dimenA=#1
  \!dimenB=#2
  \!countB=\!dimenB
  \!dimenC=#3
  \advance\!dimenC -\!dimenA
  \!countA=\!dimenC
  \divide\!countA \!countB
  \ifdim\!dimenC>\!zpt
    \!dimenD=\!countA\!dimenB
    \ifdim\!dimenD<\!dimenC
      \advance\!countA 1 
    \fi
  \fi
  \!dimenC=\!countA\!dimenB
    \advance\!dimenC \!dimenA
  #4=\!countA
  #5=\!dimenC
  \ignorespaces}

\def\!qshade#1(#2,#3,#4)#5(#6,#7,#8){%
  \!xM=\!M{#2}\!!xunit
  \!ybM=\!M{#3}\!!yunit
  \!ytM=\!M{#4}\!!yunit
  \!xE=\!M{#6}\!!xunit
  \!ybE=\!M{#7}\!!yunit
  \!ytE=\!M{#8}\!!yunit
  \!getcoeffs\!xS\!ybS\!xM\!ybM\!xE\!ybE\!ybB\!ybC
  \!getcoeffs\!xS\!ytS\!xM\!ytM\!xE\!ytE\!ytB\!ytC
  \def\!getylimits{\!qgetylimits}%
  \!shade{#1}\ignorespaces}

\def\!lshade#1(#2,#3,#4){%
  \!xE=\!M{#2}\!!xunit
  \!ybE=\!M{#3}\!!yunit
  \!ytE=\!M{#4}\!!yunit
  \!dimenE=\!xE  \advance \!dimenE -\!xS
  \!dimenC=\!ytE \advance \!dimenC -\!ytS
  \!divide\!dimenC\!dimenE\!ytB
  \!dimenC=\!ybE \advance \!dimenC -\!ybS
  \!divide\!dimenC\!dimenE\!ybB
  \def\!getylimits{\!lgetylimits}%
  \!shade{#1}\ignorespaces}

\def\!getcoeffs#1#2#3#4#5#6#7#8{%
  \!dimenC=#4\advance \!dimenC -#2
  \!dimenE=#3\advance \!dimenE -#1
  \!divide\!dimenC\!dimenE\!dimenF
  \!dimenC=#6\advance \!dimenC -#4
  \!dimenH=#5\advance \!dimenH -#3
  \!divide\!dimenC\!dimenH\!dimenG
  \advance\!dimenG -\!dimenF
  \advance \!dimenH \!dimenE
  \!divide\!dimenG\!dimenH#8
  \!removept#8\!t
  #7=-\!t\!dimenE
  \advance #7\!dimenF
  \ignorespaces}

\def\!shade#1{%
  \!getshrinkages#1<,,,>\!nil
  \advance \!dimenE \!xS
  \!lattice\!!xshade\!dshade\!dimenE
    \!parity\!xpos
  \!dimenF=-\!dimenF
    \advance\!dimenF \!xE
  \!loop\!not{\ifdim\!xpos>\!dimenF}
    \!shadecolumn%
    \advance\!xpos \!dshade
    \advance\!parity 1
  \repeat
  \!xS=\!xE
  \!ybS=\!ybE
  \!ytS=\!ytE
  \ignorespaces}

\def\!vgetshrinkages#1<#2,#3,#4,#5>#6\!nil{%
  \!override\!lshrinkage{#2}\!dimenE
  \!override\!rshrinkage{#3}\!dimenF
  \!override\!bshrinkage{#4}\!dimenG
  \!override\!tshrinkage{#5}\!dimenH
  \ignorespaces}
\def\!hgetshrinkages#1<#2,#3,#4,#5>#6\!nil{%
  \!override\!lshrinkage{#2}\!dimenG
  \!override\!rshrinkage{#3}\!dimenH
  \!override\!bshrinkage{#4}\!dimenE
  \!override\!tshrinkage{#5}\!dimenF
  \ignorespaces}

\def\!shadecolumn{%
  \!dxpos=\!xpos
  \advance\!dxpos -\!xS
  \!removept\!dxpos\!dx
  \!getylimits
  \advance\!ytpos -\!dimenH
  \advance\!ybpos \!dimenG
  \!yloc=\!!yshade
  \ifodd\!parity
     \advance\!yloc \!dshade
  \fi
  \!lattice\!yloc{2\!dshade}\!ybpos%
    \!countA\!ypos
  \!dimenA=-\!shadexorigin \advance \!dimenA \!xpos
  \loop\!not{\ifdim\!ypos>\!ytpos}
    \!setshadelocation
    \!rotateaboutpivot\!xloc\!yloc%
    \!dimenA=-\!shadexorigin \advance \!dimenA \!xloc
    \!dimenB=-\!shadeyorigin \advance \!dimenB \!yloc
    \kern\!dimenA \raise\!dimenB\copy\!shadesymbol \kern-\!dimenA
    \advance\!ypos 2\!dshade
  \repeat
  \ignorespaces}

\def\!qgetylimits{%
  \!dimenA=\!dx\!ytC
  \advance\!dimenA \!ytB
  \!ytpos=\!dx\!dimenA
  \advance\!ytpos \!ytS
  \!dimenA=\!dx\!ybC
  \advance\!dimenA \!ybB
  \!ybpos=\!dx\!dimenA
  \advance\!ybpos \!ybS}

\def\!lgetylimits{%
  \!ytpos=\!dx\!ytB
  \advance\!ytpos \!ytS
  \!ybpos=\!dx\!ybB
  \advance\!ybpos \!ybS}

\def\!vsetshadelocation{
  \!xloc=\!xpos
  \!yloc=\!ypos}
\def\!hsetshadelocation{
  \!xloc=\!ypos
  \!yloc=\!xpos}

\def\!axisticks {%
  \def\!nextkeyword##1 {%
    \expandafter\ifx\csname !ticks##1\endcsname \relax
      \def\!next{\!fixkeyword{##1}}%
    \else
      \def\!next{\csname !ticks##1\endcsname}%
    \fi
    \!next}%
  \!axissetup
    \def\!axissetup{\relax}%
  \edef\!ticksinoutsign{\!ticksinoutSign}%
  \!ticklength=\longticklength
  \!tickwidth=\linethickness
  \!gridlinestatus
  \!setticktransform
  \!maketick
  \!tickcase=0
  \def\!LTlist{}%
  \!nextkeyword}

\def\ticksout{%
  \def\!ticksinoutSign{+}}

\ticksout

\def\nogridlines{%
  \def\!gridlinestatus{\!gridlinestoofalse}}
\nogridlines

\def\loggedticks{%
  \def\!setticktransform{\let\!ticktransform=\!logten}}
\def\unloggedticks{%
  \def\!setticktransform{\let\!ticktransform=\!donothing}}
\def\!donothing#1#2{\def#2{#1}}
\unloggedticks

\expandafter\def\csname !ticks/\endcsname{%
  \!not {\ifx \!LTlist\empty}
    \!placetickvalues
  \fi
  \def\!tickvalueslist{}%
  \def\!LTlist{}%
  \expandafter\csname !axis/\endcsname}

\def\!maketick{%
  \setbox\!boxA=\hbox{%
    \beginpicture
      \!setdimenmode
      \setcoordinatesystem point at {\!zpt} {\!zpt}
      \linethickness=\!tickwidth
      \ifdim\!ticklength>\!zpt
        \putrule from {\!zpt} {\!zpt} to
          {\!ticksinoutsign\!tickxsign\!ticklength}
          {\!ticksinoutsign\!tickysign\!ticklength}
      \fi
      \if!gridlinestoo
        \putrule from {\!zpt} {\!zpt} to
          {-\!tickxsign\!xaxislength} {-\!tickysign\!yaxislength}
      \fi
    \endpicturesave <\!Xsave,\!Ysave>}%
    \wd\!boxA=\!zpt}

\def\!ticksin{%
  \def\!ticksinoutsign{-}%
  \!maketick
  \!nextkeyword}

\def\!ticksout{%
  \def\!ticksinoutsign{+}%
  \!maketick
  \!nextkeyword}

\def\!tickslength<#1> {%
  \!ticklength=#1\relax
  \!maketick
  \!nextkeyword}

\def\!tickslong{%
  \!tickslength<\longticklength> }

\def\!ticksshort{%
  \!tickslength<\shortticklength> }

\def\!tickswidth<#1> {%
  \!tickwidth=#1\relax
  \!maketick
  \!nextkeyword}

\def\!ticksandacross{%
  \!gridlinestootrue
  \!maketick
  \!nextkeyword}

\def\!ticksbutnotacross{%
  \!gridlinestoofalse
  \!maketick
  \!nextkeyword}

\def\!tickslogged{%
  \let\!ticktransform=\!logten
  \!nextkeyword}

\def\!ticksunlogged{%
  \let\!ticktransform=\!donothing
  \!nextkeyword}

\def\!ticksunlabeled{%
  \!tickcase=0
  \!nextkeyword}

\def\!ticksnumbered{%
  \!tickcase=1
  \!nextkeyword}

\def\!tickswithvalues#1/ {%
  \edef\!tickvalueslist{#1! /}%
  \!tickcase=2
  \!nextkeyword}

\def\!ticksquantity#1 {%
  \ifnum #1>1
    \!updatetickoffset
    \!countA=#1\relax
    \advance \!countA -1
    \!ticklocationincr=\!axisLength
      \divide \!ticklocationincr \!countA
    \!ticklocation=\!axisstart
    \loop \!not{\ifdim \!ticklocation>\!axisend}
      \!placetick\!ticklocation
      \ifcase\!tickcase
          \relax 
        \or
          \relax 
        \or
          \expandafter\!gettickvaluefrom\!tickvalueslist
          \edef\!tickfield{{\the\!ticklocation}{\!value}}%
          \expandafter\!listaddon\expandafter{\!tickfield}\!LTlist%
      \fi
      \advance \!ticklocation \!ticklocationincr
    \repeat
  \fi
  \!nextkeyword}

\def\!ticksat#1 {%
  \!updatetickoffset
  \edef\!Loc{#1}%
  \if /\!Loc
    \def\next{\!nextkeyword}%
  \else
    \!ticksincommon
    \def\next{\!ticksat}%
  \fi
  \next}

\def\!ticksfrom#1 to #2 by #3 {%
  \!updatetickoffset
  \edef\!arg{#3}%
  \expandafter\!separate\!arg\!nil
  \!scalefactor=1
  \expandafter\!countfigures\!arg/
  \edef\!arg{#1}%
  \!scaleup\!arg by\!scalefactor to\!countE
  \edef\!arg{#2}%
  \!scaleup\!arg by\!scalefactor to\!countF
  \edef\!arg{#3}%
  \!scaleup\!arg by\!scalefactor to\!countG
  \loop \!not{\ifnum\!countE>\!countF}
    \ifnum\!scalefactor=1
      \edef\!Loc{\the\!countE}%
    \else
      \!scaledown\!countE by\!scalefactor to\!Loc
    \fi
    \!ticksincommon
    \advance \!countE \!countG
  \repeat
  \!nextkeyword}

\def\!updatetickoffset{%
  \!dimenA=\!ticksinoutsign\!ticklength
  \ifdim \!dimenA>\!offset
    \!offset=\!dimenA
  \fi}

\def\!placetick#1{%
  \if!xswitch
    \!xpos=#1\relax
    \!ypos=\!axisylevel
  \else
    \!xpos=\!axisxlevel
    \!ypos=#1\relax
  \fi
  \advance\!xpos \!Xsave
  \advance\!ypos \!Ysave
  \kern\!xpos\raise\!ypos\copy\!boxA\kern-\!xpos
  \ignorespaces}

\def\!gettickvaluefrom#1 #2 /{%
  \edef\!value{#1}%
  \edef\!tickvalueslist{#2 /}%
  \ifx \!tickvalueslist\!endtickvaluelist
    \!tickcase=0
  \fi}
\def\!endtickvaluelist{! /}

\def\!ticksincommon{%
  \!ticktransform\!Loc\!t
  \!ticklocation=\!t\!!unit
  \advance\!ticklocation -\!!origin
  \!placetick\!ticklocation
  \ifcase\!tickcase
    \relax 
  \or 
    \ifdim\!ticklocation<-\!!origin
      \edef\!Loc{$\!Loc$}%
    \fi
    \edef\!tickfield{{\the\!ticklocation}{\!Loc}}%
    \expandafter\!listaddon\expandafter{\!tickfield}\!LTlist%
  \or 
    \expandafter\!gettickvaluefrom\!tickvalueslist
    \edef\!tickfield{{\the\!ticklocation}{\!value}}%
    \expandafter\!listaddon\expandafter{\!tickfield}\!LTlist%
  \fi}

\def\!separate#1\!nil{%
  \!ifnextchar{-}{\!!separate}{\!!!separate}#1\!nil}
\def\!!separate-#1\!nil{%
  \def\!sign{-}%
  \!!!!separate#1..\!nil}
\def\!!!separate#1\!nil{%
  \def\!sign{+}%
  \!!!!separate#1..\!nil}
\def\!!!!separate#1.#2.#3\!nil{%
  \def\!arg{#1}%
  \ifx\!arg\!empty
    \!countA=0
  \else
    \!countA=\!arg
  \fi
  \def\!arg{#2}%
  \ifx\!arg\!empty
    \!countB=0
  \else
    \!countB=\!arg
  \fi}

\def\!countfigures#1{%
  \if #1/%
    \def\!next{\ignorespaces}%
  \else
    \multiply\!scalefactor 10
    \def\!next{\!countfigures}%
  \fi
  \!next}

\def\!scaleup#1by#2to#3{%
  \expandafter\!separate#1\!nil
  \multiply\!countA #2\relax
  \advance\!countA \!countB
  \if -\!sign
    \!countA=-\!countA
  \fi
  #3=\!countA
  \ignorespaces}

\def\!scaledown#1by#2to#3{%
  \!countA=#1\relax
  \ifnum \!countA<0 
    \def\!sign{-}
    \!countA=-\!countA
  \else
    \def\!sign{}%
  \fi
  \!countB=\!countA
  \divide\!countB #2\relax
  \!countC=\!countB
    \multiply\!countC #2\relax
  \advance \!countA -\!countC
  \edef#3{\!sign\the\!countB.}
  \!countC=\!countA 
  \ifnum\!countC=0 
    \!countC=1
  \fi
  \multiply\!countC 10
  \!loop \ifnum #2>\!countC
    \edef#3{#3\!zero}%
    \multiply\!countC 10
  \repeat
  \edef#3{#3\the\!countA}
  \ignorespaces}

\def\!placetickvalues{%
  \advance\!offset \tickstovaluesleading
  \if!xswitch
    \setbox\!boxA=\hbox{%
      \def\\##1##2{%
        \!dimenput {##2} [B] (##1,\!axisylevel)}%
      \beginpicture
        \!LTlist
      \endpicturesave <\!Xsave,\!Ysave>}%
    \!dimenA=\!axisylevel
      \advance\!dimenA -\!Ysave
      \advance\!dimenA \!tickysign\!offset
      \if -\!tickysign
        \advance\!dimenA -\ht\!boxA
      \else
        \advance\!dimenA  \dp\!boxA
      \fi
    \advance\!offset \ht\!boxA
      \advance\!offset \dp\!boxA
    \!dimenput {\box\!boxA} [Bl] <\!Xsave,\!Ysave> (\!zpt,\!dimenA)
  \else
    \setbox\!boxA=\hbox{%
      \def\\##1##2{%
        \!dimenput {##2} [r] (\!axisxlevel,##1)}%
      \beginpicture
        \!LTlist
      \endpicturesave <\!Xsave,\!Ysave>}%
    \!dimenA=\!axisxlevel
      \advance\!dimenA -\!Xsave
      \advance\!dimenA \!tickxsign\!offset
      \if -\!tickxsign
        \advance\!dimenA -\wd\!boxA
      \fi
    \advance\!offset \wd\!boxA
    \!dimenput {\box\!boxA} [Bl] <\!Xsave,\!Ysave> (\!dimenA,\!zpt)
  \fi}

\normalgraphs
\catcode`!=12 

\catcode`@=11 \catcode`!=11

\let\!pictexendpicture=\endpicture
\let\!pictexframe=\frame
\let\!pictexlinethickness=\linethickness
\let\!pictexmultiput=\multiput
\let\!pictexput=\put

\def\beginpicture{%
  \setbox\!picbox=\hbox\bgroup%
  \let\endpicture=\!pictexendpicture
  \let\frame=\!pictexframe
  \let\linethickness=\!pictexlinethickness
  \let\multiput=\!pictexmultiput
  \let\put=\!pictexput
  \let\input=\@@input   
  \!xleft=\maxdimen
  \!xright=-\maxdimen
  \!ybot=\maxdimen
  \!ytop=-\maxdimen}

\let\frame=\!latexframe

\let\pictexframe=\!pictexframe

\let\linethickness=\!latexlinethickness
\let\pictexlinethickness=\!pictexlinethickness

\let\\=\@normalcr
\catcode`@=12 \catcode`!=12

\newtheorem{prop}{Proposition}[section]
\newtheorem{lem}[prop]{Lemma}

\newtheorem{cor}[prop]{Corollary}
\newtheorem{them}[prop]{Theorem}

\newtheorem{question}[prop]{Question}
\newtheorem{conjecture}[prop]{Conjecture}

\theoremstyle{definition}

\newtheorem{defn}[prop]{Definition}

\newtheorem{numrmk}[prop]{Remark}
\newtheorem{numrmks}[prop]{Remarks}

\theoremstyle{remark}

\newtheorem{note}{Note}

\newtheorem{example}{Example}
\newtheorem{examples}{Examples}
\newtheorem{rmk}{Remark}
\newtheorem{rmks}{Remarks}

            {\nolinebreak $\Box$ \end{trivlist}}

\newcommand{\noprint}[1]{}

\newcommand{\rtext}[1]{\text{\normalshape #1}}

\renewcommand{\tilde}{\widetilde}

\newcommand{\cart}{{\mbox{\tiny cart}}}

\newcommand{\et}{\mbox{\tiny \'{e}t}}

\newcommand{\topo}{\mbox{\tiny top}}

\newcommand{\stab}{\mbox{\tiny stab}}
\newcommand{\virt}{\mbox{\tiny virt}}

\newcommand{\upst}{^{\ast}}
\newcommand{\op}{{\text{\normalshape op}}}
\newcommand{\upsh}{^{!}}
\newcommand{\lst}{_{\ast}}

\newcommand{\com}{^{\scriptscriptstyle\bullet}}

\renewcommand{\AA}{{\frak A}}

\newcommand{\TT}{{\frak T}}

\newcommand{\GG}{{\frak G}}

\newcommand{\MM}{{\frak M}}

\newcommand{\VV}{{\frak V}}
\newcommand{\WW}{{\frak W}}

\newcommand{\Aa}{{\frak a}}

\newcommand{\CC}{{\frak C}}
\newcommand{\zz}{{\Bbb Z}}
\newcommand{\hh}{{\Bbb H}}

\renewcommand{\ll}{{\Bbb L}}
\newcommand{\qq}{{\Bbb Q}}
\newcommand{\pp}{{\Bbb P}}
\newcommand{\cc}{{\Bbb C}}

\newcommand{\qql}{{{\Bbb Q}_\ell}}

\newcommand{\tT}{{\cal T}}

\renewcommand{\O}{{\cal O}}

\newcommand{\fF}{{\cal F}}
\newcommand{\lL}{{\cal L}}

\newcommand{\sS}{{\cal S}}

\newcommand{\del}{\partial}
\newcommand{\resto}{{ \mid }}
\newcommand{\st}{\mathrel{\mid}}

\newcommand{\rk}{\operatorname{\rm rk}}
\newcommand{\ev}{\operatorname{\rm ev}}

\newcommand{\cl}{\operatorname{\rm cl}\nolimits}

\newcommand{\Mor}{\operatorname{\rm Mor}\nolimits}

\newcommand{\ob}{\operatorname{ob}}

\newcommand{\chr}{\operatorname{\rm char}\nolimits}

\newcommand{\spec}{\operatorname{\rm Spec}\nolimits}

\newcommand{\proj}{\operatorname{\rm Proj}\nolimits}
\newcommand{\id}{\operatorname{\rm id}}
\newcommand{\Hom}{\operatorname{\rm Hom}\nolimits}

\newcommand{\Isomu}{\operatorname{\underline{\rm Isom}}\nolimits}

\newcommand{\Pic}{\operatorname{\rm Pic}\nolimits}

\renewcommand{\projlim}{\operatorname{{\lim\limits_{
\textstyle\longleftarrow}}}\limits}
\newcommand{\comp}{\mathbin{{\scriptstyle\circ}}}
\newcommand{\ol}{\overline}
\newcommand{\ul}{\underline}

\newcommand{\ldiag}[1]%
       {\makebox[0cm]{${\scriptstyle#1}\downarrow\phantom{\scriptstyle#1}$}}
\newcommand{\ldiagup}[1]%
       {\makebox[0cm]{${\scriptstyle#1}\uparrow\phantom{\scriptstyle#1}$}}
\newcommand{\rdiag}[1]%
       {\makebox[0cm]{$\phantom{\scriptstyle#1}\downarrow{\scriptstyle#1}$}}
\newcommand{\sediagr}[1]%
       {\makebox[0cm]{$\phantom{\scriptstyle#1}\searrow{\scriptstyle#1}$}}
\newcommand{\nediagr}[1]%
       {\makebox[0cm]{$\phantom{\scriptstyle#1}\nearrow{\scriptstyle#1}$}}
\newcommand{\rdiagup}[1]%
       {\makebox[0cm]{$\phantom{\scriptstyle#1}\uparrow{\scriptstyle#1}$}}
\newcommand{\swdiag}[1]%
       {\makebox[0cm]{$\phantom{\scriptstyle#1}\swarrow{\scriptstyle#1}$}}
\newcommand{\sediag}[1]%
       {\makebox[0cm]{${\scriptstyle#1}\searrow\phantom{\scriptstyle#1}$}}
\newcommand{\nediag}[1]%
       {\makebox[0cm]{${\scriptstyle#1}\nearrow\phantom{\scriptstyle#1}$}}

\newcommand{\iso}{\stackrel{\sim}{\rightarrow}}

\newcommand{\doublearrowstack}[2]%
{{{{\scriptstyle#1}\atop{\textstyle \longrightarrow}}\atop{{\textstyle
\longrightarrow}\atop{ \scriptstyle#2}}}}
\newcommand{\rightleftarrowstack}[2]%
{{{{\scriptstyle#1}\atop{\textstyle \longrightarrow}}\atop{{\textstyle
\longleftarrow}\atop{ \scriptstyle#2}}}}
\newcommand{\leftrightarrowstack}[2]%
{{{{\scriptstyle#1}\atop{\textstyle \longleftarrow}}\atop{{\textstyle
\longrightarrow}\atop{ \scriptstyle#2}}}}

\newcommand{\comdia}[9]{%
\begin{array}{ccc}
#1 & \stackrel{#2}{\longrightarrow} & #3 \\
\ldiag{#4} & #5 & \rdiag{#6} \\
#7 & \stackrel{#8}{\longrightarrow} & #9
\end{array}}
\newcommand{\comdiaback}[9]{%
\begin{array}{ccc}
#1 & \stackrel{#2}{\longleftarrow} & #3 \\
\ldiag{#4} & #5 & \rdiag{#6} \\
#7 & \stackrel{#8}{\longleftarrow} & #9
\end{array}}

\newcommand{\comtri}[6]{%
\begin{array}{ccc}
#1 & \stackrel{#2}{\longrightarrow} & #3         \\
   & \sediag{#4}                    & \rdiag{#5} \\
   &                                & #6
\end{array}}

\newcommand{\overtoparrow}%
{\makebox[0cm]{\beginpicture
\setcoordinatesystem units <.8cm,.4cm> point at 0 0
\setplotarea x from -3 to 3, y from 0 to 1
\setquadratic
\plot -3 0 0 1 3 0 /
\put{\vector(3,-1){0}}[Bl] at 3 0
\endpicture}}

\newcommand{\underbottomarrow}%
{\makebox[0cm]{\beginpicture
\setcoordinatesystem units <.8cm,.4cm> point at 0 0
\setplotarea x from -3 to 3, y from 0 to 1
\setquadratic
\plot -3 1 0 0 3 1 /
\put{\vector(3,1){0}}[Bl] at 3 1
\endpicture}}

\renewcommand{\t}{_{\tau}}
\newcommand{\s}{_{\sigma}}

\setcounter{secnumdepth}{1}
\setcounter{tocdepth}{2}

\begin{document}
\maketitle
\begin{abstract}
We construct the motivic tree-level system of Gromov-Witten invariants
for convex varieties.
\end{abstract}

\setcounter{section}{-1}
\section{Introduction}

Let $V$ be a projective algebraic manifold. In \cite{KM}, Sec.\ 2,
Gromov-Witten invariants of $V$ were described
axiomatically as a collection of linear maps
\[I^{V}_{g,n,\beta}:\ H\upst(V)^{\otimes n}\longrightarrow
H\upst(\ol{M}_{g,n},\qq),\quad \beta\in H_2(V,\zz)\]
satisfying certain axioms, and a program to construct them by algebro-geometric
(as opposed to symplectic) techniques
was suggested. The program is based upon Kontsevich's notion of a stable map
$(C,x_1,\dots ,x_n,f)$, $f:C\rightarrow V$.
This data consists of an algebraic curve $C$ with $n$ labeled points on it and
a map $f$ such that if an irreducible
component of $C$ is contracted by $f$ to a point, then this component together
with its special points is
Deligne-Mumford stable.  For more details, see \cite{K} and below.

\smallskip

The construction consists of three major steps.

\medskip

A. Construct an orbispace (or rather a stack) of stable maps
$\ol{M}_{g,n}(V,\beta)$ such that $g=\text{genus of }C$,
$f\lst{([C])}=\beta$, and its two morphisms to $V^n$ and $\ol{M}_{g,n}$. On the
level of points, these morphisms are
given respectively by
\begin{eqnarray*}
p:(C,x_1,\dots ,x_n,f) & \longmapsto & (f(x_1),\dots ,f(x_n)), \\
q:(C,x_1,\dots ,x_n,f) & \longmapsto & [(C,x_1,\dots ,x_n)]^{\stab},
\end{eqnarray*}
where the last expression means the  stabilization of $(C,x_1,\dots ,x_n)$.

\medskip

B. Construct a ``virtual fundamental class'' $[\ol{M}_{g,n}(V,\beta)]_{\virt}$,
or ``orientation'' (see
Definition~\ref{domb} below) and use it to define a correspondence in the Chow
ring $C^V_{g,n,\beta}\in A(V^n\times
\ol{M}_{g,n})$.

\smallskip

This step suggested in \cite{K} is quite subtle and has not been spelled out in
full detail. It can be bypassed for
$g=0$ and $V=G/P$ (generalized flag spaces) where the virtual class coincides
with the usual one (see \cite{KM}).

\smallskip

In general, it involves a definition of a new $\zz$-graded supercommutative
structure sheaf on $\ol{M}_{g,n}(V,\beta)$.
The virtual class is obtained as a product of the class of this sheaf and the
inverse Todd class of an appropriate
tangent complex. Geometrically, it serves as a general position argument
furnishing the Dimension Axiom of \cite{KM} and
replacing the deformation of the complex structure used in the symplectic
context.

\medskip

C. Use $C^V_{g,n,\beta}$ in order to construct the induced maps
$I^V_{g,n,\beta}$ on any cohomology satisfying some
version of the standard properties making it functorial on the category of
correspondences.

\medskip

In this approach, the main features of $I^V_{g,n,\beta}$ axiomatized in
\cite{KM} reflect functorial properties of
$\ol{M}_{g,n}(V,\beta)$ and the cotangent complex with respect to degenerations
of stable maps. In particular, the key
``Splitting Axiom'' (or Associativity Equations for $g=0$) expresses the
compatibility between the divisors at infinity
of $\ol{M}_{g,n}(V,\beta)$ and $\ol{M}_{g,n}.$

\smallskip

A neat way to organize this information is to introduce the category of marked
stable modular graphs indexing
degeneration types of stable maps and to treat various modular stacks
$\ol{M}_{g,n}(V,\beta)$ as values of this modular
functor on the simplest one-vertex graphs. Then the check of the axioms in
\cite{KM} essentially boils down to a
calculation of this functor on a family of generating morphisms and objects in
the graph category.

\smallskip

The degeneration type of $(C,x_1,\dots ,x_n, f)$ is described by the graph
whose vertices are the irreducible components
of $C$, edges are singular points of $C$, and tails (``one-vertex edges'') are
$x_1,\dots ,x_n.$ In addition, each
vertex is marked by the homology class in $V$ which is the $f$-image of the
fundamental class of the respective
component of $C$ and by the genus of the normalization of this component. The
description of morphisms is somewhat more
delicate, cf.\ Sec.\ 1 below.

\smallskip

This philosophy is an extension of the operadic picture which already gained
considerable importance from various
viewpoints.  In turn, it leads to a new notion of a $\Gamma$-operad as a
monoidal functor on an appropriate category
$\Gamma$ of graphs, and an algebra over an operad as a morphism of such
functors.  This approach will be developed
elsewhere (see \cite{KapMan}). It clarifies the origin of the proliferation of
the types of operads considered recently
(May's, Markl's, modular, cyclic, ...)

\smallskip

In Part~I of the present paper we treat in this way Step~A, stressing the
functoriality not only with respect to the
degeneration types with fixed $V$ but also with respect to $V$, expressed by
the change of the marking semigroup of
abstract non-negative homology classes. We hope also that our approach will
help to introduce quantum cohomology with
coefficients and to understand better the K\"unneth formula for quantum
cohomology from \cite{KM2}.

\smallskip

Part~II is devoted to Steps~B and~C for $g=0$ and convex manifolds $V$. The
formalism of orientation classes is
introduced axiomatically, but we did not attempt to justify the relevant claims
of \cite{K} in general.

\smallskip

A word of warning and apology is due. The reader will meet several different
categories of marked graphs in this paper
of which the most important are $\GG_s$ (cf.\ Definition~\ref{dmmsg}),
$\tilde{\GG}_s(A)$ (cf.\ Definition~\ref{ecisg}
and the preceding discussion) and $\tilde{\GG}_s(V)_{\cart}$ (cf.\
Definition~\ref{ceicv}). They differ mainly by their
classes of morphisms. Specifically, certain elementary arrows which are
combinatorially ``the same'', run in opposite
directions in different categories, which affects the whole structure of the
morphism semigroups. The reason is that
functorial properties of moduli stacks of maps considered {\em by themselves }
are different form the functorial
properties of their virtual fundamental classes treated {\em as
correspondences}. Since graphs are used mainly as a
bookkeeping device, their categorical properties must reflect this distinction.

\medskip

\subsection{Acknowledgements}

The authors are grateful to the Max-Planck-Institut f\"ur Mathematik in Bonn
where this work was done for support and
stimulating atmosphere.

\vfill\eject
\part{Stacks of Stable Maps}

\section{Graphs} \label{graphs}

\begin{defn} \label{dog}
A {\em graph }$\tau$ is a quadruple $(F\t,V\t,j\t,\del\t)$, where $F\t$ and
$V\t$ are finite sets,
$\del\t:F\t\rightarrow V\t$ is a map and $j\t:F\t\rightarrow F\t$ an
involution. We call $F\t$ the set of {\em flags},
$V\t$ the set of {\em vertices}, $S\t=\{f\in F\t\st j\t f=f\}$ the set of {\em
tails } and $E\t=\{\{f_1,f_2\}\subset
F\t\st\text{$f_2=j\t f_1$}\}$ the set of edges of $\tau$. For $v\in V\t$ let
$F\t(v)=\del^{-1}\t(v)$ and $|v|=\#
F\t(v)$, the {\em valence }of $v$.
\end{defn}

\begin{defn} \label{geomr}
Let $\tau$ be a graph. We define the {\em geometric realization $|\tau|$ of
$\tau$} as follows. We start with the
disjoint union of closed intervals and singletons
\[\coprod_{f\in F\t}[0,{\textstyle{1\over2}}]\amalg\coprod_{v\in V\t}\{|v|\}.\]
We denote the real number $x\in [0,{1\over2}]$ in the component indexed by
$f\in F\t$ by $x_f$.  Then for every $v\in
V\t$ we identify all elements of $\{0_f\st f\in F\t(v)\}$ with $|v|$ and for
every edge $\{f_1,f_2\}$ of $\tau$, we
identify ${1\over2}_{f_1}$ and ${{1\over2}}_{f_2}$. Finally, we remove for
every tail $f\in S\t$ the point
${1\over2}_f$. We consider $|\tau|$ as a topological space with base points
given by $\{|v|\st v\in V\t\}$, the {\em
vertices } of $|\tau|$. It should always be clear from the context whether
$|v|$ denotes the geometric realization of a
vertex or it's valence.
\end{defn}

\begin{defn} \label{doc}
Let $\tau$ and $\sigma$ be graphs. A {\em contraction
}$\phi:\tau\rightarrow\sigma$ is a pair of maps
$\phi^F:F\s\rightarrow F\t$ and $\phi_V:V\t\rightarrow V\s$ such that the
following conditions are satisfied.
\begin{enumerate}
\item $\phi^F$ is injective and $\phi_V$ is surjective.
\item The diagram
\[\begin{array}{ccc}
F\t & \stackrel{\del\t}{\longrightarrow} & V\t \\
\ldiagup{\phi^F} & & \rdiag{\phi_V} \\
F\s & \stackrel{\del\s}{\longrightarrow} & V\s
\end{array}\] commutes.
\item $\phi^F\comp j\s =j\t\comp \phi^F$, so that $\phi$ induces injections
$\phi^S:S\s\rightarrow S\t$ and
$\phi^E:E\s\rightarrow E\t$ on tails and edges.
\item $\phi^S$ is a bijection, so $F\t-\phi^F(F\s)$ consists entirely of edges,
the edges being contracted.
\item Call two vertices $v,w\in V\t$ {\em equivalent}, if there exists an $f\in
F\t-\phi^F(F\s)$ such that $f\in F\t(v)$
and $j\t f\in F\t(w)$. Then pass to the associated equivalence relation on
$V\t$. The map $\phi_V:V\t\rightarrow V\s$
induces a bijection $V\t/\mathop{\sim}\rightarrow V\s$.
\end{enumerate}
For a vertex $v\in V\s$ the graph whose set of flags is
\[\{f\in F\t\st \rtext{$f\not\in\phi^F(F\s)$ and $\phi_V(\del\t f)=v$}\},\]
whose set of vertices is $\phi_V^{-1}(v)$ and whose $j$ and $\del$ are obtained
from $j\t$ and $\del\t$ by restriction,
is called the {\em graph being contracted onto $v$}. If the graphs being
contracted have together exactly one edge, we
call $\phi$ an {\em elementary contraction}.
\end{defn}

\begin{numrmks} \label{rmc}
\begin{enumerate}
\item It is clear how to compose contractions, and that the composition of
contractions is a contraction.
\item If $\phi:\tau\rightarrow\sigma$ and $\phi':\tau\rightarrow\sigma'$ are
contractions with the same set of edges
being contracted, then there exists a unique isomorphism
$\psi:\sigma\rightarrow\sigma'$ such that
$\psi\comp\phi=\phi'$.
\item Every contraction is a composition of elementary contractions.
\item \label{rmcone} To carry out a construction for contractions of graphs,
which is compatible with composition of
contractions, it suffices to perform this construction for elementary
contractions and check that the construction is
independent of the order in which it is realized for two elementary
contractions.
\end{enumerate}
\end{numrmks}

\begin{defn} \label{domg}
A {\em modular graph } is a graph $\tau$ endowed with a map
$g\t:V\t\rightarrow\zz_{\geq0};v\mapsto g(v)$. The number
$g(v)$ is called the {\em genus }of the vertex $v$.
\end{defn}

We say that a semigroup $A$ has {\em indecomposable zero}, if $a+b=0$ implies
$a=0$ and $b=0$, for any two elements
$a,b\in A$.

\begin{defn} \label{doasmg}
Let $\tau$ be a modular graph and $A$ a semigroup with indecomposable zero. An
{\em $A$-structure } on $\tau$ is a map
$\alpha:V_{\tau}\rightarrow A$. The element $\alpha(v)$ is called the {\em
class } of the vertex $v$. The pair
$(\tau,\alpha)$ is called a {\em modular graph with $A$-structure }(or {\em
$A$-graph}, by abuse of language).

A {\em marked graph } is a pair $(A,\tau)$, where $A$ is a semigroup with
indecomposable zero and $\tau$ an $A$-graph.
\end{defn}

\begin{defn} \label{commor}
Let $(\sigma,\alpha)$ and $(\tau,\beta)$ be $A$-graphs. A {\em combinatorial
morphism
$a:(\sigma,\alpha)\rightarrow(\tau,\beta)$} is a pair of maps
$a_F:F\s\rightarrow F\t$ and $a_V:V\s\rightarrow V\t$,
satisfying the following conditions.
\begin{enumerate}
\item \label{commor1} The diagram
\[\comdia{F\s}{\del\s}{V\s}{a_F}{}{a_V}{F\t}{\del\t}{V\t}\]
commutes. In particular, for every $v\in V\s$, letting $w=a_V(v)$, we get an
induced map $a_{V,v}:F\s(v)\rightarrow
F\t(w)$.
\item With the notation of (\ref{commor1}), for every $v\in V\s$ the map
$a_{V,v}:F\s(v)\rightarrow F\t(w)$ is
injective.
\item \label{commor3} Let $f\in F\s$ and $\ol{f}=j\s(f)$. If $f\not=\ol{f}$,
there exists an $n\geq1$ and $2n$ (not
necessarily distinct) flags $f_1,\ldots,f_n,\ol{f}_1,\ldots,\ol{f}_n\in F\t$
such that
\begin{enumerate}
\item $f_1=a_F(f)$ and $\ol{f}_n=a_F(\ol{f})$,
\item $j\t(f_i)=\ol{f}_i$, for all $i=1,\ldots,n$,
\item $\del\t(\ol{f}_i)=\del\t(f_{i+1})$ for all $i=1,\ldots,n-1$,
\item for all $i=1,\ldots,n-1$ we have
\[\ol{f}_i\not= f_{i+1}\Longrightarrow\text{$g({v_i})=0$ and $\beta(v_i)=0$},\]
where $v_i=\del(\ol{f}_i)=\del(f_{i+1})$,
\end{enumerate}
\item for every $v\in V\s$ we have $\alpha(v)=\beta(a_V(v))$,
\item for every $v\in V\s$ we have $g({v})=g({a_V(v)})$.
\end{enumerate}

A {\em combinatorial morphism of marked graphs
}$(B,\sigma,\beta)\rightarrow(A,\tau,\alpha)$ is a pair $(\xi,a)$, where
$\xi:A\rightarrow B$ is a homomorphism of semigroups and
$a:(\sigma,\beta)\rightarrow(\tau,\xi\comp\alpha)$ is a
combinatorial morphism of $B$-graphs.

Usually, we will suppress the subscripts of $a$.
\end{defn}

\begin{rmks}
\begin{enumerate}
\item The composition of two combinatorial morphisms is again a combinatorial
morphism.
\item We say that a combinatorial morphism $a:\sigma\rightarrow\tau$ is {\em
complete}, if for every $v\in V\s$ the map
$a_{V,v}:F\s(v)\rightarrow F\t(a(v))$ is bijective. Examples of complete
combinatorial morphism are
\begin{enumerate}
\item the inclusion of a connected component,
\item the morphism $\sigma\rightarrow\tau$, where $\sigma$ is obtained from
$\tau$ by {\em cutting an edge}, i.e.\
changing $j$ in such a way as to turn a two element orbit into two one element
orbits.
\end{enumerate}
\item Let $\tau$ be an $A$-graph and $f\in S\t$ a tail of $\tau$. Let
$F\s=F\t-\{f\}$, $V\s=V\t$ and define $\del\s$ and
$j\s$ by restricting $\del\t$ and $j\t$. Then $\sigma$ is naturally an
$A$-graph called {\em obtained from $\tau$ by
forgetting the tail $f$}. There is a canonical combinatorial morphism
$\sigma\rightarrow\tau$.
\item Every combinatorial morphism $a:\sigma\rightarrow\tau$ is a composition
$a=b\comp c$, where $b$ is complete and
$c$ is a finite composition of morphisms forgetting tails. If $\sigma$ and
$\tau$ are stable (Definition~\ref{dosag}),
all intermediate graphs in such a factorization are stable.
\item Condition~(\ref{commor3}) of Definition~\ref{commor} can be rephrased in
a more geometric way---see the remark
after Proposition~\ref{cmmgpe}.
\end{enumerate}
\end{rmks}

\begin{defn}
A {\em contraction }$\phi:(\tau,\alpha)\rightarrow(\sigma,\beta)$ of $A$-graphs
is a contraction of graphs
$\phi:\tau\rightarrow\sigma$ such that for every $v\in V\s$ we have
\begin{enumerate}
\item \[g(v)=\sum_{w\in\phi_V^{-1}(v)}g(w)+\dim H^1(|\tau_v|),\]
where $\tau_v$ is the graph being contracted onto $v$,
\item \[\beta(v)=\sum_{w\in\phi_V^{-1}(v)}\alpha(w).\]
\end{enumerate}
\end{defn}

\begin{defn} \label{dosag}
A vertex $v$ of a modular graph with $A$-structure $(\tau,\alpha)$ is called
{\em stable}, if $\alpha(v)=0$ implies
$2g(v)+|v|\geq3$. Otherwise, $v$ is called {\em unstable}. The $A$-graph $\tau$
is called {\em stable}, if all its
vertices are stable.
\end{defn}

We now come to an important construction which we shall call {\em stable
pullback}. Consider the following setup.  We
suppose given a homomorphism of semigroups $\xi:A\rightarrow B$, a contraction
of $A$-graphs
$\phi:\sigma\rightarrow\tau$ and a combinatorial morphism
$a:(B,\rho)\rightarrow(A,\tau)$ of marked graphs.  Moreover,
we assume that $\rho$ is a {\em stable }$B$-graph. We shall construct a stable
$B$-graph $\pi$, together with a
contraction of $B$-graphs $\psi:\pi\rightarrow\rho$ and a combinatorial
morphism of marked graphs
$b:\pi\rightarrow\sigma$. This $B$-graph $\pi$ will be called the {\em stable
pullback }of $\rho$ under $\phi$.
\[\begin{array}{ccccc}
B & & \pi & \stackrel{\psi}{\longrightarrow} & \rho \\
\ldiagup{\xi} & \phantom{\longrightarrow} & \ldiag{b} &  & \rdiag{a} \\
A & & \sigma & \stackrel{\phi}{\longrightarrow} & \tau
\end{array}\]

According with Remark~\ref{rmc}(\ref{rmcone}), we shall assume that $\phi$ is
elementary and contracts the edge
$\{f,\ol{f}\}$ of $\sigma$. Let $v_1=\del\s(f)$, $v_2=\del\s(\ol{f})$ and
$v_0=\phi(v_1)=\phi(v_2)$.

{\em Case I (Contracting a loop). } In this case $v_1=v_2$. Let
$w_1,\ldots,w_n$ be the vertices of $\rho$ that map to
$v_0$ under $a$. Note that $g(w_i)\geq1$, since $g(v_0)\geq1$. Let $\pi$ be
equal to $\rho$ with a loop
$\{f_i,\ol{f}_i\}$ attached to $w_i$, for each $i=1,\ldots,n$ and
$g_\pi(w_i)=g_\rho(w_i)-1$. Clearly, $\pi$ is stable,
the drop in certain genera is made up for by the addition of flags. The
morphism $b:\pi\rightarrow\sigma$ is the obvious
combinatorial morphism mapping every one of the loops $\{f_i,\ol{f}_i\}$ to
$\{f,\ol{f}\}$. The contraction
$\psi:\pi\rightarrow\rho$ simply contracts all the added loops.

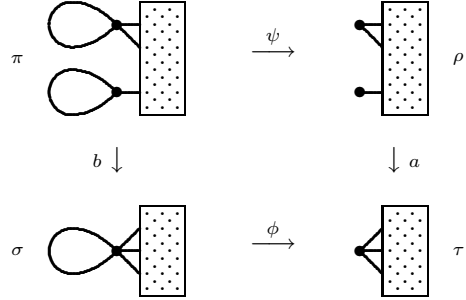
\begin{equation*}
\begin{array}{ccc}
\beginpicture
\setcoordinatesystem units <.3cm,.3cm> point at 3 2.5
\setplotarea x from 0 to 6, y from 0 to 5
\plot 3 4 4 4 /
\plot 3 4 4 3 /
\plot 3 1 4 1 /
\setquadratic
\plot 3 4 2 4.75 1 5 /
\plot 3 4 2 3.25 1 3 /
\plot 3 1 2 1.75 1 2 /
\plot 3 1 2 .25 1 0 /
\circulararc 180 degrees from 1 2 center at 1 1
\circulararc 180 degrees from 1 5 center at 1 4
\shaderectangleson
\setshadegrid span <1mm>
\putrectangle corners at 4 0 and 6 5
\put {\circle*{4}} [Bl] at 3 1
\put {\circle*{4}} [Bl] at 3 4
\axis left invisible label {$\scriptstyle\pi$} /
\axis right invisible label {\phantom{$\scriptstyle\pi$}} /
\axis bottom invisible label {\phantom{.}} /
\endpicture & \stackrel{\psi}{\longrightarrow} &
\beginpicture
\setcoordinatesystem units <.3cm,.3cm> point at 1.5 2.5
\setplotarea x from 0 to 3, y from 0 to 5
\plot 0 1 1 1 /
\plot 0 4 1 4 /
\plot 0 4 1 3 /
\shaderectangleson
\setshadegrid span <1mm>
\putrectangle corners at 1 0 and 3 5
\put {\circle*{4}} [Bl] at 0 1
\put {\circle*{4}} [Bl] at 0 4
\axis left invisible label {\phantom{$\scriptstyle\rho$}} /
\axis right invisible label {$\scriptstyle\rho$} /
\axis bottom invisible label {\phantom{.}} /
\endpicture \\
\ldiag{b} & & \rdiag{a} \\
\beginpicture
\setcoordinatesystem units <.3cm,.3cm> point at 3 2
\setplotarea x from 0 to 6, y from 0 to 4
\plot 3 2 4 3 /
\plot 3 2 4 2 /
\plot 3 2 4 1 /
\setquadratic
\plot 3 2 2 2.75 1 3 /
\plot 3 2 2 1.25 1 1 /
\circulararc 180 degrees from 1 3 center at 1 2
\shaderectangleson
\setshadegrid span <1mm>
\putrectangle corners at 4 0 and 6 4
\put {\circle*{4}} [Bl] at 3 2
\axis left invisible label {$\scriptstyle\sigma$} /
\axis right invisible label {\phantom{$\scriptstyle\sigma$}} /
\axis top invisible label {\phantom{.}} /
\endpicture & \stackrel{\phi}{\longrightarrow} &
\beginpicture
\setcoordinatesystem units <.3cm,.3cm> point at 1.5 2
\setplotarea x from 0 to 3, y from 0 to 4
\plot 0 2 1 3 /
\plot 0 2 1 2 /
\plot 0 2 1 1 /
\shaderectangleson
\setshadegrid span <1mm>
\putrectangle corners at 1 0 and 3 4
\put {\circle*{4}} [Bl] at 0 2
\axis left invisible label {\phantom{$\scriptstyle\tau$}} /
\axis right invisible label {$\scriptstyle\tau$} /
\axis top invisible label {\phantom{.}} /
\endpicture
\end{array}
\end{equation*}

{\em Case II (Contracting a non-looping edge). } In this case $v_1\not=v_2$.
Again, let $w_1,\ldots, w_n$ be the
vertices of $\rho$ that map to $v_0$ under $a$. First we shall construct an
intermediate graph $\pi'$. Let us fix an
$i=1,\ldots,n$. Construct $\pi'$ from $\rho$ by replacing $w_i$ with two
vertices $w_i'$ and $w_i''$, connected by an
edge $\{f_i,\ol{f}_i\}$, such that $\del(f_i)=w_i'$ and $\del(\ol{f}_i)=w_i''$.
Let $f$ be a flag of $\rho$ such that
$\del_{\rho}(f)=w_i$. If $\del\s\phi^F(a(f))=v_1$ we attach $f$ to $w_i'$ and
if $\del\s\phi^F(a(f))=v_2$ we attach $f$
to $w_i''$. Set $g(w_i')=g(v_1)$, $g(w_i'')=g(v_2)$,
$\beta(w_i')=\xi(\alpha(v_1))$ and
$\beta(w_i'')=\xi(\alpha(v_2))$. This defines the $B$-graph $\pi'$. The problem
with $\pi'$ is that it might not be
stable. So to construct $\pi$ we proceed as follows. Fix an $i=1,\ldots,n$. If
$w_i'$ and $w_i''$ are stable vertices of
$\pi'$ we do not change anything. If either of $w_i'$ or $w_i''$ is unstable,
we go back to where we started, by
contracting $\{f_i,\ol{f}_i\}$ again, obtaining the stable vertex $w_i$. This
finally finishes the construction of
$\pi$. The contraction $\psi:\pi\rightarrow\rho$ is defined by contracting all
the edges that were just inserted into
$\rho$ to construct $\pi$. There is an obvious combinatorial morphism
$b':\pi'\rightarrow\sigma$ mapping the edge
$\{f_i,\ol{f}_i\}$ to $\{f,\ol{f}\}$. Moreover, we define a combinatorial
morphism $c:\pi\rightarrow\pi'$ as follows. If
$i=1,\ldots,n$ is an index such that either of $w_i'$ or $w_i''$ is an unstable
vertex of $\pi'$, we map the vertex
$w_i$ of $\pi$ to the stable one of the two, say $w_i'$, to fix notation. If
$f$ is a flag of $\rho$ such that
$\del{f}=w_i$, then $f$ is also considered as a flag of $\pi$ and $\pi'$, and
under $c$ we map $f$ to itself, if
$\del_{\pi'}(f)=w_i'$, and to $f_i$, otherwise. Finally,
$b:\pi\rightarrow\sigma$ is defined as the composition
$b=b'\comp c$.

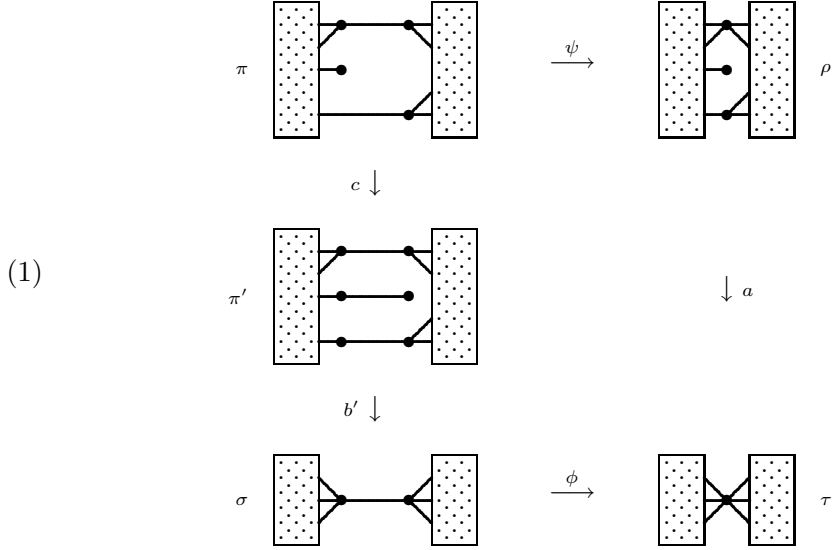
\begin{equation} \label{tadsp}
\begin{array}{ccc}
\beginpicture
\setcoordinatesystem units <.3cm,.3cm> point at 4.5 3
\setplotarea x from 0 to 9, y from 0 to 6
\plot 2 5 7 5 /
\plot 2 4 3 5 /
\plot 6 5 7 4 /
\plot 2 3 3 3 /
\plot 2 1 7 1 /
\plot 6 1 7 2 /
\shaderectangleson
\setshadegrid span <1mm>
\putrectangle corners at 0 0 and 2 6
\putrectangle corners at 7 0 and 9 6
\put {\circle*{4}} [Bl] at 3 5
\put {\circle*{4}} [Bl] at 6 5
\put {\circle*{4}} [Bl] at 3 3
\put {\circle*{4}} [Bl] at 6 1
\axis left invisible label {$\scriptstyle\pi$} /
\axis right invisible label {\phantom{$\scriptstyle\pi$}} /
\axis bottom invisible label {\phantom{.}} /
\endpicture & \stackrel{\psi}{\longrightarrow} &
\beginpicture
\setcoordinatesystem units <.3cm,.3cm> point at 3 3
\setplotarea x from 0 to 6, y from 0 to 6
\plot 2 5 4 5 /
\plot 2 4 3 5 /
\plot 3 5 4 4 /
\plot 3 3 2 3 /
\plot 2 1 4 1 /
\plot 3 1 4 2 /
\shaderectangleson
\setshadegrid span <1mm>
\putrectangle corners at 0 0 and 2 6
\putrectangle corners at 4 0 and 6 6
\put {\circle*{4}} [Bl] at 3 5
\put {\circle*{4}} [Bl] at 3 3
\put {\circle*{4}} [Bl] at 3 1
\axis left invisible label {\phantom{$\scriptstyle\rho$}} /
\axis right invisible label {$\scriptstyle\rho$} /
\axis bottom invisible label {\phantom{.}} /
\endpicture \\
\ldiag{c} & & \\
\beginpicture
\setcoordinatesystem units <.3cm,.3cm> point at 4.5 3
\setplotarea x from 0 to 9, y from 0 to 6
\plot 2 5 7 5 /
\plot 2 4 3 5 /
\plot 6 5 7 4 /
\plot 2 3 6 3 /
\plot 2 1 7 1 /
\plot 6 1 7 2 /
\shaderectangleson
\setshadegrid span <1mm>
\putrectangle corners at 0 0 and 2 6
\putrectangle corners at 7 0 and 9 6
\put {\circle*{4}} [Bl] at 3 5
\put {\circle*{4}} [Bl] at 6 5
\put {\circle*{4}} [Bl] at 3 3
\put {\circle*{4}} [Bl] at 6 3
\put {\circle*{4}} [Bl] at 3 1
\put {\circle*{4}} [Bl] at 6 1
\axis left invisible label {$\scriptstyle\pi'$} /
\axis right invisible label {\phantom{$\scriptstyle\pi'$}} /
\axis bottom invisible label {\phantom{.}} /
\axis top invisible label {\phantom{.}} /
\endpicture & & \rdiag{a} \\
\ldiag{b'} & & \\
\beginpicture
\setcoordinatesystem units <.3cm,.3cm> point at 4.5 2
\setplotarea x from 0 to 9, y from 0 to 4
\plot 2 2 7 2 /
\plot 2 3 3 2 2 1 /
\plot 7 3 6 2 7 1 /
\shaderectangleson
\setshadegrid span <1mm>
\putrectangle corners at 0 0 and 2 4
\putrectangle corners at 7 0 and 9 4
\put {\circle*{4}} [Bl] at 3 2
\put {\circle*{4}} [Bl] at 6 2
\axis left invisible label {$\scriptstyle\sigma$} /
\axis right invisible label {\phantom{$\scriptstyle\sigma$}} /
\axis top invisible label {\phantom{.}} /
\endpicture & \stackrel{\phi}{\longrightarrow} &
\beginpicture
\setcoordinatesystem units <.3cm,.3cm> point at 3 2
\setplotarea x from 0 to 6, y from 0 to 4
\plot 2 2 4 2 /
\plot 2 3 3 2 2 1 /
\plot 4 3 3 2 4 1 /
\shaderectangleson
\setshadegrid span <1mm>
\putrectangle corners at 0 0 and 2 4
\putrectangle corners at 4 0 and 6 4
\put {\circle*{4}} [Bl] at 3 2
\axis left invisible label {\phantom{$\scriptstyle\tau$}} /
\axis right invisible label {$\scriptstyle\tau$} /
\axis top invisible label {\phantom{.}} /
\endpicture
\end{array}
\end{equation}

Iterating this construction leads to the construction of a stable pullback for
arbitrary contractions of $A$-graphs.

\begin{rmk}
Let
\[\begin{array}{ccccc}
B & & \pi & \stackrel{\psi}{\longrightarrow} & \rho \\
\ldiagup{\xi} & \phantom{\longrightarrow} & \ldiag{b} &  & \rdiag{a} \\
A & & \sigma & \stackrel{\phi}{\longrightarrow} & \tau
\end{array}\]
be a stable pullback.
\begin{enumerate}
\item The diagram
\[\comdia{V_\pi}{\psi_V}{V_\rho}{b}{}{a}{V\s}{\phi_V}{V\t}\]
commutes.
\item The diagram
\[\comdiaback{F_\pi}{\psi^F}{F_\rho}{b}{}{a}{F\s}{\phi^F}{F\t}\]
does {\em not }commute (except for special cases, e.g.\ if the $B$-graph $\pi'$
constructed above is stable).
\end{enumerate}
\end{rmk}

\begin{prop} \label{cspcn}
Stable pullback is independent of the order in which $\phi$ is decomposed into
elementary contractions.
Moreover, if
\[\begin{array}{ccccc}
B & & \pi & \stackrel{\psi}{\longrightarrow} & \rho \\
\ldiagup{\xi} & \phantom{\longrightarrow} & \ldiag{b} &  & \rdiag{a} \\
A & & \sigma & \stackrel{\phi}{\longrightarrow} & \tau
\end{array}\]
and
\[\begin{array}{ccccc}
B & & \pi' & \stackrel{\psi'}{\longrightarrow} & \pi \\
\ldiagup{\xi} & \phantom{\longrightarrow} & \ldiag{b'} &  & \rdiag{b} \\
A & & \sigma' & \stackrel{\phi'}{\longrightarrow} & \sigma
\end{array}\]
are stable pullbacks, then
\[\begin{array}{ccccc}
B & & \pi' & \stackrel{\psi\comp\psi'}{\longrightarrow} & \rho \\
\ldiagup{\xi} & \phantom{\longrightarrow} & \ldiag{b'} &  & \rdiag{a} \\
A & & \sigma' & \stackrel{\phi\comp\phi'}{\longrightarrow} & \tau
\end{array}\]
is a stable pullback, too.
\end{prop}
\begin{pf}
To check that stable pullback is well-defined, it suffices by
Remark~\ref{rmc}(\ref{rmcone}) to check that the above
construction yields the same result for both orders in which two elementary
contractions can be composed. This is a
straightforward, though maybe slightly tedious calculation.  The compatibility
of stable pullback with compositions of
contractions follows trivially from the definition.
\end{pf}

\begin{prop} \label{cspcb}
If
\[\begin{array}{ccccc}
B & & \pi & \stackrel{\psi}{\longrightarrow} & \rho \\
\ldiagup{\xi} & \phantom{\longrightarrow} & \ldiag{b} &  & \rdiag{a} \\
A & & \sigma & \stackrel{\phi}{\longrightarrow} & \tau
\end{array}\]
and
\[\begin{array}{ccccc}
C & & \pi' & \stackrel{\chi}{\longrightarrow} & \rho' \\
\ldiagup{\eta} & \phantom{\longrightarrow} & \ldiag{b'} &  & \rdiag{a'} \\
B & & \pi & \stackrel{\psi}{\longrightarrow} & \rho
\end{array}\]
are stable pullbacks, then
\[\begin{array}{ccccc}
C & & \pi' & \stackrel{\chi}{\longrightarrow} & \rho' \\
\ldiagup{\eta\comp\xi} & \phantom{\longrightarrow} & \ldiag{b\comp b'} &  &
\rdiag{a\comp a'} \\
A & & \sigma & \stackrel{\phi}{\longrightarrow} & \tau
\end{array}\]
is a stable pullback, too.
\end{prop}
\begin{pf}
Of course, it suffices to consider the case that $\phi$ is an elementary
contraction. Then the claim follows immediately
from the construction.
\end{pf}

We are now ready to define the notion of {\em morphism } of marked stable
graphs.

\begin{defn} \label{dmmsg}
Let $(A,\tau)$ and $(B,\sigma)$ be marked stable graphs. A {\em morphism } from
$(A,\tau)$ to $(B,\sigma)$ is a
quadruple $(\xi,a,\tau',\phi)$, where $\xi:A\rightarrow B$ is a homomorphism of
semigroups, $\tau'$ is a stable
$B$-graph, $a:\tau'\rightarrow\tau$ makes $(\xi,a)$ a combinatorial morphism of
marked graphs, and
$\phi:\tau'\rightarrow\sigma$ is a contraction of $B$-graphs. We also say that
$(a,\tau',\phi)$ is a morphism of marked
stable graphs {\em covering }$\xi$.

Let $(\xi,a,\tau',\phi):(A,\tau)\rightarrow(B,\sigma)$ and
$(\eta,b,\sigma',\psi):(B,\sigma)\rightarrow(C,\rho)$ be
morphisms of stable marked graphs. Then we define the {\em composition }
$(\eta,b,\sigma',\psi)\comp(\xi,a,\tau',\phi):(A,\tau)\rightarrow(C,\rho)$ to
be $(\eta\xi,a c,\tau'',\psi\chi)$, where
$(c,\tau'',\chi)$ is the stable pullback of $\sigma'$ under $\phi$.
\[\begin{array}{ccccccc}
C & & \tau'' & \stackrel{\chi}{\longrightarrow} & \sigma' &
\stackrel{\psi}{\longrightarrow} & \rho \\
\ldiagup{\eta} & \phantom{\longrightarrow} & \ldiag{c} &  & \rdiag{b} & & \\
B & & \tau' & \stackrel{\phi}{\longrightarrow} & \sigma & & \\
\ldiagup{\xi} & & \ldiag{a} & & & & \\
A & & \tau & & & &
\end{array}\]
\end{defn}

\begin{rmks}
\begin{enumerate}
\item In reality a morphism is an isomorphy class of quadruples as in this
definition. But we shall always stick to the
abuse of language begun here.
\item The composition of morphisms is associative by Propositions~\ref{cspcn}
and~\ref{cspcb}.
\item Every combinatorial morphism of marked graphs whose source and target are
stable defines a morphism of marked
stable graphs, but in the {\em opposite } direction.
\item Every contraction of $A$-graphs whose source (and hence target) is stable
defines a morphism of marked stable
graphs (in the same direction).
\end{enumerate}
\end{rmks}

The category of stable marked graphs shall be denoted by $\GG_s$. Let $\AA$ be
the category of (additive) semigroups with
indecomposable zero element. By projecting onto the first component, we get a
functor $\Aa:\GG_s\rightarrow\AA$. For
$A\in\ob\AA$ let $\GG_s(A)$ be the fiber of $\Aa$ over $A$, i.e.\ the category
of stable $A$-graphs.

\begin{prop} \label{ppes}
Let $\tau$ be an $A$-graph. Then there exists a stable $A$-graph $\tau^s$,
together with a combinatorial morphism
$\tau^s\rightarrow\tau$, such that every combinatorial morphism
$\sigma\rightarrow\tau$, where $\sigma$ is a stable
$A$-graph, factors uniquely through $\tau^s$. We call $\tau^s$ the {\em
stabilization  }of $\tau$.
\end{prop}
\begin{pf}
let $\alpha$ denote the $A$-structure of $\tau$.

{\em Case I\@. }Assume that $\tau$ has a vertex $v_0$ such that
 $g(v_0)=0$,
 $\alpha(v_0)=0$,
 $v_0$ has a unique flag $f_1$, and
 $f_2:=j(f_1)\not=f_1$.
Let $\tau'\rightarrow\tau$ be the `subgraph' defined by
 $F_{\tau'}=F\t-\{f_1\}$,
 $V_{\tau'}=V\t-\{v_0\}$,
 $\del_{\tau'}=\del\t\resto F_{\tau'}$,
 $j_{\tau'}\resto F_{\tau'}-\{f_2\}=j\t\resto F_{\tau'}-\{f_2\}$ and
$j_{\tau'}(f_2)=f_2$.

{\em Case II\@. }Assume that $\tau$ has a vertex $v_0$ such that
 $g(v_0)=0$,
 $\alpha(v_0)=0$,
 $v_0$ has exactly two flags, $f_1$ and $f_2$,
 $f_1$ is a tail of $\tau$ and
 $f_3:=j(f_2)\not=f_2$.
Let $\tau'\rightarrow\tau$ be the `subgraph' defined by
 $F_{\tau'}=F\t-\{f_1,f_2\}$,
 $V_{\tau'}=V\t-\{v_0\}$,
 $\del_{\tau'}=\del\t\resto F_{\tau'}$,
 $j_{\tau'}\resto F_{\tau'}-\{f_3\}=j\t\resto F_{\tau'}-\{f_3\}$ and
$j_{\tau'}(f_3)=f_3$.

{\em Case III\@. }Assume that $\tau$ has a vertex $v_0$ such that
 $g(v_0)=0$,
 $\alpha(v_0)=0$,
 $v_0$ has exactly two flags, $f_1$ and $f_2$,
 $\ol{f}_1:=j(f_1)\not=f_1$ and $\ol{f}_2:=j(f_2)\not=f_2$.
Let $\tau'\rightarrow\tau$ be the `subgraph' defined by
 $F_{\tau'}=F\t-\{f_1,f_2\}$,
 $V_{\tau'}=V\t-\{v_0\}$,
 $\del_{\tau'}=\del\t\resto F_{\tau'}$,
 $j_{\tau'}\resto F_{\tau'}-\{\ol{f}_1,\ol{f}_2\}=j\t\resto
F_{\tau'}-\{\ol{f}_1,\ol{f}_2\}$ and $j_{\tau'}(\ol{f}_1)=\ol{f}_2$.

{\em Case IV\@. }Assume that $\tau$ has a vertex $v_0$ such that
 $2g(v_0)+|v_0|<3$,
 $\alpha(v_0)=0$ and
 $F\t(v_0)$ is a union of orbits of $j\t$.
Let $\tau'\rightarrow\tau$ be the `subgraph' defined by
 $F_{\tau'}=F\t-F\t(v_0)$,
 $V_{\tau'}=V\t-\{v_0\}$,
 $\del_{\tau'}=\del\t\resto F_{\tau'}$ and
 $j_{\tau'}=j\t\resto F_{\tau'}$.

In each of these four cases every combinatorial morphism
$\sigma\rightarrow\tau$, with $\sigma$ stable factors uniquely
through $\tau'$. By induction on the number of vertices of $\tau$, the graph
$\tau'$ has a stabilization, which is thus
also a stabilization of $\tau$. If $\tau$ has no vertices $v_0$ of the kind
covered by the above four cases, $\tau$ is stable
and $\tau$ itself may serve as stabilization of $\tau$.
\end{pf}

See Section~10 in \cite[Exp.\ ~VI]{sga1}) for the definition of {\em
cofibration } of categories.

\begin{prop} \label{goscf}
The functor $\Aa:\GG_s\rightarrow\AA$ is a cofibration.
\end{prop}
\begin{pf}
Let $\xi:A\rightarrow B$ be a homomorphism in $\AA$, and $(\tau,\alpha)$ a
stable $A$-graph. We need to construct a
stable $B$-graph $\sigma=\xi\lst\tau$, together with a morphism
$(a,\tau',\phi):(A,\tau)\rightarrow(B,\sigma)$ covering
$\xi$, with the following universal mapping property. Whenever
$\eta:B\rightarrow C$ is another homomorphism in $\AA$,
$\rho$ is a stable $C$-graph and $(b,\tau'',\psi):(A,\tau)\rightarrow(C,\rho)$
is a morphism covering $\eta\comp\xi$,
there exists a unique morphism $(c,\sigma',\chi):(B,\sigma)\rightarrow(C,\rho)$
covering $\eta$, such that
$(c,\sigma',\chi)\comp(a,\tau',\phi)=(b,\tau'',\psi)$, i.e.\ such that $\tau''$
is the stable pullback of $\sigma'$
under $\phi$.

In fact, it is not difficult to see that the stabilization of
$(\tau,\xi\comp\alpha)$ satisfies this universal mapping
property.
\end{pf}

\begin{numrmk} \label{pfabs}
Choosing a {\em clivage normalis\'e }(see Definition~7.1 in \cite[Exp.\
{}~VI]{sga1}) of
$\GG_s$ over $\AA$ amounts to choosing a pushforward functor
$\xi\lst:\GG_s(A)\rightarrow\GG_s(B)$ for any homomorphism
$\xi:A\rightarrow B$ in $\AA$. We may call $\xi\lst$ {\em stabilization } with
respect to $\xi$. If $B=\{0\}$, we speak
of {\em absolute stabilization }(or simply {\em stabilization}, if no confusion
seems likely to arise).
\end{numrmk}

\section{Prestable Curves}

We recall the definition of prestable curves. A morphism of prestable curves is
defined in such a way that it has degree
at most one and contracts at most rational components.

\begin{defn} \label{dopsc}
A {\em prestable curve }over the scheme $T$ is a flat proper morphism
$\pi:C\rightarrow T$ of schemes such that the
geometric fibers of $\pi$ are reduced, {\em connected}, one-dimensional and
have at most ordinary double points (nodes)
as singularities.  The {\em genus } of a prestable curve $C\rightarrow T$ is
the map $t\mapsto\dim H^1(C_t,\O_{C_t})$,
which is a locally constant function $g:T\rightarrow\zz_{\geq0}$. If $L$ is a
line bundle on $C$, then the {\em degree }
of $L$ is the locally constant function $\deg L:T\rightarrow\zz_{\geq0}$ given
by $t\mapsto\chi(L_t)+g-1$.

A {\em morphism }$p:C\rightarrow D$ of prestable curves over $T$ is a
$T$-morphism of schemes, such that for every
geometric point $t$ of $T$ we have
\begin{enumerate}
\item if $\eta$ is the generic point of an irreducible component of $D_t$, then
the fiber of $p_t$ over $\eta$ is a finite
$\eta$-scheme of degree at most one,
\item if $C'$ is the normalization of an irreducible component of $C_t$, then
$p_t(C')$ is a single point only if $C'$
is rational.
\end{enumerate}
\end{defn}

If $V$ is a scheme and $f:C\rightarrow V$ a morphism, then $L\mapsto\deg f\upst
L$ defines a locally constant function
$T\rightarrow\Hom_\zz(\Pic V,\zz)$ which we shall call the {\em homology class
} of $f$, by abuse of language, denoted
$f\lst[C]$.

If $V$ is a scheme admitting an ample invertible sheaf let
\[H_2(V)^+=\{\alpha\in\Hom_\zz(\Pic V,\zz)\st \text{$\alpha(L)\geq0$ whenever
$L$ is ample}\}.\]
Note that $H_2(V)^+$ is a semigroup with indecomposable zero. This is because
if $V$ admits an ample invertible sheaf
then $\Pic V$ is generated by ample invertible sheaves (see Remarque~4.5.9 in
\cite{ega2}). So if $f:C\rightarrow V$ is
a morphism from a prestable curve into $V$, then the homology class is a
locally constant function $T\rightarrow
H_2(V)^+$.

\begin{lem} \label{tloq}
Let $f:X\rightarrow Y$ be a proper surjective morphism of $T$-schemes such that
$f\lst\O_X=\O_Y$. Let $g:X\rightarrow U$
be another morphism of $T$-schemes, such that for every geometric point $t$ of
$T$ the map $g_t:X_t\rightarrow U_t$ is
constant (as a map of underlying Zariski topological spaces) on the fibers of
$f_t:X_t\rightarrow Y_t$. Then $g$ factors
uniquely through $f$.
\end{lem}
\begin{pf}
This follows easily, for example, from Lemma~8.11.1 in \cite{ega2}.
\end{pf}

\begin{cor} \label{zc}
Let $C$ be a prestable curve over $T$ and $f:C\rightarrow V$ a morphism, where
$V$ is a scheme admitting
an ample invertible sheaf. Then $f\lst[C]=0$ if and only if $f$ factors through
$T$.  \qed
\end{cor}

We shall need the following two results about gluing marked prestable curves at
the marks.
\begin{prop} \label{glue}
Let $T$ be a scheme and $C_1$, $C_2$ two prestable curves over $T$. Let $x_1\in
C_1(T)$ and $x_2\in C_2(T)$ be sections
such that for every geometric point $t$ of $T$ we have that $x_1(t)$ and
$x_2(t)$ are in the smooth locus of $C_{1,t}$
and $C_{2,t}$, respectively. Then there exists a prestable curve $C$ over $T$,
together with $T$-morphisms
$p_1:C_1\rightarrow C$ and $p_2:C_2\rightarrow C$, such that
\begin{enumerate}
\item $p_1(x_1)=p_2(x_2)$,
\item $C$ is universal among all $T$-schemes with this property.
\end{enumerate}
The curve $C$ is uniquely determined (up to unique isomorphism) and will be
called {\em obtained by gluing $C_1$ and
$C_2$ along the sections $x_1$ and $x_2$}, notation
\[C=C_1\amalg_{x_1,x_2}C_2.\]

If $u:S\rightarrow T$ is a morphism of schemes, then $C_S$ is the curve
obtained by gluing $C_{1,S}$ and $C_{2,S}$ along
$x_{1,S}$ and $x_{2,S}$. If $g_i$ is the genus of $C_i$, for $i=1,2$, then for
the genus $g$ of $C$ we have $g=g_1+g_2$.
If, for $i=1,2$, $f_i:C_i\rightarrow V$ is a morphism into a scheme such that
$f_1(x_1)=f_2(x_2)$, and $f:C\rightarrow
V$ is the induced morphism, we have $f\lst[C]={f_1}\lst[C_1]+{f_2}\lst[C_2]$ in
$\Hom_\zz(\Pic V,\zz)$.  \qed
\end{prop}
\begin{prop} \label{glue1}
Let $T$ be a scheme and $C$ a prestable curve over $T$. Let $x_1\in C(T)$ and
$x_2\in C(T)$ be sections
such that for every geometric point $t$ of $T$ we have that $x_1(t)$ and
$x_2(t)$ are in the smooth locus of $C_{t}$ and
$x_1(t)\not= x_2(t)$. Then there exists a prestable curve $\tilde{C}$ over $T$,
together with a $T$-morphism
$p:C\rightarrow \tilde{C}$, such that
\begin{enumerate}
\item $p(x_1)=p(x_2)$,
\item $\tilde{C}$ is universal among { all }$T$-schemes with this property.
\end{enumerate}
The curve $\tilde{C}$ is uniquely determined (up to unique isomorphism) and
will be called {\em obtained by gluing $C$
with itself along the sections $x_1$ and $x_2$}, notation
\[\tilde{C}=C/{x_1\sim x_2}.\]

If $u:S\rightarrow T$ is a morphism of schemes, then $(\tilde{C})_S$ is the
curve obtained by gluing $C_{S}$ with itself along
$x_{1,S}$ and $x_{2,S}$. If $g$ is the genus of $C$, then for the genus
$\tilde{g}$ of $\tilde{C}$ we have $\tilde{g}=g+1$.
If $f:C\rightarrow V$ is a morphism into a scheme such that $f(x_1)=f(x_2)$,
and $\tilde{f}:\tilde{C}\rightarrow
V$ is the induced morphism, we have $\tilde{f}\lst[\tilde{C}]={f}\lst[C]$ in
$\Hom_\zz(\Pic V,\zz)$. \qed
\end{prop}

\begin{defn} \label{mpc}
Let $\tau$ be a modular graph. A {\em $\tau$-marked prestable curve over $T$}
is a pair $(C,x)$, where $C=(C_v)_{v\in
V\t}$ is a family of prestable curves $\pi_v:C_v\rightarrow T$ and
$x=(x_i)_{i\in F\t}$ is a family of sections
$x_i:T\rightarrow C_{\del\t(i)}$, such that for every geometric point $t$ of
$T$ we have
\begin{enumerate}
\item $x_i(t)$ is in the smooth locus of $C_{\del\t(i),t}$, for all $i\in F\t$,
\item $x_i(t)\not=x_j(t)$, if $i\not=j$, for $i,j\in F\t$,
\item $g(C_{v,t})=g(v)$ for all $v\in V\t$.
\end{enumerate}

We define a {\em marked prestable curve over $T$} to be a triple $(\tau,C,x)$,
where $\tau$ is a modular graph and
$(C,x)$ a $\tau$-marked prestable curve over $T$.
\end{defn}

\section{Stable Maps}

We now come to the definition of stable maps, the central concept of this work,
which is due to Kontsevich.

Fix a field $k$ and let $\VV$ be the category of smooth projective (not
necessarily connected) varieties over $k$.
Consider the covariant functor
\begin{eqnarray*}
H_2^+:\VV & \longrightarrow & \AA \\
V & \longmapsto & H_2(V)^+,
\end{eqnarray*}
where $\AA$ is the category of semigroups with indecomposable zero (see
Section~\ref{graphs}).
Define the category $\VV\GG_s$ as the fibered product (see Section~3 in
\cite[Exp.\ ~VI]{sga1})
\[\comdia{\VV\GG_s}{}{\GG_s}{}{\Box}{\Aa}{\VV}{H_2^+}{\AA.}\]
To spell this definition out, we have
\begin{enumerate}
\item objects of $\VV\GG_s$ are pairs $(V,\tau)$, where $V$ is a smooth
projective variety over $k$
and $\tau$ is a stable $H_2(V)^+$-graph,
\item a morphism $(V,\tau)\rightarrow (W,\sigma)$ is a quadruple
$(\xi,a,\tau',\phi)$, where
$\xi:V\rightarrow W$ is a morphism of $k$-varieties and
$(H_2^+(\xi),a,\tau',\phi)$ is a morphism in $\GG_s$ as defined in
Definition~\ref{dmmsg}.
\end{enumerate}

\begin{numrmk} \label{ugscf}
By Corollary~6.9 of \cite[Exp.\ ~VI]{sga1} and Proposition~\ref{goscf} the
category $\VV\GG_s$ is a cofibered category
over $\VV$.
\end{numrmk}

\begin{defn} \label{dsm}
Let $(V,\tau,\alpha)$ be an object of $\VV\GG_s$ and $T$ a $k$-scheme. A {\em
stable $(V,\tau,\alpha)$-map over $T$} is
a triple $(C,x,f)$, where $(C,x)$ is a $\tau$-marked prestable curve over $T$
and $f=(f_v)_{v\in V\t}$ is a family of
$k$-morphisms $f_v:C_v\rightarrow V$, such that the following conditions are
satisfied.
\begin{enumerate}
\item For every $i\in F\t$ we have
$f_{\del\t(i)}(x_i)=f_{\del\t(j\t(i))}(x_{j\t(i)})$ as $k$-morphisms from $T$
to $V$.
\item For all $v\in V\t$ we have that ${f_v}\lst[C_{v}]=\alpha(v)$ in
$H_2(V)^+$.
\item For every geometric point $t$ of $T$ and every $v\in V\t$ the {\em
stability condition }is satisfied. This means that if
$C'$ is the normalization of a component of $C_{v,t}$ that maps to a point
under $f_{v,t}:C_{v,t}\rightarrow V_t$, then
\begin{enumerate}
\item if the genus of $C'$ is zero, then $C'$ has at least three special
points,
\item if the genus of $C'$ is one, then $C'$ has at least one special point.
\end{enumerate}
Here, a point of $C'$ is called {\em special}, if it maps in $C_{v,t}$ to a
marked point or a node.
\end{enumerate}

We define a {\em stable map over $T$} to be a sextuple $(V,\tau,\alpha,C,x,f)$,
where $(V,\tau,\alpha)$ is an object
of $\VV\GG_s$ and $(C,x,f)$ is a stable $(V,\tau,\alpha)$-map over $T$.

A {\em morphism }$(V,\tau,\alpha,C,x,f)\rightarrow (W,\sigma,\beta,D,y,h)$ of
stable maps over $T$ is a quintuple
$(\xi,a,\tau',\phi,p)$, where $(\xi,a,\tau',\phi):(V,\tau,\alpha)\rightarrow
(W,\sigma,\beta)$ is a morphism in
$\VV\GG_s$ and $p=(p_v)_{v\in V_{\tau'}}$ is a family of morphisms of prestable
curves $p_v:C_{a(v)}\rightarrow
D_{\phi_V(v)}$, such that the following are true.
\begin{enumerate}
\item \label{mpc1}For every $i\in F\s$ we have
$p_{\del(\phi^F(i))}(x_{a\phi^F(i)})=y_i$,
\item \label{pre} If $\{i_1,i_2\}$ is an edge of $\tau'$ which is being
contracted by $\phi$, then
$p_{v_1}(x_{a(i_1)})=p_{v_2}(x_{a(i_2)})$, where $v_1=\del i_1$ and $v_2=\del
i_2$. So, in particular, if $v_1\not=v_2$
there exists an induces morphism
\[p_{12}:C_{a(v_1)}\amalg_{x_{a(i_1)},x_{a(i_2)}}C_{a(v_2)}\rightarrow D_w,\]
where $w=\phi(v_1)=\phi(v_2)$.
\item \label{post} With the notation of (\ref{pre}), if $v_1\not=v_2$, the
morphism $p_{12}$ is a morphism of prestable
curves.
\item For every $v\in V_{\tau'}$ the diagram
\[\comdia{C_{a(v)}}{f_{a(v)}}{V}{p_v}{}{\xi}{D_{\phi(v)}}{h_{\phi(v)}}{W}\]
commutes.
\end{enumerate}
In this situation we also say that $p:(C,x,f)\rightarrow(D,y,h)$ is a morphism
of stable maps {\em covering
}the morphism $(\xi,a,\tau',\alpha)$ in $\VV\GG_s$.

To define the {\em composition } of morphisms, let
$(\xi,a,\tau',\phi,p):(V,\tau,\alpha,C,x,f)\rightarrow(W,\sigma,\beta,D,y,h)$
and
$(\eta,b,\sigma',\psi,q):(W,\sigma,\beta,D,y,h)\rightarrow(U,\rho,
\gamma,E,z,e)$ be morphisms of stable maps over
$T$. We already know how to compose the morphisms $(\xi,a,\tau',\phi)$ and
$(\eta,b,\sigma',\psi)$ in $\VV\GG_s$. Use
notation as in Definition~\ref{dmmsg}. Then this composition is
$(\eta\xi,ac,\tau'',\psi\chi)$. Define the family
$r=(r_u)_{u\in V_{\tau''}}$ of morphisms of prestable curves
$r_u:C_{ac(u)}\rightarrow E_{\psi\chi(u)}$ as
$r_u=q_{\chi(u)}\comp p_{c(u)}$, which is well-defined, since
$\phi_Vc(u)=a\chi_V(u)$. Then we define our composition as
\[(\eta,b,\sigma',\psi,q)\comp(\xi,a,\tau',\phi,p) =
(\eta\xi,ac,\tau'',\psi\chi,r).\]
\end{defn}

\begin{prop} \label{cmsmm}
The composition of morphisms of stable maps is a morphism of stable maps.
\end{prop}
\begin{pf}
The proof will be given at the same time as the proof of Theorem~\ref{mbfc}
below.
\end{pf}

\begin{defn} \label{sovgs}
Let $V\in\ob\VV$ be a variety, $\beta\in H_2(V)^+$ a homology class and
$g,n\geq0$ integers. Then $(V,g,n,\beta)$ shall
denote the object $(V,\tau,\beta)$ of $\VV\GG_s$ whose modular graph $\tau$ is
given by $F\t=\ul{n}$,
$V\t=\{\varnothing\}$, $\del\t:F\t\rightarrow V\t$ the unique map,
$j\t=\id_{\ul{n}}$ and $g(\varnothing)=g$. The
$H_2(V)^+$-structure on $\tau$ is given by $\beta(\varnothing)=\beta$. A stable
$(V,g,n,\beta)$-map is also called a
{\em stable map from an $n$-pointed curve (of genus $g$) to $V$ (of class
$\beta$)}. Here we use the notation
$\ul{n}=\{1,\ldots,n\}$.
\end{defn}

\begin{lem} \label{lcsn}
Over an algebraically closed field, let $(C,x,f)$ be a stable map from an
$n$-pointed curve of genus $g$ to $V$ of class
$\beta$ and let $(D,y,h)$ be a stable map from an $m$-pointed curve of genus
$g$ to $V$ of class $\beta$, where $m\leq
n$. Let $p:C\rightarrow D$ be a morphism such that $p(x_i)=y_i$ for $i\leq m$
and $hp=f$. If $C'\subset C$ is a
subcurve (a connected union of irreducible components), such that
\begin{enumerate}
\item letting $C''$ be the closure of the complement of $C'$ in $C$, the curves
$C'$ and $C''$ have exactly one node in
common,
\item $g(C')=0$,
\item $f(C')$ is a point,
\item for $i\leq m$ the $x_i$ do not lie on $C'$ except for at most one of
them,
\end{enumerate}
then $p$ maps $C'$ to a point in $D$. \qed
\end{lem}

Let us denote the category of stable maps over $T$ by $\ol{M}(T)$. It comes
together with a functor
\[\ol{M}(T)\longrightarrow\VV\GG_s,\]
defined by projecting onto the first components. For a morphism $u:S\rightarrow
T$, pulling back defines a
$\VV\GG_s$-functor
\[u\upst:\ol{M}(T)\longrightarrow \ol{M}(S).\]

\begin{them} \label{mbfc}
For every $k$-scheme $T$ the functor $\ol{M}(T)\rightarrow\VV\GG_s$ is a
cofibration, whose fibers are groupoids. In
other words, $\ol{M}(T)$ is cofibered in groupoids over $\VV\GG_s$.

For every base change $u:S\rightarrow T$ the $\VV\GG_s$-functor
$u\upst:\ol{M}(T)\rightarrow \ol{M}(S)$ is cocartesian.
\end{them}
\begin{pf}
To prove that $\ol{M}(T)\rightarrow\VV\GG_s$ is a cofibration, we need to prove
the following. Let
$(\xi,a,\tau',\phi):(V,\tau)\rightarrow(W,\sigma)$ be a morphism in $\VV\GG_s$
and $(C,x,f)$ a stable $(V,\tau)$-map
over $T$. Then there exists a {\em pushforward }$(D,y,h)$ of $(C,x,f)$ under
$(\xi,a,\tau',\phi)$. This pushforward
comes with a morphism $p:(C,x,f)\rightarrow(D,y,h)$ of stable maps covering
$(\xi,a,\tau',\phi)$ and is characterized by
the following universal mapping property. Whenever
$(\eta,b,\sigma',\psi):(W,\sigma)\rightarrow(U,\rho)$ is another
morphism in $\VV\GG_s$, $(E,z,e)$ a stable $(U,\rho)$-map over $T$ and
$r:(C,x,f)\rightarrow(E,z,e)$ a morphism of
stable maps covering $(\eta\xi,ac,\tau'',\psi\chi):(V,\tau)\rightarrow(U,\rho)$
(in the notation of
Definition~\ref{dmmsg}), there exists a unique morphism of stable maps
$q:(D,y,h)\rightarrow(E,z,e)$ covering
$(\eta,b,\sigma',\psi):(W,\sigma)\rightarrow(U,\rho)$ such that $r=q\comp p$.
\begin{equation} \label{umprq}
\begin{array}{ccccc}
 & & \stackrel{r}{\overtoparrow} & & \\
(C,x,f) & \stackrel{p}{\longrightarrow} & (D,y,h) &
\stackrel{q}{\longrightarrow} & (E,z,e) \\
\mid & & \mid & & \mid \\
(V,\tau) & \stackrel{(\xi,a,\tau',\phi)}{\longrightarrow} & (W,\sigma) &
\stackrel{(\eta,b,\sigma',\psi)}{\longrightarrow} & (U,\rho) \\
 & & {\underbottomarrow\atop(\eta\xi,ac,\tau'',\psi\chi)} &  &
\end{array}
\end{equation}
To prove that $u\upst:\ol{M}(T)\rightarrow\ol{M}(S)$ is always cocartesian, we
need to prove that this pushforward
commutes with base change.

Recall that we also wish to prove Proposition~\ref{cmsmm}, i.e.\ that if
morphisms of stable maps
$p:(C,x,f)\rightarrow(D,y,h)$ and $q:(D,y,h)\rightarrow(E,z,e)$ are given as in
(\ref{umprq}), then the composition
$r:(C,x,f)\rightarrow(E,z,e)$ is also a morphism of stable maps.

Purely formal considerations tell us that to prove these three facts, we may
decompose the morphism
$(\xi,a,\tau',\phi):(V,\tau)\rightarrow(W,\sigma)$ into a composition of other
morphisms in any way we wish and prove the
three facts for the factors of this decomposition. We shall thus consider the
following five cases.

{\em Case I (Changing $V$). } In this case $\sigma=\xi\lst\tau$. This means
that $\sigma$ is the pushforward of $\tau$
under $\xi:V\rightarrow W$, using the fact that $\VV\GG_s\rightarrow\VV$ is a
cofibration (Remark~\ref{ugscf}). In other
words, $\sigma$ is the stabilization of $\tau$ with respect to the induced
$H_2(W)^+$-structure (Proposition~\ref{ppes}).
Thus $\tau'=\sigma$ and $\phi=\id\s$.

In all other cases $W=V$ and $\xi$ is the identity. In the next two cases
$a=\id\t$ and $\tau'=\tau$.

{\em Case II (Contracting and edge). } The contraction
$\phi:\tau\rightarrow\sigma$ contracts exactly one edge
$\{i_1,i_2\}\subset F\t$ and we have $v_1\not=v_2$, where $v_1=\del(i_1)$ and
$v_2=\del(i_2)$. To fix notation, let
$v_0=\phi(v_1)=\phi(v_2)$.

{\em Case III (Contracting a loop). } This is the same as Case~II, except that
we have $v_1=v_2$.

In the last two cases $\tau'=\sigma$ and $\phi=\id\s$.

{\em Case IV (Complete combinatorial). } The combinatorial morphism
$a:\sigma\rightarrow\tau$ has the property that
$a:F\s(v)\rightarrow F\t(a(v))$ is a bijection, for all $v\in V\s$.

{\em Case V (Removing a tail). } In this case, $V\s=V\t$, we have given a
vertex $v_0\in V\t$ and a {\em tail }$i_0\in
F\t(v_0)$ of $\tau$ and we have
\begin{enumerate}
\item $F\s= F\t-\{i_0\}$,
\item $\del\s=\del\t\resto F\s$,
\item $j\s=j\t\resto F\s$.
\end{enumerate}
\renewcommand{\qed}{}\end{pf}

Note that the proof of Proposition~\ref{cmsmm} is only interesting (if at all)
for Case~II, since only in this case
carrying out the composition of $(\xi,a,\tau',\phi)$ and
$(\eta,b,\sigma',\psi)$ involves the second case of the
construction of stable pullback (Section~\ref{graphs}).

{\em Case I\@. } First we note the following trivial lemma.

\begin{lem} \label{tcvw}
Assume that $\tau$ is stable with respect to the induced $H_2(W)^+$-structure,
so that $\sigma=\tau$ and $a=\id\t$. Then
if $(C,x,\xi\comp f)$ satisfies the stability condition it may serve as
pushforward of $(C,x,f)$ under $\xi$. \qed
\end{lem}

We shall now reduce Case~I to Cases~IV and~V\@. By the claimed compatibility
with base change, we may construct the
pushforward locally, and pass to an \'etale cover of $T$, whenever desirable.
Thus we add tails to $\tau$, obtaining
$\tilde{\tau}$, and corresponding sections of $C$, obtaining $(C,\tilde{x})$
until $\tilde{\tau}$ with the induced
$H_2(W)^+$-structure is stable and $(C,\tilde{x},\xi\comp f)$ satisfies the
stability condition. Then we have the
commutative diagram
\begin{equation} \label{divgs}
\comdia{(V,\tilde{\tau})}{}{(V,\tau)}{}{}{}{(W,\tilde{\tau})}{}{(W,\sigma)}
\end{equation}
in $\VV\GG_s$. The top row of (\ref{divgs}) is covered by
$(C,\tilde{x},f)\rightarrow(C,x,f)$, and clearly $(C,x,f)$ is
the pushforward of $(C,\tilde{x},f)$ under
$(V,\tilde{\tau})\rightarrow(V,\tau)$ (see also Case~V). The first column of
(\ref{divgs}) is covered by $(C,\tilde{x},f)\rightarrow(C,\tilde{x},\xi\comp
f)$, which is a pushforward by
Lemma~\ref{tcvw}. Now the pushforward of $(C,\tilde{x},\xi\comp f)$ under
$(W,\tilde{\tau})\rightarrow(W,\sigma)$ will
also be the sought after pushforward of $(C,x,f)$ under
$(V,\tau)\rightarrow(W,\sigma)$. But
$(W,\tilde{\tau})\rightarrow(W,\sigma)$ is covered by Cases~IV and~V, achieving
the reduction. \qed

{\em Case II\@. } The diagram defining the composition of $\phi$ and
$(b,\sigma',\psi)$ is
\[\begin{array}{ccccc}
\tau' & \stackrel{\chi}{\longrightarrow} & \sigma' &
\stackrel{\psi}{\longrightarrow} & \rho \\
\ldiag{c} & & \rdiag{b} & & \\
\tau & \stackrel{\phi}{\longrightarrow} & \sigma. & &
\end{array}\]
Let us first deal with the proof of Proposition~\ref{cmsmm}.

\begin{lem} \label{lc2}
For every $i\in F_{\sigma'}$ we have
\[q_{\del(i)}p_{\del
c\chi^F(i)}(x_{c\chi^F(i)})=q_{\del(i)}p_{\del\phi^Fb(i)}(x_{\phi^Fb(i)}).\]
\end{lem}
\begin{pf}
Assume that $c\chi^F(i)\not=\phi^Fb(i)$, since otherwise there is nothing to
prove. In this case, necessarily,
$c\chi^F(i)$ is being contracted by $\phi$. Without loss of generality, let
$c\chi^F(i)=i_1$, so the situation is as in
the following diagram (cf.\ ~(\ref{tadsp})).
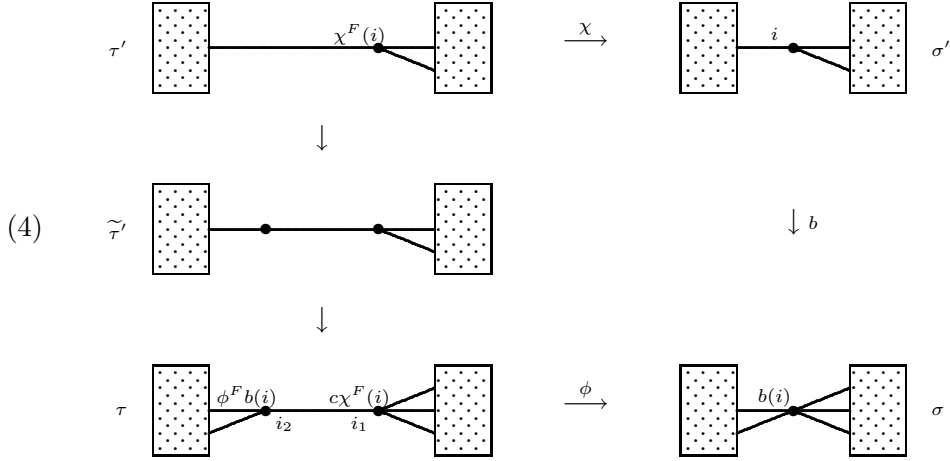
\begin{equation} \label{fdgc}
\begin{array}{ccc}
\beginpicture
\setcoordinatesystem units <.5cm,.3cm> point at 4.5 2
\setplotarea x from 0 to 9, y from 0 to 4
\plot 1.5 2 7.5 2 /
\plot 6 2 7.5 1 /
\shaderectangleson
\setshadegrid span <1mm>
\putrectangle corners at 0 0 and 1.5 4
\putrectangle corners at 7.5 0 and 9 4
\put {\circle*{4}} [Bl] at 6 2
\put {$\scriptstyle\chi^F(i)$} at 5.5 2.6
\axis left invisible label {$\scriptstyle\tau'$} /
\axis right invisible label {\phantom{$\scriptstyle\tau'$}} /
\axis bottom invisible label {\phantom{.}} /
\endpicture & \stackrel{\chi}{\longrightarrow} &
\beginpicture
\setcoordinatesystem units <.5cm,.3cm> point at 3 2
\setplotarea x from 0 to 6, y from 0 to 4
\plot 1.5 2 4.5 2 /
\plot 3 2 4.5 1 /
\shaderectangleson
\setshadegrid span <1mm>
\putrectangle corners at 0 0 and 1.5 4
\putrectangle corners at 4.5 0 and 6 4
\put {\circle*{4}} [Bl] at 3 2
\put {$\scriptstyle i$} at 2.5 2.6
\axis left invisible label {\phantom{$\scriptstyle\sigma'$}} /
\axis right invisible label {$\scriptstyle\sigma'$} /
\axis bottom invisible label {\phantom{.}} /
\endpicture \\
\ldiag{} & & \\
\beginpicture
\setcoordinatesystem units <.5cm,.3cm> point at 4.5 2
\setplotarea x from 0 to 9, y from 0 to 4
\plot 1.5 2 7.5 2 /
\plot 6 2 7.5 1 /
\shaderectangleson
\setshadegrid span <1mm>
\putrectangle corners at 0 0 and 1.5 4
\putrectangle corners at 7.5 0 and 9 4
\put {\circle*{4}} [Bl] at 3 2
\put {\circle*{4}} [Bl] at 6 2
\axis left invisible label {$\scriptstyle\tilde{\tau}'$} /
\axis right invisible label {\phantom{$\scriptstyle\tilde{\tau}'$}} /
\axis top invisible label {\phantom{.}} /
\axis bottom invisible label {\phantom{.}} /
\endpicture & & \rdiag{b} \\
\ldiag{} & & \\
\beginpicture
\setcoordinatesystem units <.5cm,.3cm> point at 4.5 2
\setplotarea x from 0 to 9, y from 0 to 4
\plot 1.5 2 7.5 2 /
\plot 3 2 1.5 1 /
\plot 7.5 3 6 2 7.5 1 /
\shaderectangleson
\setshadegrid span <1mm>
\putrectangle corners at 0 0 and 1.5 4
\putrectangle corners at 7.5 0 and 9 4
\put {\circle*{4}} [Bl] at 3 2
\put {\circle*{4}} [Bl] at 6 2
\put {$\scriptstyle\phi^Fb(i)$} at 2.5 2.6
\put {$\scriptstyle c\chi^F(i)$} at 5.5 2.6
\put {$\scriptstyle i_1$} at 5.5 1.4
\put {$\scriptstyle i_2$} at 3.5 1.4
\axis left invisible label {$\scriptstyle\tau$} /
\axis right invisible label {\phantom{$\scriptstyle\tau$}} /
\axis top invisible label {\phantom{.}} /
\endpicture & \stackrel{\phi}{\longrightarrow} &
\beginpicture
\setcoordinatesystem units <.5cm,.3cm> point at 3 2
\setplotarea x from 0 to 6, y from 0 to 4
\plot 1.5 2 4.5 2 /
\plot 3 2 1.5 1 /
\plot 4.5 3 3 2 4.5 1 /
\shaderectangleson
\setshadegrid span <1mm>
\putrectangle corners at 0 0 and 1.5 4
\putrectangle corners at 4.5 0 and 6 4
\put {\circle*{4}} [Bl] at 3 2
\put {$\scriptstyle b(i)$} at 2.5 2.6
\axis left invisible label {\phantom{$\scriptstyle\sigma$}} /
\axis right invisible label {$\scriptstyle\sigma$} /
\axis top invisible label {\phantom{.}} /
\endpicture
\end{array}
\end{equation}
Here, ${\tau}'$ is the stable pullback and $\tilde{\tau}'$ the intermediate
graph used in the construction of
$\tau'$. Using the fact that $p$ is a morphism of stable maps we get a morphism
$p_{12}:C_{12}\rightarrow
D_{v_0}$ of prestable curves, where
\[C_{12}=C_{v_1}\amalg_{x_{i_1},x_{i_2}}C_{v_2}.\]
Compose this with $q_{\del(i)}:D_{v_0}\rightarrow E_{\psi\del(i)}$. Let
$f_{12}:C_{12}\rightarrow V$ be the map induced
from $f_{v_1}$ and $f_{v_2}$ and $\tilde{x}=x\resto F\t(v_1)\cup
F\t(v_2)-\{i_1,i_2\}$. Then $(C_{12},\tilde{x},f_{12})$
is a stable map and
\[q_{\del(i)}\comp
p_{12}:(C_{12},\tilde{x},f_{12})\longrightarrow(E_{\psi\del(i)},z\resto
F_\rho(\psi\del(i)),e_{\psi\del(i)})\] is a morphism of stable maps to which
Lemma~\ref{lcsn} applies, with $C'=C_{v_2}$
and $x_{\phi^Fb(i)}\in C'$ being the only marked point coming from
$F_\rho(\psi\del(i))$, if there exists such a point
at all (this is because $\tau'\not=\tilde{\tau}'$). So by Lemma~\ref{lcsn}
$q_{\del(i)}p_{v_2}(C_{v_2})$ is a point in
$E_{\psi\del(i)}$. To be precise, this holds if $T$ is the spectrum of an
algebraically closed field. For the general
case, applying Lemma~\ref{tloq} yields that $q_{\del(i)}\comp p_{v_2}$ factors
thought $T$. In particular,
\[q_{\del(i)}p_{v_2}(x_{\phi^Fb(i)})= q_{\del(i)}p_{v_2}(x_{i_2})=
q_{\del(i)}p_{v_1}(x_{i_1}),\]
which is what we set out to prove.
\end{pf}

Let us check that $r:(C,x,f)\rightarrow(E,z,e)$ is a morphism of stable maps,
i.e.\ satisfies Properties~(1) through~(4)
from Definition~\ref{dsm}.

{\em Property (1). }  Let $i\in F_\rho$. The we have
\begin{align*}
r_{\del\chi^F\psi^F(i)}(x_{c\chi^F\psi^F(i)}) &=  q_{\del\psi^F(i)}p_{\del
c\chi^F\psi^F(i)}(x_{c\chi^F\psi^F(i)}) \\
\intertext{by Definition~\ref{dsm},}
 &=  q_{\del\psi^F(i)}p_{\del\phi^Fb\psi^F(i)}(x_{\phi^Fb\psi^F(i)}) \\
\intertext{by Lemma~\ref{lc2},}
 &=  q_{\del\psi^F(i)}(y_{b\psi^F(i)}) \\
 &=  z_i,
\end{align*}
since $p$ and $q$ are morphisms of stable maps.

{\em Property (2). }  Let $\{j_1,j_2\}$ be an edge of $\tau'$ which is being
contracted by $\psi\chi$. Let
$u_1=\del j_1$ and $u_2=\del j_2$.

{\em Case 1. } Let $\{j_1,j_2\}$ be contracted by $\chi$. Then
$\{c(j_1),c(j_2)\}$ is being contracted by $\phi$. So
without loss of generality $c(j_1)=i_1$ and $c(j_2)=i_2$. Then
\begin{align*}
r_{u_1}(x_{i_1}) & =  q_{\chi(u_1)}p_{v_1}(x_{i_1}) \\
 & =  q_{\chi(u_2)}p_{v_2}(x_{i_2}) \\
 & =  r_{u_2}(x_{i_2}),
\end{align*}
since $p$ is a morphism of stable maps and $\chi(u_1)=\chi(u_2)$.

{\em Case 2. } If $\{j_1,j_2\}$ is not contracted by $\chi$, then there exists
a unique edge $\{j_1',j_2'\}$ of
$\sigma'$ being contracted by $\psi$, such that $j_1=\chi^F(j_1')$ and
$j_2=\chi^F(j_2')$. Then
\begin{align*}
r_{u_1}(x_{c(j_1)}) & =  q_{\chi(u_1)}p_{c(u_1)}(x_{c(j_1)}) \\
 & =  q_{\chi(u_1)}p_{\del\phi^Fb(j_1')}(x_{\phi^Fb(j_1')}) \\
\intertext{by Lemma~\ref{lc2},}
 & =  q_{\chi(u_1)}(y_{b(j_1')}) \\
 & =  q_{\chi(u_2)}(y_{b(j_2')}) \\
\intertext{since $q$ is a morphism of stable maps,}
 & =  r_{u_2}(x_{c(j_2)}),
\end{align*}
by symmetry.

{\em Property (3). } This follows from the fact that the composition of
morphisms of prestable curves is again
a morphism of prestable curves.

{\em Property (4). } Straightforward.

This finishes the proof of Proposition~\ref{cmsmm} in Case~II\@. Let us now
construct the pushforward $(D,y,h)$ of
$(C,x,f)$ under $\phi$.

Let $w\in V\s$. If $w\not=v_0$, let $v$ be the unique vertex $v\in V\t$ such
that $\phi_V(v)=w$ and set $D_w=C_v$. If
$w=v_0$ set
\[D_{v_0}=C_{v_1}\amalg_{x_{i_1},x_{i_2}}C_{v_2}.\]
This defines a family of prestable curves $D$. For every $v\in V\t$ let
$p_v:C_v\rightarrow D_{\phi(v)}$ be the canonical
map. Define sections $y_i$, for $i\in F\s$, by
\[y_i=p_{\del\phi^F(i)}(x_{\phi^F(i)}).\]
Finally, define for every $w\in V\s$ a map $g_w:D_w\rightarrow V$ from $f$ (by
using Proposition~\ref{glue}, if
$w=v_0$). Essentially by definition, $(D,y,h)$ is a stable $(V,\sigma)$-map and
$p:(C,x,f)\rightarrow(D,y,h)$ is a
morphism of stable maps covering $\phi:(V,\tau)\rightarrow(V,\sigma)$. It
remains to check that $(D,y,h)$ satisfies the
universal mapping property of a pushforward under $\phi$. So let
$r:(C,x,f)\rightarrow(E,z,e)$ as in
Diagram~(\ref{umprq}) be given.

Let $u\in V_{\sigma'}$. We need to define a unique morphism
$q_u:D_{b(u)}\rightarrow E_{\psi(u)}$ such that for every
$u'\in V_{\tau'}$, satisfying $\chi(u')=u$, the diagram
\[\begin{array}{ccc}
C_{c(u')} & & \\
\ldiag{p_{c(u')}} & \sediagr{r_{u'}} & \\
D_{b(u)} & \stackrel{q_u}{\longrightarrow} & E_{\psi(u)}
\end{array}\]
commutes. If $b(u)\not=v_0$, necessarily, $q_u=r_{u'}$. So let $b(u)=v_0$. If
there are two vertices $u_1'$ and $u_2'$
such that $\chi(u_1')=\chi(u_2')=u$, then we have two maps
$r_{u_1'}:C_{v_1}\rightarrow E_{\psi(u)}$ and
$r_{u_2'}:C_{v_2}\rightarrow E_{\psi(u)}$ giving rise to a unique map
$q_u:D_{v_0}\rightarrow E_{\psi(u)}$. If there is
only one vertex $u_1'$ of $\tau'$ such that $\chi(u_1')=u$, then we are in a
situation as in Diagram~(\ref{fdgc}), and
by Lemma~\ref{lcsn} $q_u$ has to map $C_{v_2}\subset D_{v_0}$ to a single point
of $E_{\psi(u)}$ and
$r_{u_1'}:C_{v_1}\rightarrow E_{\psi(u)}$ suffices to determine
$q_u:D_{v_0}\rightarrow E_{\psi(u)}$ uniquely. This
defines $q:D\rightarrow E$ satisfying all properties required of a morphism of
stable maps, as some routine
considerations reveal. This finishes Case~II\@. \qed

{\em Case III\@. } This case is similar to Case~II, but much simpler, because
the construction of the composition of
$\phi$ and $(b,\sigma,\psi)$ is simpler, and thus for every $i\in F_{\sigma'}$
we have $c\chi^F(i)=\phi^Fb(i)$. We use
Proposition~\ref{glue1} instead of Proposition~\ref{glue} to construct the
pushforward of $(C,x,f)$ under $\phi$, gluing
the two sections $x_{i_1}$ and $x_{i_2}$ of $C_{v_1}=C_{v_2}$ to obtain
$D_{v_0}$. \qed

{\em Case IV\@. } To construct the pushforward, set $D_v=C_{a(v)}$,
$p_v:C_{a(v)}\rightarrow D_v$ the identity and
$h_v=f_{a(v)}$, for every $v\in V\s$. Moreover, for $i\in F\s$ set
$y_i=x_{a(i)}$. To check that $(D,y,h)$ is a stable
map and $p:(C,x,f)\rightarrow(D,y,h)$ a morphism of stable maps, the only fact
to check is that for every $i\in F\s$ we
have $h_{\del(i)}(y_i)=h_{\del(j(i))}(y_{j(i)})$, in other words
\[f_{\del a(i)}(x_{a(i)})=f_{\del(j(i))}(x_{a(j(i))}).\]
Here, Condition~(\ref{commor3}) in the definition of combinatorial morphism of
$A$-graphs (Definition~\ref{commor})
enters in. It implies this claim together with Corollary~\ref{zc}. The
universal mapping property of $(D,y,h)$ is easily
verified. \qed

{\em Case V\@. }  Before we can treat this case, we need a few preparations.

\begin{prop} \label{mplsm}
Let $(C,x_1,\ldots,x_n,f)$ be a stable map over a field from a curve of genus
$g$ to $V$ and $M$ an ample invertible
sheaf on $V$. Then
\[L=\omega_C(x_1+\ldots+x_n)\otimes f\upst M^{\otimes3}\]
is ample on $C$. Here $\omega_C$ is the dualizing sheaf of $C$.
\end{prop}
\begin{pf}
Let us first consider the case that $C$ has no nodes, so that $C$ is
irreducible and non-singular. Then is suffices to
prove that $\deg L>0$.

{\em Case 1. }The image $f(C)$ is a point. Then $\deg f\upst M=0$ and we have
\[\deg L=\deg\omega_C+n=2g-2+n\geq1,\]
by the stability condition.

{\em Case 2. }The image $f(C)$ is not a point. Then $\deg f\upst M\geq1$ and so
\[\deg L=2g-2+n+3\deg f\upst M\geq 2g-2+n+3=2g+n+1>0.\]

So suppose now that $C$ has a node $P$. Let $q:C'\rightarrow C$ be the curve
obtained by blowing up $P$ and let
$P_1,P_2\in C'$ be the two points lying over $P$. Let $L'=q\upst L$ and
$f'=f\comp q$.

{\em Case 1. }The curve $C'$ is connected. Then
$(C',x_1,\ldots,x_n,P_1,P_2,f')$ is a stable map and
\[L'=\omega_{C'}(x_1+\ldots+x_n+P_1+P_2)\otimes {f'}\upst M^{\otimes3}.\]

{\em Case 2. }The curve $C'$ is disconnected. Let $C_1'$ and $C_2'$ be the two
components of $C'$ and $L_1',L_2'$ the
restriction of $L'$ to $C_1'$ and $C_2'$, respectively. Let
$f_i':C_i'\rightarrow V$ for $i=1,2$ be the map induced by
$f'$. Without loss of generality assume that $x_1,\ldots,x_r\in  C_1'$ and
$x_{r+1},\ldots,x_n\in C_2'$, for some $0\leq
r\leq n$ and $P_1\in C_1'$, $P_2\in C_2'$. Then
$(C_1',x_1,\ldots,x_r,P_1,f_1')$ and $(C_2',x_{r+1},\ldots,x_n,P_2,f_2')$
are stable maps and
\begin{eqnarray*}
L_1' & = & \omega_{C_1'}(x_1+\ldots +x_r+P_1)\otimes {f_1'}\upst M^{\otimes3}\\
L_2' & = & \omega_{C_2'}(x_{r+1}+\ldots+x_n+P_2)\otimes {f_2'}\upst
M^{\otimes3}.
\end{eqnarray*}

Thus by induction on the the number of nodes we may assume that $L'$ is ample
on $C'$. Let $\fF$ be a coherent sheaf on
$C$ and $\fF'=q\upst\fF$. Then there exists an $n_0$ such that for all $n\geq
n_0$ we have that $\fF'\otimes{L'}^{\otimes
n}(-P_1)$, $\fF'\otimes{L'}^{\otimes n}(-P_2)$ and $\fF'\otimes{L'}^{\otimes
n}(-P_1-P_2)$ are generated by global
sections. This implies that $\fF\otimes L^{\otimes n}$ is generated by global
sections. So $L$ is ample.
\end{pf}

We will now consider the following setup. Let $(C,x_1,\ldots,x_{n+1},f)$ be a
stable map over $T$ from an
$(n+1)$-pointed curve of genus $g$ to $V$ of class $\beta\in H_2(V)^+$, where
$2g+n\geq3$ if $\beta=0$ (otherwise,
$n\geq0$).

If $t$ is a geometric point of $T$ and $C'$ a component of $C_t$, then we say
that $C'$ {\em is to be contracted}, if,
after removing $x_{n+1}$, the normalization of $C'$ violates the stability
condition. Equivalently,
\begin{enumerate}
\item $C'$ is rational without self intersection (so that $C'$ is equal to its
normalization),
\item $x_{n+1}\in C'$,
\item $C'$ has exactly two special points besides $x_{n+1}$, at least one of
which is not a marked point, but a node,
\item $f_t(C')$ is a single point of $V$.
\end{enumerate}
Pictorially, the only two possible components to be contracted look as follows.
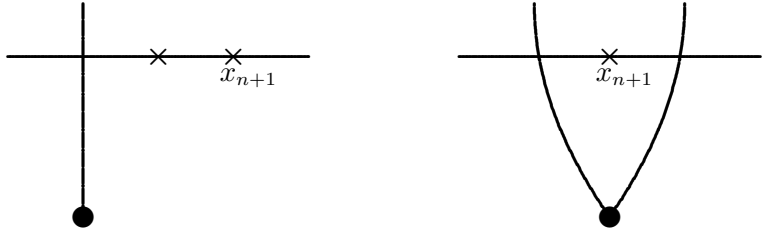
\[\beginpicture

\setcoordinatesystem units <1cm,.71cm> point at 4 0
\setplotarea x from -1 to 3, y from -3 to 1
\plot 0 1 0 -3 /
\plot -1 0 3 0 /
\put {\circle*{8}} [Bl] at 0 -3
\put {$\boldsymbol{\times}$} at 1 0
\put {$\boldsymbol{\times}$} at 2 0
\put {$x_{n+1}$} at 2.2 -.4

\setcoordinatesystem units <1cm,.71cm> point at -3 0
\setplotarea x from -2 to 2, y from -3 to 1
\plot -2 0 2 0 /
\setquadratic
\plot 0 -3 -.75 -1 -1 1 /
\plot 0 -3  .75 -1  1 1 /
\put {\circle*{8}} [Bl] at 0 -3
\put {$\boldsymbol{\times}$} at 0 0
\put {$x_{n+1}$} at 0.2 -.4

\endpicture\]
Note that every geometric fiber of $\pi:C\rightarrow T$ has at most one
component to be contracted.

We say a $T$-morphism $q:C\rightarrow U$, for a $T$-scheme $U$, {\em contracts
the components to be contracted}, if for
every geometric point $t$ of $T$ the map (of underlying Zariski topological
spaces) $q_t:C_t\rightarrow U_t$ maps every
component to be contracted to a single point in $U_t$. For example,
$f:C\rightarrow V_T$ contracts the components to be
contracted.

\begin{prop} \label{cae}
There exists a universal morphism $p:C\rightarrow\tilde{C}$ contracting the
components to be contracted. Let
$\tilde{f}:\tilde{C}\rightarrow V$ be the unique map given by the universal
mapping property of $(\tilde{C},p)$. Then
$(\tilde{C},p(x_1),\ldots, p(x_n),\tilde{f})$ is a stable map from an
$n$-pointed curve of genus $g$ to $V$ of class
$\beta$.
\end{prop}
\begin{pf}
This is a variation on Section~1 of \cite{knudsen}.  Let us first prove the
proposition for the case that
$T$ is the spectrum of an algebraically closed field. Let $C'$ be a component
to be contracted.

{\em Case 1. }The component $C'$ has one node. We define $\tilde{C}=C-(C'-C)$
and let $p:C\rightarrow \tilde{C}$ be the
obvious map. Clearly, $\O_{\tilde{C}}=p\lst\O_C$, so $\tilde{C}$ satisfies the
universal mapping property by
Lemma~\ref{tloq}. The rest is trivial.

{\em Case 2. }The component $C'$ has two nodes. We define $\ol{C}=C-(C'-C)$ and
let $P_1$ and $P_2$ be the two points of
$\ol{C}$ that intersect $C'$. Then we set $\tilde{C}=\ol{C}/P_1\sim P_2$, i.e.\
we identify the two points $P_1$ and
$P_2$. We then proceed similarly as in Case~1.

\begin{lem} \label{clocct}
Let $T$ be the spectrum of an algebraically closed field and let $\tilde{C}$ be
the universal curve contracting the
components of $C$ to be contracted.  Choose an ample invertible sheaf $M$ on
$V$. Let
\[L=\omega_C(x_1+\ldots+x_n)\otimes f\upst M^{\otimes3}\]
and
\[\tilde{L}=\omega_{\tilde{C}}(p(x_1)+\ldots+p(x_n))\otimes \tilde{f}\upst
M^{\otimes3}.\]
Then for all $k\geq0$ we have
\begin{enumerate}
\item $\tilde{L}^{\otimes k}=p\lst L^{\otimes k}$,
\item $p\upst \tilde{L}^{\otimes k}={L}^{\otimes k}$,
\item $R^1f\lst  L^{\otimes k}=0$,
\item $H^i(\tilde{C},\tilde{L}^{\otimes k})=H^i(C,L^{\otimes k})$, for $i=0,1$.
\end{enumerate}
\end{lem}
\begin{pf}
This is analogous to Lemma~1.6 of \cite{knudsen}.
\end{pf}

\begin{lem} \label{tfll}
Let $T$ be the spectrum of an algebraically closed field and let $L$ be defined
as in Lemma~\ref{clocct}. Define the
open subset $U$ of $C$ by
\[U=\{x\in C\st \rtext{$x$ is smooth and $x$ is not in any component to be
contracted}\}.\]
Then for $k$ sufficiently large we have
\begin{enumerate}
\item $L^{\otimes k}$ is generated by global sections, \label{genglob}
\item $H^1(C,L^{\otimes k})=0$,
\item $L^{\otimes k}$ is normally generated, \label{normgen}
\item $L^{\otimes k}(-P)$ is generated by global sections for all $P\in U$.
\label{mpglob}
\end{enumerate}
(The sheaf $L$ is normally generated if $\Gamma(C,L)^{\otimes
k}\rightarrow\Gamma(C,L^{\otimes k})$ is surjective, for
all $k\geq1$.)
\end{lem}
\begin{pf}
Let $\tilde{C}$ and $\tilde{L}$ be as in Lemma~\ref{clocct}. Note that one can
apply Proposition~\ref{mplsm} to
$\tilde{C}$ and $\tilde{L}$. Then the results are implied by
Lemma~\ref{clocct}.
\end{pf}

We can now proceed with the proof of Proposition~\ref{cae} for general base
$T$. Choose an ample invertible sheaf $M$ on
$V$ and consider on $C$ the invertible sheaf
\[L=\omega_{C/T}(x_1+\ldots+x_n)\otimes f\upst M^{\otimes3},\]
where $\omega_{C/T}$ is the relative dualizing sheaf of $C$ over $T$.  Then
form
\[\sS=\bigoplus_{k\geq0}\pi\lst L^{\otimes k},\]
where $\pi:C\rightarrow T$ is the structure map, and let
\[\tilde{C}=\proj\sS.\]

{\em Claim 1. }The formation of $\tilde{C}$ commutes with base change.
\begin{pf}
Clearly, the formation of $L^{\otimes k}$ commutes with base change.  That the
formation of $\pi\lst L^{\otimes k}$
commutes with base change for $k$ sufficiently large follows from the fact that
$H^1(C_t,L^{\otimes k})=0$, for all
$t\in T$, by Lemma~\ref{tfll}. For $k=0$ this is trivially true. Thus the
formation of
\[\sS^{(d)}=\bigoplus_{d\mid k}\sS_k\]
commutes with base change, for a suitable $d>0$. This implies the claim for
$\tilde{C}$, since
\[\tilde{C}=\proj\sS=\proj\sS^{(d)}.\qed\]
\renewcommand{\qed}{}\end{pf}

{\em Claim 2. }The structure map $\tilde{\pi}:\tilde{C}\rightarrow T$ is flat
and projective.
\begin{pf}
The flatness of $\proj\sS^{(d)}$ follows from the fact that $\pi\lst L^{\otimes
k}$ is locally free, for $d\mid k$,
which follows from the fact that its formation commutes with base change. By
passing to a larger $d$ if necessary, we
may assume that for every $k\geq0$ the homomorphism
\[\pi\lst(L^{\otimes d})^{\otimes k}\longrightarrow\pi\lst(L^{\otimes dk}) \]
is surjective. This may be checked on fibers and thus follows from
Lemma~\ref{tfll}(\ref{normgen}). So $\sS^{(d)}$ is
generated by $\sS^{(d)}_1$ and hence $\proj\sS^{(d)}$ is projective by
Proposition~5.5.1 in \cite{ega2}.
\end{pf}

{\em Claim 3. }The canonical morphism from an open subset of $C$ to $\tilde{C}$
is everywhere defined, proper and
surjective.
\begin{pf}
This canonical morphism is defined by $\pi\upst\sS\rightarrow\bigoplus_k
L^{\otimes k}$, or equivalently by
$\sS\rightarrow\bigoplus_k\pi\lst L^{\otimes k}$ (see Section~3.7 in
\cite{ega2}). For it to be everywhere defined, it
suffices to prove that $\pi\upst\pi\lst L^{\otimes k}\rightarrow L^{\otimes k}$
is an epimorphism, for $k$ sufficiently
large. This may be checked on fibers and thus follows from
Lemma~\ref{tfll}(\ref{genglob}). That the canonical morphism
is dominant follows immediately, since $\sS\rightarrow\bigoplus\pi\lst
L^{\otimes k}$ is an isomorphism. That it is
proper, is clear. So it has to be surjective.
\end{pf}

We call this canonical morphism $p:C\rightarrow\tilde{C}$.

{\em Claim 4. }Let $x$ be a geometric point of $\tilde{C}$ and $p^{-1}(x)$ the
fiber of $p$ over $x$. Then either the
cardinality of $p^{-1}(x)$ is one or $p^{-1}(x)$ is a component of
$C_{\tilde{\pi}(x)}$ to be contracted.
\begin{pf}
Without loss of generality we may assume that $T$ is the spectrum of an
algebraically closed field. Then with the
notation of Lemma~\ref{tfll} and by Property~(\ref{mpglob}) of the same lemma,
we have that $p\resto
U:U\rightarrow\tilde{C}$ is an open immersion. If $C'$ is to be contracted,
then $L\resto C'\cong\O_{C'}$, and so $C'$
is mapped to a point in $\tilde{C}$. These facts clearly imply Claim~4.
\end{pf}

{\em Claim 5. }We have $p\lst\O_C=\O_{\tilde{C}}$.
\begin{pf}
With the notation of Claim~4, note that
\[H^1(p^{-1}(x),\O_C\otimes_{\O_{\tilde{C}}}\kappa(x))=0,\]
since $p^{-1}(x)$ is rational if it is one and not zero dimensional. So by
Corollary~1.5 in \cite{knudsen}, the
formation of $p\lst\O_C$ commutes with base change in $T$. So to prove the
claim, we may assume that $T$ is the
spectrum of an algebraically closed field, but then it is clear.
\end{pf}

Now by Lemma~\ref{tloq} the last three claims imply that
$p:C\rightarrow\tilde{C}$ is a universal morphism contracting
the components to be contracted. In particular, we get a unique morphism
$\tilde{f}:C\rightarrow V$ such that
$\tilde{f}\comp p=f$. The fact that
$(\tilde{C},p(x_1),\ldots,p(x_n),\tilde{f})$ is a stable map from an
$n$-pointed
curve of genus $g$ to $V$ of class $\beta$ may now be checked on fibers, which
has already been done.  This finishes the
proof of Proposition~\ref{cae}.
\end{pf}

We now proceed with the proof of Theorem~\ref{mbfc} in Case~V\@. Let $n=\#
F\s(v_0)$. Choose an identification
$\ul{n+1}\rightarrow F\t(v_0)$ mapping $n+1$ to $i_0$, the flag being removed.
Then
$(C_{v_0},x_1,\ldots,x_{n+1},f_{v_0})$ is a stable map to which
Proposition~\ref{cae} applies and we let
$p_{v_0}:C_{v_0}\rightarrow D_{v_0}$ be the universal morphism contracting the
components to be contracted. For
$v\not=v_0$ we let $D_v=C_v$ and $p_v:C_v\rightarrow D_v$ be the identity. It
is then clear how to define $y$ and $h$ to
get a stable map $(D,y,h)$ satisfying the universal mapping property of a
pushforward under the graph morphism
$\tau\rightarrow \sigma$ given by $a:\sigma\rightarrow \tau$. This completes
the proof of Case~V\@. \qed

To complete the proof of Theorem~\ref{mbfc} we need to show that if
$(V,\tau,\alpha)$ is an object of $\VV\GG_s$ and
$p:(C,x,f)\rightarrow(D,y,h)$ is a morphism of stable $(V,\tau,\alpha)$-maps
(covering the identity of
$(V,\tau,\alpha)$), then $p$ is an isomorphism.

This is immediately reduced to the case that $(V,\tau,\alpha)=(V,g,n,\alpha)$
and using Lemma~\ref{tloq} to the case that
$T$ is the spectrum of an algebraically closed field. Then it follows from the
stability condition that $p$ cannot
contract any rational components, so it is injective. To prove that $p$ is
surjective use induction on the number of
nodes of $D$. So let $P$ be a node of $D$ and let $D'$ be the curve obtained
from $D$ by blowing up $P$ and let
$p':C'\rightarrow D'$ by the pullback of $p:C\rightarrow D$ under
$D'\rightarrow D$.

{\em Case 1. } The curve $D'$ is disconnected, $D'=D_1'\amalg D_2'$. Then
$C'=C_1'\amalg C_2'$ with induced maps
$p_i':C_i'\rightarrow D_i'$, for $i=1,2$. Let $g_i=g(D_i')$ and
$\alpha_i=f\lst[D_i']$, for $i=1,2$. Then $g=g_1+g_2$
and $\alpha=\alpha_1+\alpha_2$. Now $g_i(C_i')\leq g_i(D_i')$ and
$f\lst[C_i']\leq f\lst[D_i']$ imply that
$g_i(C_i')=g_i$ and $f\lst[C_i']=\alpha_i$ and thus we may apply the induction
hypothesis to $p_1'$ and $p_2'$, proving
the surjectivity of $p$.

{\em Case 2. } The curve $D'$ is connected. Then $C'$ is connected, since
otherwise we would have two
curves contradicting the induction hypothesis. So me may apply the induction
hypothesis to $p':C'\rightarrow D'$.

This finally completes the proof of Theorem~\ref{mbfc}. \qed

\begin{defn} \label{doosf}
For a given object $(V,\tau)$ of $\VV\GG_s$, we let $\ol{M}(V,\tau)(T)$ be the
fiber of $\ol{M}(T)$ over
$(V,\tau)$ under the cofibration of Theorem~\ref{mbfc}.

Letting $T$ vary we get a stack $\ol{M}(V,\tau)$ on the category of
$k$-schemes with the fppf-topology.

For $(V,\tau)=(V,g,n,\beta)$ we denote $\ol{M}(V,\tau)$ by
$\ol{M}_{g,n}(V,\beta)$.
\end{defn}

If $\chr k\not=0$, let $L$ be a very ample invertible sheaf on $V$. Consider
only stable $V$-graphs for which
$\beta(v)(L)<\chr k$, for all $v\in V\t$. If this condition is satisfied, we
say that $\tau$ or $(V,\tau)$ is {\em
bounded by the characteristic}. If $\chr k=0$, we call {\em every }$(V,\tau)$
bounded by the characteristic, so that we
have uniform terminology.

\begin{them} \label{mgnvas}
For every $(V,\tau)$, bounded by the characteristic, the stack $\ol{M}(V,\tau)$
is a proper algebraic
Deligne-Mumford stack over $k$.
\end{them}
\begin{pf}
The proof will be postponed to a later section (see Corollary~\ref{mgnvass}).
\end{pf}

Every time we refer to $\ol{M}(V,\tau)$ as a Deligne-Mumford stack, we shall
tacitly assume that $(V,\tau)$ is bounded
by the characteristic.

\begin{rmk}
Theorems~\ref{mbfc} and~\ref{mgnvas} give rise to a functor
\begin{eqnarray*}
\ol{M}:\VV\GG_s & \longrightarrow & (\rtext{\normalshape proper algebraic
DM-stacks over $k$})\\
(V,\tau) & \longmapsto & \ol{M}(V,\tau),
\end{eqnarray*}
by choosing for every $k$-scheme $T$ a {\em clivage normalis\'e }(see
Definition~7.1 in \cite[Exp.\ ~VI]{sga1}) of the
cofibered category $\ol{M}(T)$ over $\VV\GG_s$. Of course, this functor
$\ol{M}$ is essentially independent of the
choice of the {\em clivage normalis\'e}.

Another way of stating this would be to construct a fibered category $\ol{M}$
over
$\VV\GG^{\op}_s\times(\text{$k$-schemes})$, such that $\ol{M}(V,\tau)(T)$ is
the fiber of $\ol{M}$ over
$(V,\tau,T)$ and $\ol{M}(T)$ is the fiber of $\ol{M}$ over $T$.
\end{rmk}

\section{Further Study of $\ol{M}$}

\begin{prop} \label{drfu}
Let $(V,\tau)$ be an object of $\VV\GG_s$, bounded by the characteristic of
$k$. Then the diagonal
\[\Delta:\ol{M}(V,\tau)\longrightarrow\ol{M}(V,\tau)\times\ol{M}(V,\tau)\]
is representable, finite and unramified.
\end{prop}
\begin{pf}
The assumption that $(V,\beta)$ is bounded by the characteristic implies that
all stable maps of class $\beta$ are
separable. So by Lemma~\ref{rsmsc} we may reduce the case of stable maps to the
case of stable curves, which is
well-known.
\end{pf}

\begin{lem} \label{rsmsc}
Let $(C,x,f)$ and $(D,y,h)$ be $n$-pointed stable maps to $V$ over the base
$T$, and $t\in T(K)$ a geometric point of
$T$. Assume that $f_t:C_t\rightarrow V$ and $h_t:D_t\rightarrow V$ are
separable morphisms. Then there exists an \'etale
neighborhood $S\rightarrow T$ of $t$, an integer $N$, markings
$x'=(x_1',\ldots,x_N')$ of $C_S$ and
$y'=(y_1',\ldots,y_N')$ of $D_S$ such that $(C_S,x_S,x')$ and $(D_S,y_S,y')$
are stable marked curves over $S$ and a
closed immersion of sheaves on $(\rtext{$S$-schemes})$
\[\Isomu((C,x,f),(D,y,h))_S\longrightarrow\Isomu((C_S,x_S,x'),(D_S,y_S,y')).\]
\end{lem}
\begin{pf}
Without loss of generality assume that $C$ and $D$ have the same genus $g$ and
$f$ and $h$ have the same class $\beta$.
Choose an embedding $\mu:V\hookrightarrow\pp^r$, let $d=\mu\lst\beta$ and
reduce to the case $V=\pp^r$ and $d=\beta$.
Let $N=d(r+1)$. Choose linearly independent hyperplanes $H_0,\ldots,H_r$ in
$\pp^r$ such that for each $i=0,\ldots,r$
\begin{enumerate}
\item no special point of $C_t$ or $D_t$ is mapped into $H_{i,K}$ under $f_t$
or $g_t$,
\item $f_t$ and $g_t$ are transversal to $H_{i,K}$.
\end{enumerate}
Then there exists an \'etale neighborhood $S\rightarrow T$ of $t$ such that
\begin{enumerate}
\item for each $i=0,\ldots,r$
\begin{enumerate}
\item $H_{i,S}\cap C_S$ gives rise to $d$ sections $x_{di+1}',\ldots,x_{di+d}'$
of $C_S$ over $S$,
\item $H_{i,S}\cap D_S$ gives rise to $d$ sections $y_{di+1}',\ldots,y_{di+d}'$
of $D_S$ over $S$,
\end{enumerate}
\item $(C_S,x_S,x')$ and $(D_S,y_S,y')$ are marked prestable curves.
\end{enumerate}
Then  $(C_S,x_S,x')$ and $(D_S,y_S,y')$ are in fact stable and there exists an
obvious morphism
\[\Isomu((C,x,f),(D,y,h))_S\longrightarrow\Isomu((C_S,x_S,x'),(D_S,y_S,y')),\]
which is clearly a closed immersion.
\end{pf}

\begin{lem} \label{lalss}
Let $(C,x_1,\ldots,x_{n+1},f)$ be a stable map and $(D,y_1,\ldots,y_n,h)$ the
stabilization under forgetting $x_{n+1}$.
Let $p:C\rightarrow D$ be the structure morphism. Then any section $y_0$ of $D$
making $(D,y_0,\ldots,y_n)$ a marked
prestable curve lifts uniquely to a section $x_0$ of $C$ making
$(C,x_0,\ldots,x_n)$ a marked prestable curve. If $y_0$
avoids $p(x_{n+1})$, then $(C,x_0,\ldots,x_{n+1})$ is a marked prestable curve.
\end{lem}
\begin{pf}
Let $V\subset D$ be the open subset consisting of smooth points of $D$ which
are not in the image of $y_i$, for any
$i=1,\ldots,n$. Let $U=p^{-1}(V)$. Then $p$ induces an isomorphism $p\resto
U:U\iso V$. Moreover, $U$ is smooth and
$x_{n+1}$ is the only section of $C$ which may meet $U$.
\end{pf}

\begin{prop} \label{uors}
Let $(C,x_1,\ldots,x_{n+1},f)$ and
$(\tilde{C},\tilde{x}_1,\ldots,\tilde{x}_{n+1},\tilde{f})$ be stable maps with
isomorphic stabilizations forgetting the $(n+1)$-st section. Let
$(C,y_1,\ldots,y_n,h)$ be such a stabilization, with
structure maps $p:C\rightarrow D$ and $\tilde{p}:\tilde{C}\rightarrow D$. If
$p(x_{n+1})=\tilde{p}(\tilde{x}_{n+1})$
then there exists a unique isomorphism $q:C\rightarrow\tilde{C}$ of stable maps
such that $\tilde{p}\comp q=p$.
\end{prop}
\begin{pf}
This is local over the base, so we may freely choose sections as necessary. In
fact, choose sections $z_1,\ldots,z_N$ of
$D$ in the smooth locus, avoiding $y_1,\ldots,y_n$ and
$\Delta=p(x_{n+1})=\tilde{p}(\tilde{x}_{n+1})$ and making
$$(D,z_1,\ldots,z_N,y_1,\ldots,y_n)$$ a stable marked curve. By
Lemma~\ref{lalss} these lift uniquely to sections
$w_1,\ldots,w_N$ of $C$ and $\tilde{w}_1,\ldots,\tilde{w}_N$ of $\tilde{C}$
making
$$(C,w_1,\ldots,w_N,x_1,\ldots,x_{n+1})$$ and
$$(\tilde{C},\tilde{w}_1,\ldots,\tilde{w}_N,\tilde{x}_1,
\ldots,\tilde{x}_{n+1})$$ marked prestable curves. Moreover,
these are clearly marked {\em stable } curves with a common stabilization
$$(D,z_1,\ldots,z_N,y_1,\ldots,y_n)$$ forgetting
the last section, such that $p(x_{n_1})=\tilde{p}(\tilde{x}_{n+1})$. Then they
have to be isomorphic by Knutson's theorem
(see~\cite{knudsen}) that $\ol{M}_{g,N+n+1}$ is the universal curve over
$\ol{M}_{N+n}$.
\end{pf}

\begin{prop} \label{flpuc}
Let $(C,x_1,\ldots,x_n,f)$ be a stable map and $\Delta$ a section of $C$. Then
there exists up to isomorphism a unique
stable map $(\tilde{C},\tilde{x}_1,\ldots,\tilde{x}_{n+1},\tilde{f})$ such that
$(C,x_1,\ldots,x_n,f)$ is the
stabilization of $(\tilde{C},\tilde{x}_1,\ldots,\tilde{x}_{n+1},\tilde{f})$
forgetting the $(n+1)$-st section and
$p(\tilde{x}_{n+1})=\Delta$, where $p:\tilde{C}\rightarrow C$ is the structure
map.
\end{prop}
\begin{pf}
Uniqueness follows from Proposition~\ref{uors}, hence existence is a local
question. Thus we may choose sections
$z_1,\ldots,z_N$ of $C$ such that $$(C,z_1,\ldots,z_N,x_1,\ldots,x_n)$$ is a
stable marked curve. By Knudsen's result
again, there exists a stable curve
$$(C',z_1',\ldots,z_N',x_1',\ldots,x_{n+1}')$$ whose stabilization forgetting
the last
section is $$(C,z_1,\ldots,z_N,x_1,\ldots,x_n)$$ and such that
$q(x_{n+1}')=\Delta$, where $q:C'\rightarrow C$ is the
structure map. Clearly, $$(C',z_1',\ldots,z_N',x_1',\ldots,x_{n+1}',f\comp q)$$
is a stable map. Then let
$(\tilde{C},\tilde{x}_1,\ldots,\tilde{x}_{n+1},\tilde{f})$ be the stabilization
of
$$(C',z_1',\ldots,z_N',x_1',\ldots,x_{n+1}',f\comp q)$$ forgetting the sections
$z_1',\ldots,z_N'$. By its universal
mapping property there exists a morphism $p:\tilde{C}\rightarrow C$ which makes
$(C,x_1,\ldots,x_n,f)$ the stabilization of
$(\tilde{C},\tilde{x}_1,\ldots,\tilde{x}_{n+1},\tilde{f})$ forgetting
$\tilde{x}_{n+1}$.
\end{pf}

\begin{cor} \label{esuc}
Let $C_{g,n}(V,\beta)$ be the universal curve over $\ol{M}_{g,n}(V,\beta)$.
Then the canonical morphism
$\ol{M}_{g,n+1}(V,\beta)\rightarrow C_{g,n}(V,\beta)$ induced by the $(n+1)$-st
section is an isomorphism. \qed
\end{cor}

\begin{prop}
Let $(C,x,f)$ be a stable $(V,g,n,\beta)$-map over $T$. Then the set of $t\in
T$ such that $(C,x)$ is a stable marked
curve is open in $T$.
\end{prop}
\begin{pf}
The set of such $t$ is the set of all $t\in T$ for which $(C,x)$ is isomorphic
to its stabilization. For any morphism of
schemes, the set of elements of its source at which it is an isomorphism is
always open. Finally, use properness of
prestable curves.
\end{pf}

By this proposition we may define
\[U_{g,n}(V,\beta)\subset\ol{M}_{g,n}(V,\beta)\]
to be the open substack of those stable maps, whose underlying marked curve is
stable. The canonical morphism
$U_{g,n}(V,\beta)\rightarrow \ol{M}_{g,n}$ has as fiber over the marked curve
$(C,x)$ the scheme of morphisms form $C$
to $V$ of class $\beta$. By results of Grothendieck in \cite{fgaIV} this is a
quasi-projective scheme. Hence
$U_{g,n}(V,\beta)$ is an algebraic $k$-stack of finite type. Now, for given $n$
there exists an $N>n$ such that
$U_{g,N}(V,\beta)\rightarrow\ol{M}_{g,n}(V,\beta)$  is surjective. Since this
morphism is flat by Corollary~\ref{esuc},
it is a flat epimorphism, hence a presentation of $\ol{M}_{g,n}(V,\beta)$.
Together with Proposition~\ref{drfu} this
implies that $\ol{M}_{g,n}(V,\beta)$ is a finite type separated algebraic
Deligne-Mumford stack over $k$. This is then
true for all objects of $\VV\GG_s$, bounded by the characteristic.

\begin{cor} \label{mgnvass}
Theorem~\ref{mgnvas} is true.
\end{cor}
\begin{pf}
It only remains to show properness. This is easily reduced to the case
$(V,\tau)=(\pp^r,g,n,d)$ and follows from
Proposition~3.3 of \cite{pandh}.
\end{pf}

\section{An Operadic Picture}

\begin{defn}
Let $(\tau,\alpha)$ be an $A$-graph. Let $R\t\subset F\t\times F\t$ be defined
by $(f,\ol{f})\in R\t$ if and only if one
of the conditions
\begin{enumerate}
\item $\ol{f}=j\t(f)$,
\item $\del f=\del\ol{f}$ and for $v=\del f=\del\ol{f}$ we have
$g(v)=\alpha(v)=0$
\end{enumerate}
is satisfied. Let $\sim$ be the equivalence relation on $F\t$ generated by
$R\t$ and let
\[P\t=F\t/\sim.\]
(In fact, $P_{(\tau,\alpha)}$ would be better notation, but we will stick with
the abuse of notation $P_{\tau}$.)
\end{defn}

\begin{prop} \label{cmmgpe}
Let $a:(B,\sigma)\rightarrow(A,\tau)$ be a combinatorial morphism of marked
graphs. Then $a_F:F\s\rightarrow F\t$
preserves equivalence.  \qed
\end{prop}

\begin{rmk}
In fact, Condition~(\ref{commor3}) of Definition~\ref{commor} may be replaced
by requiring $a_F$ to preserve
equivalence.
\end{rmk}

\begin{prop} \label{dafpea}
Let $\phi:\tau\rightarrow\sigma$ be a contraction of $A$-graphs. Then
$\phi^F:F\s\rightarrow F\t$ preserves equivalence.
\qed
\end{prop}

\begin{prop} \label{dpc}
If
\[\begin{array}{ccccc}
B & & \pi & \stackrel{\psi}{\longrightarrow} & \rho \\
\ldiagup{\xi} & \phantom{\longrightarrow} & \ldiag{b} &  & \rdiag{a} \\
A & & \sigma & \stackrel{\phi}{\longrightarrow} & \tau
\end{array}\]
is a stable pullback, then the induced diagram
\[\comdiaback{P_\pi}{\psi^F}{P_\rho}{b}{}{a}{P\s}{\phi^F}{P\t}\]
commutes. \qed
\end{prop}

By Propositions~\ref{cmmgpe},~\ref{dafpea} and~\ref{dpc}, we have a
contravariant functor
\[P:\GG_s\longrightarrow(\rtext{finite sets})\]
given by $P(A,\tau)=P\t$ on objects. Composing with the functor
$\VV\GG_s\rightarrow\GG_s$ we get a contravariant
functor
\begin{eqnarray*}
P:\VV\GG_s & \longrightarrow & (\rtext{finite sets}) \\
(V,\tau)   & \longmapsto     & P\t.
\end{eqnarray*}
There is an obvious functor
\begin{eqnarray*}
\VV\times(\rtext{finite sets}) & \longrightarrow & \VV \\
(V,P)                          & \longmapsto     & V^P,
\end{eqnarray*}
contravariant in the second argument, and composing with $P$ times the natural
functor $\VV\GG_s\rightarrow\VV$ gives
rise to a covariant functor
\begin{eqnarray*}
P:\VV\GG_s & \longrightarrow & \VV \\
(V,\tau)   & \longmapsto     & V^{P\t},
\end{eqnarray*}
still denoted $P$, by abuse of notation. We may consider $\VV$ as a subcategory
of the 2-category of proper algebraic
Deligne-Mumford stacks over $k$ and consider this as a functor
\[P:\VV\GG_s\longrightarrow(\rtext{proper algebraic DM-stacks over $k$}).\]

Now fix an object $(V,\tau)$ of $\VV\GG_s$. Let $(C,x,f)$ be a stable
$(V,\tau)$-map over $T$. Then $x$ and $f$ define a
morphism
\begin{eqnarray*}
f(x):T & \longrightarrow & V^{F\t} \\
t      & \longmapsto     & (f(x_i(t)))_{i\in F\t}.
\end{eqnarray*}
By Corollary~\ref{zc} this morphism $f(x)$ factors through $V^{P\t}\subset
V^{F\t}$, so we consider it as a morphism
\[f(x):T\longrightarrow V^{P\t}.\]
Thus we get a map $\ol{M}(V,\tau)(T)\rightarrow P(V,\tau)(T)$. Since it is
compatible with base change $S\rightarrow
T$, we have a morphism of $k$-stacks
\[\ev(V,\tau):\ol{M}(V,\tau)\longrightarrow P(V,\tau).\]

\begin{prop}
We have defined a natural transformation of functors from $\VV\GG_s$ to
$(\rtext{proper algebraic DM-stacks over $k$})$
\[\ev:\ol{M}\longrightarrow P,\]
called {\em evaluation}.
\end{prop}

In the general framework of $\Gamma$-operads, this allows us to consider
(appropriate subfunctors of) $\ol{M}$ and $P$
as a modular operad and a cyclic endomorphism operad, respectively. The
evaluation functor then induces a structure of
$\ol{M}$-algebra on $V$.

\vfill\eject
\part{Gromov-Witten Invariants}

\section{Isogenies}

\begin{defn}
Let $\tau$ be a stable $A$-graph.
\begin{enumerate}
\item The {\em class }of $\tau$ is
\[\beta(\tau)=\sum_{v\in V\t}\beta(v).\]
\item The {\em Euler characteristic } of $\tau$ is
\[\chi(\tau)=\chi(|\tau|)-\sum_{v\in v\t}g(v).\]
\item If $|\tau|$ is non-empty and connected the {\em genus }of $\tau$ is
\[g(\tau)=1-\chi(\tau).\]
\end{enumerate}
\end{defn}

\begin{defn}
Let $\tau$ be a stable $V$-graph, where $V$ is of pure dimension.
\begin{enumerate}
\item The {\em dimension } of $\tau$ is
\[\dim(V,\tau)=\chi(\tau)(\dim V-3)-\beta(\tau)(\omega_V)+\# S\t-\# E\t,\]
where $\omega_V$ is the canonical line bundle on $V$.
\item The {\em degree } of $\tau$ is
\begin{multline*}
\deg(V,\tau)=\\
\beta(\tau)(\omega_V)+(\dim V-3)(\chi(\tau^s)-\chi(\tau))+(\# S_{\tau^s}-\#
S\t)-(\# E_{\tau^s}-\# E\t),
\end{multline*}
where $\tau^s$ is the absolute stabilization of $\tau$.
\end{enumerate}
\end{defn}

Note that
\[\dim(V,\tau)-\dim(\tau^s)=\chi(\tau^s)\dim V-\deg(V,\tau).\]

\begin{defn}
The stable $A$-graph with one vertex of genus and class zero and three tails
(no edges) shall be called the
{\em $A$-tripod}, or simply a {\em tripod}.
\end{defn}

\begin{defn}
Let $a:\tau'\rightarrow\tau$ be a combinatorial morphism of stable $A$-graphs.
A {\em tail map } for $a$ is a map
$m:S_{\tau'}\rightarrow S\t$. Let $(a,\tau',\phi):\tau\rightarrow\sigma$ be a
morphism of stable $A$-graphs. A {\em tail
map } for $(a,\tau',\phi)$ is a tail map for $a$.
\end{defn}

Since for a contraction $\phi:\tau'\rightarrow\sigma$ the map
$\phi^F:F\s\rightarrow F_{\tau'}$ induces a bijection
$\phi^S:S\s\rightarrow S_{\tau'}$, there is an obvious way to compose morphisms
with tail map of stable $A$-graphs. By
abuse we will also call the composition $m\comp\phi^S$ the tail map of
$(a,\tau',\phi,m)$.

\begin{defn} \label{dsft}
Let $a:\tau'\rightarrow \tau$ be a combinatorial morphism of stable $A$-graphs
with tail map $m:S_{\tau'}\rightarrow
S\t$. We say that $(a,m)$ is of type {\em stably forgetting a tail}, or that
$\tau'$ is obtained form $\tau$ by {\em
stably forgetting a tail}, if there exists a tail $f$ of $\tau$ such that
\begin{enumerate}
\item $\tau'$ is the stabilization of $\tau''$, where $\tau''$ is obtained from
$\tau$ by forgetting the tail $f$,
\item for all $h\in S_{\tau'}$, such that $a(h)$ is a tail, we have
$m(h)=a(h)$,
\item $f\not\in m(S_{\tau'})$.
\end{enumerate}
\end{defn}

Note that the tail $f$ is uniquely determined by $(a,m)$ and $f$ determines
$(a,m)$ uniquely. We say that $(a,m)$ stably
forgets the tail $f$.

\begin{numrmk} \label{sft}
Every combinatorial morphism (with tail map) of type stably forgetting a tail
is of one of the following types (notation
of Definition~\ref{dsft}).

{\em Type I (Incomplete case). } No stabilization is needed, i.e.\
$\tau'=\tau''$.
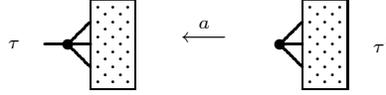
\[
\beginpicture
\setcoordinatesystem units <.3cm,.3cm> point at 3 2
\setplotarea x from 2 to 6, y from 0 to 4
\plot 3 2 4 3 /
\plot 2 2 4 2 /
\plot 3 2 4 1 /
\shaderectangleson
\setshadegrid span <1mm>
\putrectangle corners at 4 0 and 6 4
\put {\circle*{4}} [Bl] at 3 2
\axis left invisible label {$\scriptstyle\tau$} /
\axis right invisible label {\phantom{$\scriptstyle\tau$}} /
\endpicture  \stackrel{a}{\longleftarrow}
\beginpicture
\setcoordinatesystem units <.3cm,.3cm> point at 3 2
\setplotarea x from 3 to 6, y from 0 to 4
\plot 3 2 4 3 /
\plot 3 2 4 2 /
\plot 3 2 4 1 /
\shaderectangleson
\setshadegrid span <1mm>
\putrectangle corners at 4 0 and 6 4
\put {\circle*{4}} [Bl] at 3 2
\axis left invisible label {\phantom{$\scriptstyle\tau'$}} /
\axis right invisible label {{$\scriptstyle\tau'$}} /
\endpicture  \]

{\em Type II (Removing a tripod from a tail). }Only in this case does the tail
map $m$ contain any information. It
serves to `remember' which of the two tails in the following diagram is being
forgotten by $a$.
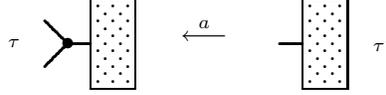
\[
\beginpicture
\setcoordinatesystem units <.3cm,.3cm> point at 3 2
\setplotarea x from 2 to 6, y from 0 to 4
\plot 2 3 3 2 /
\plot 3 2 4 2 /
\plot 2 1 3 2 /
\shaderectangleson
\setshadegrid span <1mm>
\putrectangle corners at 4 0 and 6 4
\put {\circle*{4}} [Bl] at 3 2
\axis left invisible label {$\scriptstyle\tau$} /
\axis right invisible label {\phantom{$\scriptstyle\tau$}} /
\endpicture  \stackrel{a}{\longleftarrow}
\beginpicture
\setcoordinatesystem units <.3cm,.3cm> point at 3 2
\setplotarea x from 3 to 6, y from 0 to 4
\plot 3 2 4 2 /
\shaderectangleson
\setshadegrid span <1mm>
\putrectangle corners at 4 0 and 6 4
\axis left invisible label {\phantom{$\scriptstyle\tau'$}} /
\axis right invisible label {{$\scriptstyle\tau'$}} /
\endpicture  \]

{\em Type III (Removing a tripod from an edge). }
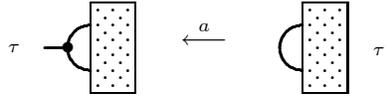
\[
\beginpicture
\setcoordinatesystem units <.3cm,.3cm> point at 3 2
\setplotarea x from 2 to 6, y from 0 to 4
\plot 2 2 3 2 /
\circulararc 180 degrees from 4 3 center at 4 2
\shaderectangleson
\setshadegrid span <1mm>
\putrectangle corners at 4 0 and 6 4
\put {\circle*{4}} [Bl] at 3 2
\axis left invisible label {$\scriptstyle\tau$} /
\axis right invisible label {\phantom{$\scriptstyle\tau$}} /
\endpicture  \stackrel{a}{\longleftarrow}
\beginpicture
\setcoordinatesystem units <.3cm,.3cm> point at 3 2
\setplotarea x from 3 to 6, y from 0 to 4
\circulararc 180 degrees from 4 3 center at 4 2
\shaderectangleson
\setshadegrid span <1mm>
\putrectangle corners at 4 0 and 6 4
\axis left invisible label {\phantom{$\scriptstyle\tau'$}} /
\axis right invisible label {{$\scriptstyle\tau'$}} /
\endpicture  \]

{\em Type IV (Forgetting a lonely tripod or a lonely elliptic component.) }
Only in this case does the number of
connected components of the geometric realization change.
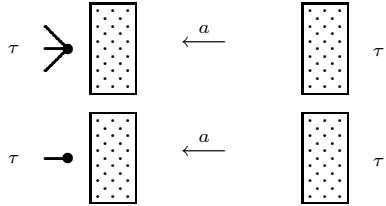
\[
\beginpicture
\setcoordinatesystem units <.3cm,.3cm> point at 3 2
\setplotarea x from 2 to 6, y from 0 to 4
\plot 2 3 3 2 /
\plot 2 2 3 2 /
\plot 2 1 3 2 /
\shaderectangleson
\setshadegrid span <1mm>
\putrectangle corners at 4 0 and 6 4
\put {\circle*{4}} [Bl] at 3 2
\axis left invisible label {$\scriptstyle\tau$} /
\axis right invisible label {\phantom{$\scriptstyle\tau$}} /
\endpicture  \stackrel{a}{\longleftarrow}
\beginpicture
\setcoordinatesystem units <.3cm,.3cm> point at 3 2
\setplotarea x from 3 to 6, y from 0 to 4
\shaderectangleson
\setshadegrid span <1mm>
\putrectangle corners at 4 0 and 6 4
\axis left invisible label {\phantom{$\scriptstyle\tau'$}} /
\axis right invisible label {{$\scriptstyle\tau'$}} /
\endpicture  \]
\[
\beginpicture
\setcoordinatesystem units <.3cm,.3cm> point at 3 2
\setplotarea x from 2 to 6, y from 0 to 4
\plot 2 2 3 2 /
\shaderectangleson
\setshadegrid span <1mm>
\putrectangle corners at 4 0 and 6 4
\put {\circle*{4}} [Bl] at 3 2
\axis left invisible label {$\scriptstyle\tau$} /
\axis right invisible label {\phantom{$\scriptstyle\tau$}} /
\endpicture  \stackrel{a}{\longleftarrow}
\beginpicture
\setcoordinatesystem units <.3cm,.3cm> point at 3 2
\setplotarea x from 3 to 6, y from 0 to 4
\shaderectangleson
\setshadegrid span <1mm>
\putrectangle corners at 4 0 and 6 4
\axis left invisible label {\phantom{$\scriptstyle\tau'$}} /
\axis right invisible label {{$\scriptstyle\tau'$}} /
\endpicture  \]
Here, the genus of the vertex displayed in the last diagram is equal to one.
\end{numrmk}

\begin{defn}
Let $(a,\tau',\phi,m):\tau\rightarrow\sigma$ be a morphism of stable $A$-graphs
with tail map. We call $(a,\tau',\phi,m)$
an {\em isogeny}, if
\begin{enumerate}
\item $(a,m)$ is a composition of morphisms of type stably forgetting a tail,
\item $\pi_0|\sigma|\rightarrow\pi_0|\tau|$ is bijective.
\end{enumerate}
We call the isogeny $\Phi:\tau\rightarrow\sigma$ an {\em elementary isogeny},
if it is an elementary contraction, or if
$\sigma$ is obtained from $\tau$ by stably forgetting a tail.
\end{defn}

\begin{note}
If $\Phi:\tau\rightarrow\sigma$ is an isogeny of stable $A$-graphs, then
$\chi(\sigma)=\chi(\tau)$.

An elementary isogeny either contracts a loop, or a non-looping edge or is of
type stably forgetting a tail~I, II,
or~III.

If we write $\Phi:\tau\rightarrow\sigma$ for an isogeny of stable $A$-graphs,
we denote the tail map of $\Phi$ by
$\Phi^S:S\s\rightarrow S\t$.
\end{note}

\begin{prop}
The composition of isogenies is an isogeny.
\end{prop}
\begin{pf}
Let
\[\begin{array}{ccccc}
A & & \pi & \stackrel{\psi}{\longrightarrow} & \rho \\
\ldiagup{\id} & \phantom{\longrightarrow} & \ldiag{b} &  & \rdiag{a} \\
A & & \sigma & \stackrel{\phi}{\longrightarrow} & \tau
\end{array}\]
be a stable pullback, where $a$ stably forgets the tail $f$ of $\tau$,
$\pi_0|\rho|\rightarrow\pi_0|\tau|$ is bijective
and $\phi$ is an elementary contraction of stable $A$-graphs. Then $b$ stably
forgets the tail $\phi^F(f)$ of $\sigma$.
Even if there is a vertex $v_0$ of $\tau$ which does not appear in $\rho$, this
vertex $v_0$ cannot be the vertex onto
which $\phi$ contracts an edge.
\end{pf}

Fix a semi-group with indecomposable zero $A$. We shall define a category
$\tilde{\GG}_s(A)$ from $\GG_s(A)$, retaining
only isogenies and morphisms of type cutting edges, but reversing the direction
of the latter, making them morphisms
{\em gluing tails}.

In fact, define the category $\tilde{\GG}_s(A)$ as follows. Objects of
$\tilde{\GG}_s(A)$ are stable $A$-graphs. A
morphism $\sigma\rightarrow\tau$ is a triple $(a,\sigma',\Phi)$, where
$a:\sigma\rightarrow\sigma'$ is a combinatorial
morphism of $A$-graphs of type cutting edges and $\Phi:\sigma'\rightarrow\tau$
is an isogeny of stable $A$-graphs. To
compose $(a,\sigma',\Phi):\sigma\rightarrow\tau$ and
$(b,\tau',\Psi):\tau\rightarrow\rho$, we need to construct a
diagram
\begin{equation}\label{dcmei}
\begin{array}{ccccc}
\sigma'' & \stackrel{\Xi}{\longrightarrow} & \tau' &
\stackrel{\Psi}{\longrightarrow} & \rho \\
\ldiagup{c} & & \rdiagup{b} &  & \\
\sigma' & \stackrel{\Phi}{\longrightarrow} & \tau & & \\
\ldiagup{a} & & & & \\
\sigma, & & & &
\end{array}
\end{equation}
where $c:\sigma'\rightarrow\sigma''$ is a combinatorial morphism of type
cutting edges and
$\Xi:\sigma''\rightarrow\tau'$ is an isogeny of stable $A$-graphs.

Let $f$ and $\ol{f}$ be two tails of $\tau$ such that $\{b(f),b(\ol{f})\}$ is
an edge of $\tau'$. Then construct
$\sigma''$ from $\sigma'$ by gluing the two tails $\Phi^S(f)$ and
$\Phi^S(\ol{f})$ to an edge. If $b$, cuts more than
one edge, iterate this process to construct $\sigma''$. This defines
composition of morphisms in $\tilde{\GG}_s(A)$,
which is clearly associative.

\begin{note}
In the situation of (\ref{dcmei}), we get a diagram in $\GG_s(A)$
\[\comdia{\sigma''}{\Xi}{\tau'}{\ol{c}}{}{\ol{b}}{\sigma'}{\Phi}{\tau},\]
which is easily seen to commute. Here, $\ol{b}$ and $\ol{c}$ are the morphisms
of stable $A$-graphs induced by $b$ and
$c$, respectively.
\end{note}

\begin{defn} \label{ecisg}
We call $\tilde{\GG}_s(A)$ the {\em extended category of isogenies of stable
$A$-graphs}, or the {\em extended isogeny
category }over $A$.

The morphisms in $\tilde{\GG}_s(A)$ are called {\em extended isogenies}. An
extended isogeny is called {\em elementary},
if it is an elementary isogeny or glues two tails to an edge.
\end{defn}

Now consider the following situation. Fix a smooth projective variety $V$ of
pure dimension. Let
$\Phi:\tau\rightarrow\sigma$ be an elementary extended isogeny of stable
modular graphs. Let $\sigma'$ be a stable
$V$-graph and $b:\sigma\rightarrow\sigma'$ a combinatorial morphism identifying
$\sigma$ as the absolute stabilization
of $\sigma'$. Note that $b$ is injective on vertices and complete, so that
$b:F\s(v)\rightarrow F_{\sigma'}(b(v))$ is
bijective, for all $v\in V\s$. Let $(a_i,\tau_i)_{i\in I}$ be a family of
pairs, where $I$ is a finite set and for each
$i\in I$ we have a combinatorial morphism $a_i:\tau\rightarrow\tau_i$
identifying $\tau$ as the absolute stabilization
of $\tau_i$. Finally, let for every $i\in I$ be given an extended isogeny of
stable $V$-graphs
$\Phi_i:\tau_i\rightarrow\sigma'$. In particular, for each $i\in I$ we have a
diagram of stable marked graphs (but note
that the horizontal and vertical morphisms live in different categories)
\[\comdia{\tau_i}{\Phi_i}{\sigma'}{\ol{a}_i}{}{\ol{b}}{\tau}{
\Phi}{\phantom{b}\sigma \phantom{b}.}\]
We shall now define what we mean by $(a_i,\tau_i,\Phi_i)_{i\in I}$ to be {\em
cartesian}, or a {\em pullback } of
$\sigma'$ under $\Phi$.  We have to distinguish six cases, according to which
kind of elementary extended isogeny $\Phi$
is.

Let us first consider the case that $\Phi$ is an elementary contraction
$\phi:\tau\rightarrow\sigma$, contracting the
edge $\{f,\ol{f}\}$ of $\tau$. As usual, let $v_1=\del f$, $v_2=\del\ol{f}$ and
$v_0=\phi(v_1)=\phi(v_2)$. Let
$w_0=b(v_0)$.

{\em Case I (Contracting a loop). } In this case $v_1=v_2$. The set $I$ has one
element, say $0$, and
$(a_0,\tau_0,\Phi_0)$ is cartesian, if $\Phi_0$ is a contraction contracting a
single loop $\{a_0(f),a_0(\ol{f})\}$ onto
$w_0$.

{\em Case II (Contracting a non-looping edge). } In this case $v_1\not=v_2$. We
require each $\Phi_i$ to contract
exactly one edge, namely $\{a_i(f),a_i(\ol{f})\}$ onto $w_0$. In particular,
this means that the only way the various
$(a_i,\tau_i,\Phi_i)$ differ is in the classes of $a_i(v_1)$ and $a_i(v_2)$. We
require that
$(\beta(a_i(v_1)),\beta(a_i(v_2)))_{i\in I}$ be a complete and non-repetitive
list of all pairs of elements of
$H_2(V)^+$ adding up to $\beta(w_0)$.

{\em Case III (Forgetting a tail). } Let us now deal with the case that
$\Phi:\tau\rightarrow\sigma$ stably forgets the
tail $f\in S\t$. We require $I$ to have one element, say $0$, and call
$(a_0,\tau_0,\Phi_0)$ cartesian if $\Phi_0$
stably forgets the tail $a_0(f)$ (and does nothing else).

{\em Case IV (Gluing two tails to an edge). } Finally, let us consider the case
that $\Phi$ is given by a combinatorial
morphism $c:\tau\rightarrow\sigma$, gluing the two tails $f$ and $\ol{f}$ of
$\tau$ to an edge $\{c(f),c(\ol{f})\}$ of
$\sigma$. Again, $I$ is required to have one element, say $0$, and
$(a_0,\tau_0,\Phi_0)$ is called cartesian if $\Phi_0$
glues two tails of $\tau_0$ to an edge of $\sigma'$ (and does nothing else).
Moreover, we require that $\Phi_0\comp
a_0=b\comp c$. An example:
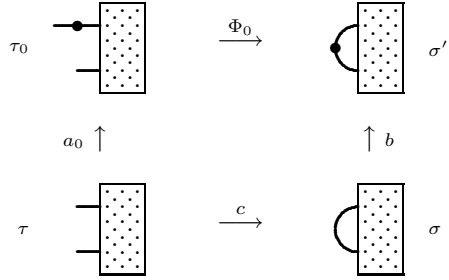
\[\begin{array}{ccc}
\beginpicture
\setcoordinatesystem units <.3cm,.3cm> point at 3 2
\setplotarea x from 2 to 6, y from 0 to 4
\plot 2 3 4 3 /
\plot 4 1 3 1 /
\shaderectangleson
\setshadegrid span <1mm>
\putrectangle corners at 4 0 and 6 4
\put {\circle*{4}} [Bl] at 3 3
\axis left invisible label {$\scriptstyle\tau_0$} /
\axis right invisible label {\phantom{$\scriptstyle\tau_0$}} /
\axis bottom invisible label {\phantom{.}} /
\endpicture  & \stackrel{\Phi_0}{\longrightarrow} &
\beginpicture
\setcoordinatesystem units <.3cm,.3cm> point at 3 2
\setplotarea x from 3 to 6, y from 0 to 4
\circulararc 180 degrees from 4 3 center at 4 2
\shaderectangleson
\setshadegrid span <1mm>
\putrectangle corners at 4 0 and 6 4
\put {\circle*{4}} [Bl] at 3 2
\axis left invisible label {\phantom{$\scriptstyle\sigma'$}} /
\axis right invisible label {{$\scriptstyle\sigma'$}} /
\axis bottom invisible label {\phantom{.}} /
\endpicture  \\
\ldiagup{a_0} & & \rdiagup{b} \\
\beginpicture
\setcoordinatesystem units <.3cm,.3cm> point at 3 2
\setplotarea x from 2 to 6, y from 0 to 4
\plot 4 1 3 1 /
\plot 3 3 4 3 /
\shaderectangleson
\setshadegrid span <1mm>
\putrectangle corners at 4 0 and 6 4
\axis left invisible label {{$\scriptstyle\tau$}} /
\axis right invisible label {\phantom{$\scriptstyle\tau$}} /
\axis top invisible label {\phantom{.}} /
\endpicture & \stackrel{c}{\longrightarrow} &
\beginpicture
\setcoordinatesystem units <.3cm,.3cm> point at 3 2
\setplotarea x from 3 to 6, y from 0 to 4
\circulararc 180 degrees from 4 3 center at 4 2
\shaderectangleson
\setshadegrid span <1mm>
\putrectangle corners at 4 0 and 6 4
\axis left invisible label {\phantom{$\scriptstyle\sigma$}} /
\axis right invisible label {{$\scriptstyle\sigma$}} /
\axis top invisible label {\phantom{.}} /
\endpicture
\end{array} \]

Note that in each case pullbacks exist, even though they are not necessarily
unique, even up to isomorphism, in the last
two cases. Note also, that for each $i\in I$ we have
$\deg(\tau_i)=\deg(\sigma')$.

We shall now define still another category, denoted $\tilde{\GG}_s(V)_{\cart}$,
called the {\em cartesian extended
isogeny category }over $V$.

\begin{defn} \label{ceicv}
Objects of $\tilde{\GG}_s(V)_\cart$ are pairs $(\tau,(a_i,\tau_i)_{i\in I})$,
where $\tau$ is a stable modular graph, $I$
is a finite set and for each $i\in I$ the pair $(a_i,\tau_i)$ is a stable
$V$-graph $\tau_i$, together with a
combinatorial morphism $a_i:\tau\rightarrow\tau_i$, identifying $\tau$ as the
absolute stabilization of $\tau_i$.

A {\em premorphism } from $(\tau,(a_i,\tau_i)_{i\in I})$ to
$(\sigma,(b_j,\sigma_j)_{j\in J})$ is a triple
$(\Phi,\lambda,(\Phi_i)_{i\in I})$, where $\Phi:\tau\rightarrow\sigma$ is an
extended isogeny of stable modular graphs,
$\lambda:I\rightarrow J$ is a map and for each $i\in I$ we have an extended
isogeny of stable $V$-graphs
$\Phi_i:\tau_i\rightarrow\sigma_{\lambda(i)}$. It is clear how to compose such
premorphisms and that composition is
associative.

A {\em morphism }in $\tilde{\GG}_s(V)_\cart$ is defined to be a premorphism
which can be factored into a composition of
elementary morphisms and isomorphisms. It is clear what an isomorphism is. An
{\em elementary morphism }is a premorphism
$(\Phi,\lambda,(\Phi_i)_{i\in I})$, as above, such that
\begin{enumerate}
\item $\Phi$ is an elementary extended isogeny,
\item for each $j\in J$ we have that
$(a_i,\tau_i,\Phi_i)_{i\in\lambda^{-1}(j)}$ is cartesian in the sense defined
in
Cases~I through~IV, above.
\end{enumerate}
\end{defn}

\begin{numrmk} \label{cfnc}
Projecting onto the first component defines a functor
\[\tilde{\GG}_s(V)_\cart\longrightarrow\tilde{\GG}_s(0).\]
Despite the notation, this is not a fibration of categories.
\end{numrmk}

We shall, in what follows, often shorten the notation $(\tau,(a_i,\tau_i)_{i\in
I})$ to $(\tau,(\tau_i)_{i\in I})$ or
even $(\tau_i)_{i\in I}$.

Call an object $(\tau_i)_{i\in I}$ of $\tilde{\GG}_s(V)_\cart$ {\em homogeneous
} of degree $n\in\zz$, if for all $i\in
I$ we have $\deg(V,\tau_i)=n$.

For a stable modular graph $\tau$, we may consider the fiber
$\tilde{\GG}_s(V)_{\cart/\tau}$ of the functor
$\tilde{\GG}_s(V)_\cart\rightarrow\tilde{\GG}_s(0)$ over $\tau$. In every such
fiber $\tilde{\GG}_s(V)_{\cart/\tau}$ we
have a functor
\[\oplus:\tilde{\GG}_s(V)_{\cart/\tau}\times\tilde{\GG}_s(V)_{\cart/\tau}
\longrightarrow\tilde{\GG}_s(V)_{\cart/\tau},\]
given by
\[(\tau_i)_{i\in I}\oplus(\sigma_j)_{j\in J}=((\tau_i)_{i\in
I},(\sigma_j)_{j\in J}),\]
where we think of the object on the right hand side as a family parametrized by
$I\amalg J$. The functor $\oplus$
satisfies some obvious properties, which we shall not list.

Every object $X=(\tau_i)_{i\in I}$ of $\tilde{\GG}_s(V)_\cart$ has a unique
decomposition $X=\bigoplus_{n\in\zz}X_n$
into homogeneous components. Every morphism in $\tilde{\GG}_s(V)_\cart$
respects this decomposition.

Finally, $\tilde{\GG}_s(V)_\cart$ is a tensor category (in the sense of
\cite{delmil}) with tensor product given by
\[\otimes:\tilde{\GG}_s(V)_{\cart}\times\tilde{\GG}_s(V)_{\cart}
\longrightarrow \tilde{\GG}_s(V)_{\cart},\]
which is defined by the formula
\[(\tau,(\tau_i)_{i\in I})\otimes(\sigma,(\sigma_j)_{j\in
J})=(\tau\times\sigma,(\tau_i\times\sigma_j)_{(i,j)\in I\times
J}).\]
For two graphs $\sigma$ and $\tau$ we denote by $\sigma\times\tau$ the graph
whose geometric realization is the
disjoint union of $|\sigma|$ and $|\tau|$. This notion extends in an obvious
way to marked graphs. The identity object
for $\otimes$ is the one element family with value the empty graph.

There are obvious compatibilities between these various structures on
$\tilde{\GG}_s(V)_\cart$. For example, if
$X=\bigoplus_nX_n$ and $Y=\bigoplus_mY_m$ are objects of
$\tilde{\GG}_s(V)_\cart$, then the decomposition of $X\otimes
Y$ into homogeneous components is given by
\[X\otimes Y=\bigoplus_r\left(\bigoplus_{n+m=r}X_n\otimes Y_m\right).\]
We summarize these properties by saying that $\tilde{\GG}_s(V)_\cart$ has
$\oplus$, $\otimes$ and $\deg$ structures.

A formally similar situation arises, for example, if we consider the category
of morphisms of an additive tensor
category $\CC$ in which all homomorphism groups are graded.  If we denote this
morphism category by $\MM\CC$, there is a
functor $\MM\CC\rightarrow\CC\times\CC$, given by source and target, whose
fibers have a graded $\oplus$-structure as
above. Also, $\MM\CC$ becomes a tensor category compatible with $\deg$ and
$\oplus$. So $\MM\CC$ has $\oplus$, $\otimes$
and $\deg$ structures. In fact, Gromov-Witten invariants may be thought of as a
functor from $\tilde{\GG}_s(V)_\cart$ to
$\MM\CC$ respecting the $\oplus$, $\otimes$ and $\deg$ structures. In this case
$\CC$ will be a category of motives.

\begin{defn} \label{doasvg}
A full subcategory $\tilde{\TT}_s(A)\subset\tilde{\GG}_s(A)$ is called {\em
admissible}, if it satisfies the following
axioms.
\begin{enumerate}
\item If $\Phi:\sigma\rightarrow\tau$ is an extended isogeny in
$\tilde{\GG}_s(A)$ and $\tau\in \ob\tilde{\TT}_s(A)$,
then $\sigma\in\ob\tilde{\TT}_s(A)$.
\item If $\sigma$ and $\tau$ are in $\tilde{\TT}_s(A)$, then so is
$\sigma\times\tau$.
\end{enumerate}
\end{defn}

For an admissible subcategory $\tilde{\TT}_s(A)\subset\tilde{\GG}_s(A)$ and a
homomorphism $\xi:A\rightarrow B$, the
full subcategory $\tilde{\TT}_s(B)\subset\tilde{\GG}_s(B)$ of graphs which are
stabilizations of objects of
$\tilde{\TT}_s(A)$ is admissible.

For a smooth projective variety $V$ of pure dimension, we may construct the
full subcategory
$\tilde{\TT}_s(V)_\cart\subset\tilde{\GG}_s(V)_\cart$, called the {\em
associated cartesian category}, which may be
characterized as the subcategory of $\tilde{\GG}_s(V)_\cart$ such that for each
object $(\tau,(a_i,\tau_i)_{i\in I})$ we
have that $\tau\in\ob\tilde{\TT}_s(0)$ and for all $i\in I$ that
$\tau_i\in\ob\tilde{\TT}_s(V)$. Note that
$\tilde{\TT}_s(V)_\cart$ inherits the $\oplus$, $\otimes$ and $\deg$ structures
from $\tilde{\GG}_s(V)_\cart$.

\begin{examples}
{\em I}. Call a marked graph $\tau$ a {\em forest}, if
\begin{enumerate}
\item $H^1(|\tau|)=0$,
\item $g(v)=0$, for all $v\in V\t$.
\end{enumerate}
Let $\tilde{\TT}_s(A)\subset\tilde{\GG}_s(A)$ be the full subcategory whose
objects are forests. Then $\tilde{\TT}_s(A)$
is an admissible subcategory, called the {\em tree level }subcategory of
$\tilde{\GG}_s(A)$.

{\em II}. Let $\tilde{\TT}_s(A)\subset\tilde{\GG}_s(A)$ be an admissible
subcategory. Let $d:A\rightarrow\zz_{\geq0}$ be
an additive map and $N>0$ an integer. Then let $\tilde{\TT}_s(A)_{d<N}$ be the
full subcategory of $\tilde{\TT}_s(A)$
given by the condition
\[\tau\in\ob\tilde{\TT}_s(A)_{d<N}\quad \Longleftrightarrow\quad
\mbox{$d(\beta(v))<N$,  for all $v\in V\t$}.\]
The subcategory $\tilde{\TT}_s(A)_{d<N}\subset\tilde{\GG}_s(A)$ is admissible.

If $A=H_2(V)^+$, then a very ample invertible sheaf $L$ on $V$ gives rise to
$d:H_2(V)^+\rightarrow\zz_{\geq0}$, by
setting $d(\beta)=\beta(L)$. If $\chr k\not=0$, we shall always pass to
$\tilde{\TT}_s(V)_{d<\chr k}$, in other words
assume that
\[\tilde{\TT}_s(V)=\tilde{\TT}_s(V)_{d<\chr k}.\]
But for emphasis, we may say that $\tilde{\TT}_s(V)$ is {\em bounded by the
characteristic}.
\end{examples}

\section{Orientations}

Fix a smooth projective variety $V$ of pure dimension.  Recall the following
five basic properties of $\ol{M}$.

{\em Property I (Mapping to a point). } Let $\tau$ be a stable $V$-graph of
class zero. Then $\tau$ is absolutely
stable. The evaluation morphism factors through $V^{\pi_0|\tau|}\subset
V^{P\t}$ and the canonical morphism
\[\ol{M}(V,\tau)\longrightarrow V^{\pi_0|\tau|}\times\ol{M}(\tau)\]
is an isomorphism. This follows immediately from Corollary~\ref{zc}. In
particular, $\ol{M}(V,\tau)$ is smooth.

Assume that $|\tau|$ is non-empty and connected. Let $(C,x)$ be the universal
family of stable marked curves over
$\ol{M}(\tau)$. Glue the $(C_v)_{v\in F\t}$ according to the edges of $\tau$ to
obtain a stable marked curve
$\pi:\tilde{C}\rightarrow\ol{M}(\tau)$ over $\ol{M}(\tau)$. Denote the vector
bundle of rank $g(\tau)\dim V$ on
$\ol{M}(V,\tau)$ given by $T_V\boxtimes R^1\pi\lst\O_{\tilde{C}}$ by
$\tT^{(1)}$.

{\em Property II (Products). } Let $\sigma$ and $\tau$ be stable $V$-graphs and
$\sigma\times\tau$ the obvious stable
$V$-graph whose geometric realization is the disjoint union of $|\sigma|$ and
$|\tau|$. There are obvious
combinatorial morphisms $\sigma\rightarrow\sigma\times\tau$ and
$\tau\rightarrow\sigma\times\tau$ giving rise to morphisms
of stable $V$-graphs $\sigma\times\tau\rightarrow\sigma$ and
$\sigma\times\tau\rightarrow\tau$ called the {\em
projections}. The induced morphism
\[\ol{M}(V,\sigma\times\tau)\longrightarrow\ol{M}(V,\sigma)
\times\ol{M}(V,\tau)\]
is an isomorphism. This follows directly from the definitions.

{\em Property III (Cutting edges). } Let $\Phi:\sigma\rightarrow\tau$ be a
morphism of stable $V$-graphs of type cutting
an edge. So $\Phi$ is induced by a combinatorial morphism
$a:\tau\rightarrow\sigma$. Let $f$ and $\ol{f}$ be the tails
of $\tau$ that come from the edge of $\sigma$ which is being cut by $\Phi$. So
this edge is $\{a(f),a(\ol{f})\}$.
The diagram of algebraic $k$-stacks
\begin{equation} \label{ceedd}
\comdia{\ol{M}(V,\sigma)}{\ev_{\{a(f),a(\ol{f})\}}}{V}{\ol{M}(\Phi)}{
}{\Delta}{\ol{M}(V, \tau)}{\ev_f\times
\ev_{\ol{f}}}{V\times V,}
\end{equation}
where the horizontal maps are evaluations at the indicated flags, is cartesian.
In particular, $\ol{M}(\Phi)$ is a
closed immersion. Again, this follows directly from the definitions.

{\em Property IV (Forgetting tails). } Let $\Phi:\sigma\rightarrow\tau$ be a
morphism of stable $V$-graphs stably
forgetting a tail. Denote the combinatorial morphism giving rise to $\Phi$ by
$a:\tau\rightarrow\sigma$ and the
forgotten tail by $f\in F\s$.

If $\Phi$ is of Type~I (i.e.\ incomplete), let $v=\del\s(f)$. Let
$\pi':C'\rightarrow\ol{M}(V,\sigma)$ be the universal
curve indexed by $v$ and $x:\ol{M}(V,\sigma)\rightarrow C'$ the universal
section given by $f$. Let
$\pi:C\rightarrow\ol{M}(V,\tau)$ be the universal curve indexed by the unique
vertex $w$ of $\tau$ such that
$a(w)=v$. Then by definition there is a commutative diagram
\[\comdia{C'}{}{C}{\pi'}{}{\pi}{\ol{M}(V,\sigma)}{\ol{M}(\Phi)}{
\ol{M}(V,\tau),}\]
and the section $x$ induces an $\ol{M}(V,\tau)$-morphism
\[\ol{M}(V,\sigma)\rightarrow C.\]
This is an isomorphism. In particular, $\ol{M}(\Phi)$ is proper and flat of
relative dimension one. This follows from
Corollary~\ref{esuc}.

If $\Phi:\sigma\rightarrow\tau$ removes a tripod, then
\[\ol{M}(\Phi):\ol{M}(V,\sigma)\rightarrow\ol{M}(V,\tau)\]
is an isomorphism. This is because
$\ol{M}(\text{$0$-tripod})=\ol{M}_{0,3}=\spec k$.

{\em Property V (Isogenies). } Let
\[(\Phi,\lambda,(\Phi_i)_{i\in I}):(\tau,(a_i,\tau_i)_{i\in
I})\longrightarrow(\sigma,(b_j,\sigma_j)_{j\in J})\]
be a morphism in $\tilde{\GG}_s(V)_\cart$, where $\Phi$ (and hence all
$\Phi_i$) is an isogeny, i.e.\ free of any tail
gluing factors.  For each $j\in J$ we have a commutative diagram
\[\comdia{\displaystyle\coprod_{i\in
I\atop\lambda(i)=j}\ol{M}(V,\tau_i)}{\amalg\ol{M}(\Phi_i)}{
\ol{M}(V,\sigma_j)}{\amalg
\ol{M}(\ol{a}_i)}{}{\ol{M}(\ol{b})}{\ol{M}(\tau)}{\ol{M}(\Phi)}{
\ol{M}(\sigma).}\]
This diagram should be considered close to being cartesian. See
Definition~\ref{domb} for a more precise statement. For
the moment let us note that the induced morphism
\[\coprod_{i\in
I\atop\lambda(i)=j}\ol{M}(V,\tau_i)\longrightarrow\ol{M}(\tau)
\times_{\ol{M}(\sigma)} \ol{M}(V,\sigma_j)\]
is surjective.

If $X$ is a separated algebraic Deligne-Mumford stack, by $A\lst(X)$ we shall
mean the rational Chow group of $X$ (see
\cite{vistoli}). If $X\rightarrow Y$ is a morphism of separated algebraic
Deligne-Mumford stacks, $A\upst(X\rightarrow
Y)$ will denote the rational bivariant intersection theory defined in
\cite{vistoli}.

\begin{defn} \label{domb}
Let $\tilde{\TT}_s(V)\subset\tilde{\GG}_s(V)$ be an admissible subcategory
(bounded by the characteristic). Let for each
$\tau\in\ob\tilde{\TT}_s(V)$ be given a cycle class
\[J(V,\tau)\in A_{\dim(V,\tau)}(\ol{M}(V,\tau)).\]
This collection of cycle classes is called an {\em orientation } of $\ol{M}$
over $\tilde{\TT}_s(V)$, if the following
axioms are satisfied.
\begin{enumerate}
\item \label{domb1} {\em (Mapping to a point). }We have
\[J(V,\tau)=c_{g(\tau)\dim V}(\tT^{(1)})\cdot[\ol{M}(V,\tau)],\]
for every stable $\tau\in\ob\tilde{\TT}_s(V)$ of class zero such that $|\tau|$
is non-empty and connected.
\item \label{domb2} {\em (Products). }In the situation of Property~II we have
\[J(V,\sigma\times\tau)=J(V,\sigma)\times J(V,\tau).\]
\item \label{domb3} {\em (Cutting edges). }In the situation of Property~III the
following is true. Let $[\ol{M}(\Phi)]\in
A^{\dim V}(\ol{M}(V,\sigma)\rightarrow\ol{M}(V,\tau))$ be the orientation class
of $\ol{M}(\Phi)$ obtained by pullback
(using Diagram~(\ref{ceedd})) from the canonical orientation $[\Delta]\in
A^{\dim V}(V\rightarrow V\times V)$. Then we
have
\[J(V,\sigma)=[\ol{M}(\Phi)]\cdot J(V,\tau).\]
In other words,
\[J(V,\sigma)=\Delta\upsh J(V,\tau),\]
where $\Delta\upsh$ is the Gysin homomorphism given by the complete
intersection morphism $\Delta$.
\item \label{domb4} {\em (Forgetting tails). }In the situation of Property~IV
the morphism $\ol{M}(\Phi)$ has a canonical
orientation $[\ol{M}(\Phi)]\in
A\upst(\ol{M}(V,\sigma)\rightarrow\ol{M}(V,\tau))$. We require that
\[J(V,\sigma)=[\ol{M}(\Phi)]\cdot J(V,\tau).\]
In other words,
\[J(V,\sigma)=\ol{M}(\Phi)\upst J(V,\tau),\]
where $\ol{M}(\Phi)\upst$ is given by flat pullback.
\item \label{domb5} {\em (Isogenies). }In the situation of Property~V, we have
for every $j\in J$ a class
\[\ol{M}(\Phi)\upsh J(V,\sigma_j)\in
A_{\dim(V,\sigma_j)}(\ol{M}(\tau)\times_{\ol{M}(\sigma)}\ol{M}(V,\sigma_j)),\]
since $\ol{M}(\Phi)$ has a canonical orientation, $\ol{M}(\tau)$ and
$\ol{M}(\sigma)$ being smooth of pure dimension. We
also have a morphism
\[h:\coprod_{\lambda(i)=j}\ol{M}(V,\tau_i)\longrightarrow\ol{M}(
\tau)\times_{\ol{M}(
\sigma)}\ol{M}(V,\sigma_j)),\]
which is proper. The requirement is that
\[h\lst(\sum_{\lambda(i)=j}J(V,\tau_i))=\ol{M}(\Phi)\upsh J(V,\tau).\]
\end{enumerate}
\end{defn}

\begin{numrmk} \label{afmsl}
To check Axiom~(\ref{domb5}), it suffices to do so for $\Phi$ an elementary
isogeny, $\#J=1$ and
$(a_i,\tau_i,\Phi_i)_{i\in I}$ a pullback. This follows from the projection
formula.
\end{numrmk}

\begin{example}
If $\tau$ is a stable $V$-graph such that $|\tau|$ is non-empty and connected,
define
\[J_0(V,\tau)=
\begin{cases}
c_{g(\tau)\dim V}(\tT^{(1)})\cdot[\ol{M}(V,\tau)] & \text{if $\beta(\tau)=0$,}
\\
0                                                 & \text{otherwise.}
\end{cases}\]
For an arbitrary stable $V$-graph $\tau$, let
$\tau=\tau_1\times\ldots\times\tau_n$, for stable $V$-graphs
$\tau_1,\ldots,\tau_n$, such that $|\tau|=|\tau_1|\amalg\ldots\amalg|\tau_n|$
is the decomposition of $|\tau|$ into
connected components. Then set
\[J_0(V,\tau)=J_0(V,\tau_1)\times\ldots\times J_0(V,\tau_n).\]
We claim that $J_0$ is an orientation of $\ol{M}$ over $\tilde{\GG}_s(V)$,
called the {\em trivial orientation}.
\end{example}

\begin{defn}
Call a smooth projective variety $V$ {\em convex}, if for every morphism
$f:\pp^1\rightarrow V$ (defined over an
extension $K$ of $k$) we have $H^1(\pp^1,f\upst T_V)=0$.
\end{defn}

\noprint{
\begin{prop}
Let $\pi:C\rightarrow T$ be a prestable curve of genus zero and $f:C\rightarrow
V$ a morphism to a convex variety. Then
$R^1\pi\lst f\upst T_V=0$.
\end{prop}
\begin{pf}
\end{pf}
}

\begin{prop}
Let $V$ be convex and $\tau$ a stable $V$-forest. Then $\ol{M}(V,\tau)$ is
smooth of dimension $\dim(V,\tau)$. Moreover,
the morphism
\[\ol{M}(V,\tau)\longrightarrow\ol{M}(\tau^s)\]
is flat of relative dimension $\chi(\tau^s)\dim V-\deg(V,\tau)$.
\end{prop}
\begin{pf}
Let us start with some general remarks. Let $\tau$ be an absolutely stable
$V$-graph. Then we define
\[U(V,\tau)\subset\ol{M}(V,\tau)\]
to be the open substack of those stable maps $(C,x,f)$, such that
$(C_v,(x_i)_{i\in F\t(v)})$ is  a stable marked curve,
for all $v\in V\t$.  Let $(C,x):T\rightarrow\ol{M}(\tau)$ be  a $T$-valued
point of $\ol{M}(\tau)$, i.e.\
$(C_v,(x_i)_{i\in F\t(v)})_{v\in V\t}$ is a family of stable marked curves
parametrized by $T$. Let
$(\tilde{C},\tilde{x})$ be the stable marked curve over $T$ obtained by gluing
the $C_v$ according to the edges of
$\tau$. The diagram
\[\comdia{\Mor_T(\tilde{C},V_T)}{}{T}{}{}{}{U(V,\tau)}{}{\ol{M}(\tau)}\]
is cartesian. In particular, by Grothendieck \cite{fgaIV}, the morphism
$U(V,\tau)\rightarrow\ol{M}(\tau)$ is
representable, separated and of finite type. Moreover, let $(C,x,f)$ be a
$K$-valued point of $U(V,\tau)$. Let
$(\tilde{C},\tilde{x})$ be the marked curve obtained by gluing the $C_v$ and
$\tilde{f}:\tilde{C}\rightarrow V$
the morphism induced by the $f_v$. If $H^1(\tilde{C},\tilde{f}\upst T_V)=0$,
then $(C,x,f)$ is a smooth point of
$U(V,\tau)\rightarrow\ol{M}(\tau)$ and we have
\[T_{U(V,\tau)/\ol{M}(\tau)}(C,x,f)=H^0(\tilde{C},\tilde{f}\upst T_V)\]
for the relative tangent space. (This is the case, if $\tau$ is a $V$-forest
and $V$ is convex.)

In this smooth case we may calculate the relative dimension of $U(V,\tau)$ over
$\ol{M}(\tau)$ at $(C,x,f)$ as
\begin{eqnarray*}
\dim_K H^0(\tilde{C},\tilde{f}\upst T_V)  & = & \chi(\tilde{f}\upst T_V) \\
& = & \deg\tilde{f}\upst T_V + \rk(\tilde{f}\upst T_V)\chi(\tilde{C}) \\
& = & -\beta(\tau)(\omega_V) + \dim V \chi(\tau) \\
& = & \dim(V,\tau)-\dim(\tau).
\end{eqnarray*}
Since $\ol{M}(\tau)$ is smooth of dimension $\dim(\tau)$, we get that
$U(V,\tau)$ is smooth of dimension $\dim(V,\tau)$
at $(C,x,f)$.

Now let $\tau$ be an arbitrary stable $V$-graph. Then there exists an
absolutely stable $V$-graph $\tau'$, together with
a morphism $\tau'\rightarrow\tau$ of type forgetting tails, such that the
morphism
\[U(V,\tau')\longrightarrow\ol{M}(V,\tau)\]
is surjective, hence a flat epimorphism of relative dimension
$\#S_{\tau'}-\#S\t$. So if $U(V,\tau')$ is smooth of
dimension $\dim(V,\tau')$, then $\ol{M}(V,\tau)$ is smooth of dimension
\[\dim(V,\tau')-\#S_{\tau'}+\#S\t=\dim(V,\tau).\]
Finally, by considering the commutative diagram
\[\comdia{U(V,\tau')}{}{\ol{M}(V,\tau)}{}{}{
}{\ol{M}(\tau')}{}{\ol{M}(\tau^s),}\]
we see that in this case $\ol{M}(V,\tau)\rightarrow\ol{M}(\tau^s)$ is flat of
relative dimension $\chi(\tau^s)\dim
V-\deg(V,\tau)$.
\end{pf}

\begin{them} \label{uosc}
Let $V$ be a convex variety and $\tilde{\TT}_s(V)\subset\tilde{\GG}_s(V)$ the
admissible subcategory of
$V$-forests bounded by the characteristic. Then the collection
\[J(V,\tau)=[\ol{M}(V,\tau)]\]
is an orientation of $\ol{M}$ over $\tilde{\TT}_s(V)$.
\end{them}
\begin{pf}
Let us check the axioms.

(1) {\em Mapping to a point. }  This follows from the fact that $g(\tau)=0$ and
hence
\[c_{g(\tau)\dim V}(\tT^{(1)})=c_0(0)=1.\]

(2) {\em Products. } In complete generality we have for smooth proper
Deligne-Mumford stacks $X$ and $Y$ that
\[[X\times Y]=[X]\times[Y]\]
in $A\lst(X\times Y)$.

(3) {\em Cutting edges. } Again we have a general fact to the following effect.
Consider the cartesian diagram of
separated Deligne-Mumford stacks
\[\comdia{X}{f}{V}{j}{}{i}{Y}{}{W,}\]
where $i$ and $j$ are regular embeddings such that for the normal bundles we
have
\[f\upst N_{V/W}=N_{X/Y}.\]
Then $i\upsh[Y]=[X]$. If all four participating stacks are smooth and $i$ and
$j$ are closed immersions of the same
codimension, then these conditions are automatically satisfied (see for example
Proposition~17.13.2 in \cite{ega4}).
Thus we may apply this fact in our case.

More generally, we have that $i\upsh[Y]=[X]$ if all participating stacks are
smooth and
\[\dim X+\dim W=\dim Y +\dim V.\]

(4) {\em Forgetting tails. } Again, there is a general fact that
$f\upsh[Y]=[X]$ if $f:X\rightarrow Y$ is a flat
morphism of smooth and proper Deligne-Mumford stacks.

(5) {\em Isogenies. } In accordance with Remark~\ref{afmsl} we assume that
$\Phi$ is an elementary isogeny, $\#J=1$ and
that $(a_i,\tau_i,\Phi_i)_{i\in I}$ is a pullback. There are five cases to
consider, according to what type of
elementary isogeny $\Phi$ is. We use notation as in the definition of pullback.

{\em Case I (Contracting a loop). } This case does not occur, since $\sigma$
and $\tau$ are forests.

{\em Case II (Contracting an edge). } We will start with some general remarks.
Let $\tau$ be a stable $V$-graph, and
$v_1,\ldots,v_n$ absolutely stable vertices of $\tau$, i.e.\ vertices $v$ such
that $2g(v)+|v|\geq3$. (To avoid
ill-defined notation we assume that $n\geq1$.) Let
\[U_{v_1,\ldots,v_n}(V,\tau)\subset\ol{M}(V,\tau)\]
be the open substack of all those stable maps $(C,x,f)\in\ol{M}(V,\tau)$ such
that $$(C_{v_{\nu}},(x_i)_{i\in
F\t(v_{\nu})})$$ is a stable marked curve, for all $\nu=1,\ldots,n$.

With this notation the diagram
\[\comdia{\displaystyle\coprod_{i\in I}
U_{a_i(v_1),a_i(v_2)}(V,\tau_i)}{}{U_{b(v_0)}(V,\sigma')}{}{
}{}{\ol{M}(\tau)}{}{\ol{M}(\sigma)}\]
is cartesian. Consider for a fixed $i\in I$ the open immersion
\[U_{a_i(v_1),a_i(v_2)}(V,\tau_i)\subset\ol{M}(V,\tau_i).\]
Let
\[Z_{a_i(v_1),a_i(v_2)}(V,\tau_i)\subset\ol{M}(V,\tau_i)\]
be the closed complement. We have
\[\dim Z_{a_i(v_1),a_i(v_2)}(V,\tau_i)<\dim \ol{M}(V,\tau_i).\]
Thus, to prove the equality of two cycles of degree $\dim(V,\tau_i)$ in
$\ol{M}(\tau)\times_{\ol{M}(\sigma)}\ol{M}(V,\sigma')$, it suffices to prove
the equality of the cycles restricted to
$\coprod_i U_{a_i(v_1),a_i(v_2)}(V,\tau_i)$. This reduces us to proving that
\[\ol{M}(\Phi)\upsh[U_{b(v_0)}(V,\sigma')]=\sum_i[U_{a_i(v_1),
a_i(v_2)}(V,\tau_i)].\]
This claim finally follows from the general fact already mentioned in the proof
of Axiom~(3).

{\em Case III (Forgetting a tail, incompletely). }
Let $f\in F\t$ be the forgotten flag, $v=\del_{\tau_0}(a_0(f))$ and $w\in
V_{\sigma'}$ the vertex of $\sigma'$
corresponding to $v$ via $\Phi_0$. We have an open immersion
\[U_v(V,\tau_0)\subset \ol{M}(V,\tau_0)\]
with closed complement
\[Z_v(V,\tau_0)\subset \ol{M}(V,\tau_0)\]
of strictly smaller dimension. Thus, as in the previous case, we may reduce to
proving that
\[\ol{M}(\Phi)\upsh[U_w(V,\sigma')]=[U_v(V,\tau_0)].\]
This follows from the fact that the diagram
\[\comdia{U_v(V,\tau_0)}{}{U_w(V,\sigma')}{}{}{
}{\ol{M}(\tau)}{\ol{M}(\Phi)}{\ol{M}(\sigma)}\]
is cartesian.

{\em Cases IV and V (Removing a tripod). } These cases are trivial, since
$\ol{M}(\Phi_0)$ and $\ol{M}(\Phi)$ are
isomorphisms.
\end{pf}

\noprint{
\begin{conjecture} \label{gmc}
Let $V$ be a Grassmannian variety. Then the collection
\[J(V,\tau)=[\ol{M}(V,\tau)]^0\]
is an orientation of $\ol{M}$ over $\tilde{\GG}_s(V)$.
\end{conjecture}
\begin{pf}
\end{pf}
}

\section{Deligne-Mumford-Chow Motives}

We shall imitate the usual construction of the category of Chow motives, as
described for example in \cite{scholl}.

Fix a ground field $k$. Let $\WW$ be the category of smooth and proper
algebraic Deligne-Mumford stacks over $k$. For an
object $X$ of $\WW$, let $A\upst(X)$ be the rational Chow ring of $X$ defined
by Vistoli \cite{vistoli}. Then $A\upst$
is a generalized cohomology theory with coefficient field $\qq$ in the sense of
\cite{kleiman}. Moreover, it is a graded
global intersection theory with Poincar\'e duality and cycle map in the
terminology of \cite{kleiman}.

If $X$ and $Y$ are objects of $\WW$ we define $S^d(Y,X)$, the group of {\em
correspondences }from $Y$ to $X$ of degree
$d$, to be
\[S^d(Y,X)=A^{n+d}(Y\times X),\]
if $Y$ is purely $n$-dimensional and
\[S^d(Y,X)=\bigoplus_i S^d(Y_i,X),\]
if $Y=\coprod_i Y_i$ is the decomposition of $Y$ into irreducible components.
Note that $S^d(Y,X)\subset A\upst(Y\times
X)$. The isomorphism $Y\times X\cong X\times Y$ exchanging components induces
an isomorphism
\[S^d(Y,X)\cong S^{d+n-m}(X,Y),\]
if $\dim Y=n$ and $\dim X=m$. We call this isomorphism {\em transpose }of
correspondences.
For objects $Z$, $Y$ and $X$ of $\WW$ we define composition of correspondences
by the usual formula
\[g\comp f={p_{13}}\lst(p_{12}\upst f\cdot p_{23}\upst g),\]
for $f\in S^d(Z,Y)$ and $g\in S^e(Y,X)$. Then $g\comp f\in S^{d+e}(Z,X)$.

The category $\ol{\WW}$ of {\em Deligne-Mumford-Chow motives } (or DMC-motives)
is now defined to be the category of
triples $(X,p,n)$, where $X\in \ob\WW$, $p\in S^0(X,X)$ such that $p^2=p$ and
$n\in\zz$. Homomorphisms are defined by
\[\Hom_{\ol{\WW}}((Y,q,m),(X,p,n))=p S^{n-m}(Y,X) q.\]
Note that $\Hom_{\ol{\WW}}((Y,q,m),(X,p,n))\subset S^{n-m}(Y,X)$. Composition
of homomorphisms in $\ol{\WW}$ is defined
as composition of correspondences.

There is a contravariant involution $\ol{\WW}\rightarrow\ol{\WW}$, denoted
$M\mapsto M^{\vee}$, defined by
$(X,p,n)^{\vee}=(X,^tp,\dim X-n)$, where $^tp$ is the transpose of $p$, on
objects and by transpose of correspondences on
homomorphisms.

\begin{prop}
The category $\ol{\WW}$ is a $\qq$-linear pseudo-abelian category. \qed
\end{prop}

Every morphism $f:X\rightarrow Y$ in $\WW$ defines a correspondence of degree
zero $\ol{f}\in S^0(Y,X)$ by
\[\ol{f}={\Gamma_f}\lst[X]\in A\upst(Y\times X),\]
where $\Gamma_f:X\rightarrow Y\times X$ is the graph of $f$.
We define the contravariant functor $h:\WW\rightarrow\ol{\WW}$ by
$h(X)=(X,\ol{\id}_X,0)$ and $h(f)=\ol{f}$. We usually
write $f\upst$ for $h(f)$ and $f\lst$ for $h(f)^{\vee}$.

Let $\ll=(\spec k,\ol{\id},-1)$ be the {\em Lefschetz motive}.  We shall use
the notation
\[M(n)=M\otimes\ll^{-n}.\]
We set
\[\Hom_{\ol{\WW}}^i(M,N)=\Hom_{\ol{\WW}}(M\otimes\ll^i,N)\]
and
\[\Hom_{\ol{\WW}}\upst(M,N)=\bigoplus_{i\in\zz}\Hom_{\ol{\WW}}^i(M,N).\]
The category with the same objects as $\ol{\WW}$, but with homomorphism groups
given by $\Hom_{\ol{\WW}}\upst(M,N)$ will be called the
category of {\em graded }DMC-motives.

For a DMC-motive $M$, define
\[A^i(M)=\Hom(\ll^i,M)\]
and
\[A\upst(M)=\bigoplus_iA^i(M).\]

\begin{prop}[Identity principle]
If $f,g:M\rightarrow N$ are two homomorphisms of DMC-motives, such that the
induced homomorphisms
\[A\upst(M\otimes h(X))\longrightarrow A\upst(N\otimes h(X))\]
agree, for all $X\in\ob\WW$, then $f=g$. \qed
\end{prop}

Let $\ol{\VV}$ be the category of Chow motives (which is defined as $\ol{\WW}$
is above, but starting with $\VV$ instead
of $\WW$). There is a natural fully faithful functor
$\ol{\VV}\rightarrow\ol{\WW}$.

\begin{question}
Is the functor $\ol{\VV}\rightarrow\ol{\WW}$ an equivalence of categories?
\end{question}

Let $H$ be a graded generalized cohomology theory on $\WW$ with a coefficient
field $\Lambda$ of characteristic zero,
possessing a cycle map such that $\pp^1$ satisfies epu (see \cite{kleiman}).
Then $H$ induces a covariant functor
(called a {\em realization functor\/})
\[\ol{H}:(\rtext{graded DMC-motives})\longrightarrow(\rtext{graded
$\Lambda$-algebras}),\]
such that for $X\in\ob\WW$ we have $\ol{H}(h(X))=H(X)$ and for a correspondence
$\xi\in S^d(Y,X)$ we have an induced
homomorphism
\begin{eqnarray*}
\ol{H}(\xi):H(Y) & \longrightarrow & H(X) \\
\alpha           & \longmapsto     &
{p_X}\lst({p_Y}\upst(\alpha)\cup\cl_{Y\times X}(\xi)).
\end{eqnarray*}
The functor $\ol{H}$ doubles the degree of a homomorphism.

The following are examples of such a cohomology theory $H$.
\begin{enumerate}
\item If $k=\cc$, consider to $X$ the associated topological stack $X^{\topo}$.
This is a stack on the category of
topological spaces with the \'etale topology. It has an associated \'etale
topos $X^{\topo}_{\et}$. Set
\[H_B(X)=H\upst(X^{\topo}_{\et},\qq)\]
and call it the {\em Betti cohomology }of $X$. Here $\Lambda=\qq$.
\item If $\ell\not=\chr k$ set
\[H_{\ell}(X)=H\upst(\ol{X}_{\et},\qql)={\displaystyle\projlim\limits_n}
H\upst(\ol{X}_{\et},\zz/\ell^n),\]
where $\ol{X}=X\times_{\spec k}\spec \ol{k}$ is the lift of $X$ to an algebraic
closure of $k$ and $\ol{X}_{\et}$
denotes the \'etale topos of $\ol{X}$. We call $H_{\ell}(X)$ the {\em
$\ell$-adic cohomology } of $X$. In this case
$\Lambda=\qql$.
\item If $\chr k=0$, let $\Omega_X\com$ be the algebraic deRham complex of $X$
and set
\[H_{dR}(X)=\hh\upst(X,\Omega_X\com).\]
We call $H_{dR}(X)$ the {\em algebraic deRham cohomology }of $X$. Here
$\Lambda=k$.
\end{enumerate}

\section{Motivic Gromov-Witten Classes}

Define the contravariant tensor functor
\[h(\ol{M}):\tilde{\GG}_s(0)\longrightarrow(\rtext{DMC-motives})\]
by $h(\ol{M})(\tau)=h(\ol{M}(\tau))$ on objects. For a morphism
$(a,\sigma',\Phi):\sigma\rightarrow\tau$ we have
$\ol{M}(\ol{a}):\ol{M}(\sigma')\rightarrow\ol{M}(\sigma)$ and
$\ol{M}(\Phi):\ol{M}(\sigma')\rightarrow\ol{M}(\tau)$.
Then let
\[h(\ol{M})(a,\sigma',\Phi)=\ol{M}(\ol{a})\lst\comp\ol{M}(\Phi)\upst.\]
This makes sense, because $\ol{M}(\ol{a})\lst$ is of degree zero,
$\ol{M}(\ol{a})$ being an isomorphism. This is also
why $h(\ol{M})$ is functorial.

Now fix a smooth projective variety $V$ of pure dimension and consider the
contravariant tensor functor
\[h(V)^{\otimes S}{(\chi\dim
V)}:\tilde{\GG}_s(0)\longrightarrow(\rtext{DMC-motives})\]
defined on objects by
\[\tau\longmapsto h(V)^{\otimes S\t}({\chi(\tau)\dim V}).\]
For a morphism $(a,\sigma',\Phi):\sigma\rightarrow\tau$ let $E$ be the set of
edges of $\sigma'$ which are cut by
$a:\sigma\rightarrow\sigma'$. Then we have
$V^{S\s}=V^{S_{\sigma'}}\times(V\times V)^E$. Let $p:V^{S_{\sigma'}}\times
V^E\rightarrow V^{S_{\sigma'}}$ be the projection,
$\Delta:V^{S_{\sigma'}}\times V^E\rightarrow V^{S_{\sigma'}}\times
(V\times V)^E=V^{S\s}$ the identity times the $E$-fold power of the diagonal.
Finally, we have an injection
$\Phi^S:S\t\rightarrow S_{\sigma'}$ giving rise to
$\Phi^S:V^{S_{\sigma'}}\rightarrow V^{S\t}$. We define the
homomorphism
\[h(V)^{\otimes S\t}{(\chi(\tau)\dim V)}\longrightarrow h(V)^{\otimes
S\s}{(\chi(\sigma)\dim V)}\]
as the composition of the three homomorphisms
\[(\Phi^S)\upst:h(V)^{\otimes S\t}{(\chi(\tau)\dim V)}\longrightarrow
h(V)^{\otimes
S_{\sigma'}}{(\chi(\sigma')\dim V)},\]
\[p\upst:h(V)^{\otimes S_{\sigma'}}{(\chi(\sigma')\dim V)}\longrightarrow
h(V)^{\otimes
S_{\sigma'}\cup E}{(\chi(\sigma')\dim V)}\]
and
\[\Delta\lst:h(V)^{\otimes S_{\sigma'}\cup E}{(\chi(\sigma')\dim
V)}\longrightarrow h(V)^{\otimes
S\s}{(\chi(\sigma)\dim V)},\] noting that $\chi(\tau)=\chi(\sigma')$ and
$\chi(\sigma')=\chi(\sigma)-\#
E$. Functoriality is a straightforward check using the identity principle.

Pulling back $h(\ol{M})$ and $h(V)^{\otimes S}{(\chi\dim V)}$ to the cartesian
extended isogeny category over $V$ via the
functor of Remark~\ref{cfnc}, we get two contravariant tensor functors
\[\tilde{\GG}_s(V)_\cart\longrightarrow(\text{graded DMC-motives}).\]

Now let $\tilde{\TT}_s(V)\subset\tilde{\GG}_s(V)$ be an admissible subcategory
(bounded by the characteristic) and $J$
an orientation of $\ol{M}$ over $\tilde{\TT}_s(V)$. For every object $\tau$ of
$\tilde{\TT}_s(V)$ we have a morphism
\[\phi_{(V,\tau)}:\ol{M}(V,\tau)\longrightarrow V^{S_{\tau^s}}\times
\ol{M}(\tau^s).\]
The first component is given by evaluation, noting that we have a map
$F_{\tau^s}\rightarrow F\t$.
Then
\begin{eqnarray*}
{\phi_{(V,\tau)}}\lst J(V,\tau) & \in &
S^{\dim(\tau^s)-\dim(V,\tau)}(V^{S_{\tau^s}},\ol{M}(\tau^s)) \\
 & = & \Hom_{\ol{\WW}}^{\deg(V,\tau)}(h(V^{S_{\tau^s}}){(\chi(\tau^s)\dim
V)},h(\ol{M}(\tau^s))).
\end{eqnarray*}

\begin{defn}
Define
\[I(V,\tau)={\phi_{(V,\tau)}}\lst J(V,\tau),\]
so that we have a homomorphism
\[I(V,\tau):h(V)^{\otimes S_{\tau^s}}{(\chi(\tau^s)\dim V)}\longrightarrow
h(\ol{M}(\tau^s)){(\deg(V,\tau))}\]
of DMC-motives over $k$. We call $I$ the system of {\em Gromov-Witten classes
}associated to the orientation $J$.
\end{defn}

Restricting the two functors $h(\ol{M})$ and $h(V)^{\otimes S}{(\chi\dim V)}$
to $\tilde{\TT}_s(V)_\cart$, we
get two contravariant tensor functors
\[\tilde{\TT}_s(V)_\cart\longrightarrow(\text{graded DMC-motives}).\]
We shall now define a natural transformation
\[I:h(V)^{\otimes S}{(\chi\dim V)}\longrightarrow h(\ol{M}).\]
So let $(\tau,(\tau_i)_{i\in I})$ be an object of $\tilde{\TT}_s(V)_\cart$, and
define
\[I(\tau,(\tau_i)_{i\in I})=\sum_{i\in I}I(V,\tau_i):h(V)^{\otimes
S_{\tau}}{(\chi(\tau)\dim V)}\longrightarrow
h(\ol{M}(\tau)).\]

\begin{them} \label{vagwc}
The Gromov-Witten transformation $I$ is a natural transformation compatible
with the $\oplus$, $\otimes$ and $\deg$
structures. Moreover,
\begin{enumerate}
\item \label{vag1} {\em (Mapping to a point). } The triangle
\[\comtri{h(V)^{\otimes S\t}{(\chi(\tau)\dim
V)}}{\rtext{mult}}{h(V){(\chi(\tau)\dim V)}}{I(V,\tau)}{c_{g(\tau)\dim
V}(\tT^{(1)})}{h(\ol{M}(\tau))}\]
commutes, for any stable $V$-graph $\tau$  of class zero in
$\tilde{\TT}_s(V)$, such that $|\tau|$ is non-empty and
connected.
\item \label{vag2} {\em (Divisor). }  Let $\lL\in\Pic(V)$ be a line bundle, so
its Chern class induces a homomorphism
$c_1(\lL):\ll\rightarrow h(V)$. Let $\Phi:\sigma\rightarrow\tau$ be a morphism
in $\tilde{\TT}_s(V)$ of type forgetting a tail,
such that the corresponding vertex of $\tau$ is absolutely stable. Then the
square
\[\begin{array}{ccc}
h(V)^{\otimes S_{\sigma^s}}{(\chi(\sigma^s)\dim V)} &
\stackrel{I(V,\sigma)}{\longrightarrow} &
h(\ol{M}(\sigma^s)){(\deg(V,\sigma))}\\
\ldiagup{c_1(\lL)} & & \rdiag{\ol{M}(\Phi)\lst} \\
h(V)^{\otimes S_{\tau^s}}{(\chi(\tau^s)\dim V)} \otimes\ll &
\stackrel{\beta(\lL)I(V,\tau)}{\longrightarrow} &
h(\ol{M}(\tau^s)){(\deg(V,\tau))}\otimes\ll
\end{array}\]
commutes.
\end{enumerate}
\end{them}
\begin{rmk}
To make this statement more precise, consider to $(\text{graded DMC-motives})$
the associated category of morphisms.
Then the natural transformation $I$ may be considered as a functor
\[I:\tilde{\TT}_s(V)_\cart\longrightarrow(\text{graded morphisms of
DMC-motives}).\]
Both categories have $\oplus$, $\otimes$ and $\deg$ structures and $I$
preserves them.  This essentially means that
\begin{enumerate}
\item $I((\tau_i)\oplus(\sigma_j))=I((\tau_i))+I((\sigma_j))$,
\item $\deg I((\tau_i))=\deg(\tau_i)$, if $(\tau_i)$ is homogeneous,
\item $I((\tau,\tau_i)\otimes(\sigma,\sigma_j))=I(\tau,\tau_i)\otimes
I(\sigma,\sigma_j)$.
\end{enumerate}
\end{rmk}
\begin{pf}
All this follows formally from Definition~\ref{domb} using the identity
principle and the bivariant formalism (as
explained for example in \cite{fulmacp}).
\end{pf}

\begin{rmks}
\begin{enumerate}
\item
Applying Theorem~\ref{uosc} we get the tree level system of Gromov-Witten
invariants for convex varieties.
\item
By applying a realization functor, we get Betti, $\ell$-adic or deRham
Gromov-Witten classes.
\item
Theorem~\ref{vagwc} implies all the axioms for Gromov-Witten classes listed in
\cite{KM}. Perhaps only Formula~(2.7) is
not quite evident. In view of its importance (it implies that the fundamental
class remains the identity with respect to
quantum multiplication), we will show that it follows from the rest of the
axioms.  In fact, assume that
\begin{equation} \label{leman}
\langle I_{0,3,\beta}\rangle (\gamma_1\otimes\gamma_2\otimes e^0) \ne 0.
\end{equation}

Choose a divisorial class $\delta$ with nonvanishing intersection with $\beta$.
In view of the Divisor Axiom, we must
then have

$$\langle I_{0,4,\beta}\rangle (\gamma_1\otimes\gamma_2\otimes\delta\otimes
e^0) \ne 0.$$

In view of (2.6), the last class is the lift of

$$\langle I_{0,3,\beta}\rangle (\gamma_1\otimes\gamma_2\otimes\delta ).$$

But this cannot be non-vanishing simultaneously with~(\ref{leman}) because the
Grading Axiom does not allow this.

More generally, this argument shows that whenever $e^0$ is among the arguments,
then $\langle I\rangle = 0$ for $\beta
\ne 0$, any genus, any $n.$ Geometrically: `if one of the points on $C$ is
unconstrained, the problem cannot have
finitely many (and non-zero) solutions.'

\noprint{The composition
\[h(V)^{\otimes2}{-\dim V}\stackrel{p\upst}{\longrightarrow}
h(V)^{\otimes3}{-\dim
V}\stackrel{I_{0,3}(V,\beta)}{\longrightarrow}\ll^{\beta(\omega_V)}\]
is zero, if $\beta\not=0$. Here $p:V^3\rightarrow V^2$ is a projection.}

\end{enumerate}
\end{rmks}

\end{document}